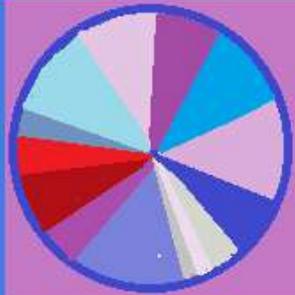
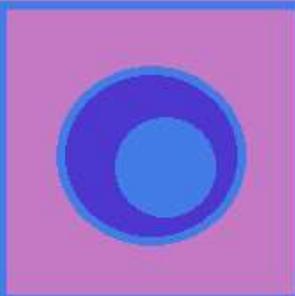
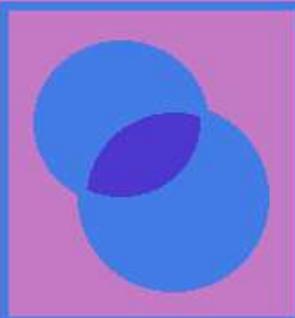
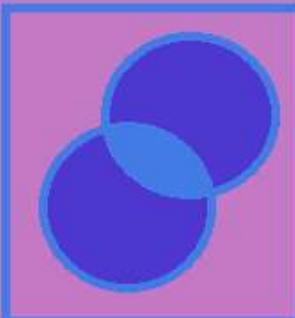
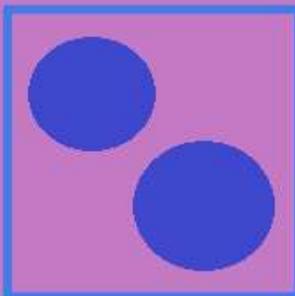
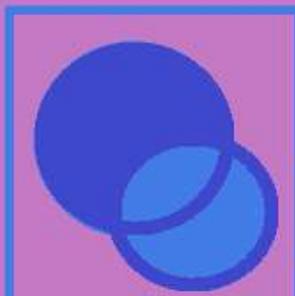

**ЧІНАРОВА Л.Л.**
**АНДРОНОВ І.Л.**

# ЕЛЕМЕНТИ ТЕОРІЇ ЙМОВІРНОСТЕЙ та МАТЕМАТИЧНОЇ СТАТИСТИКИ

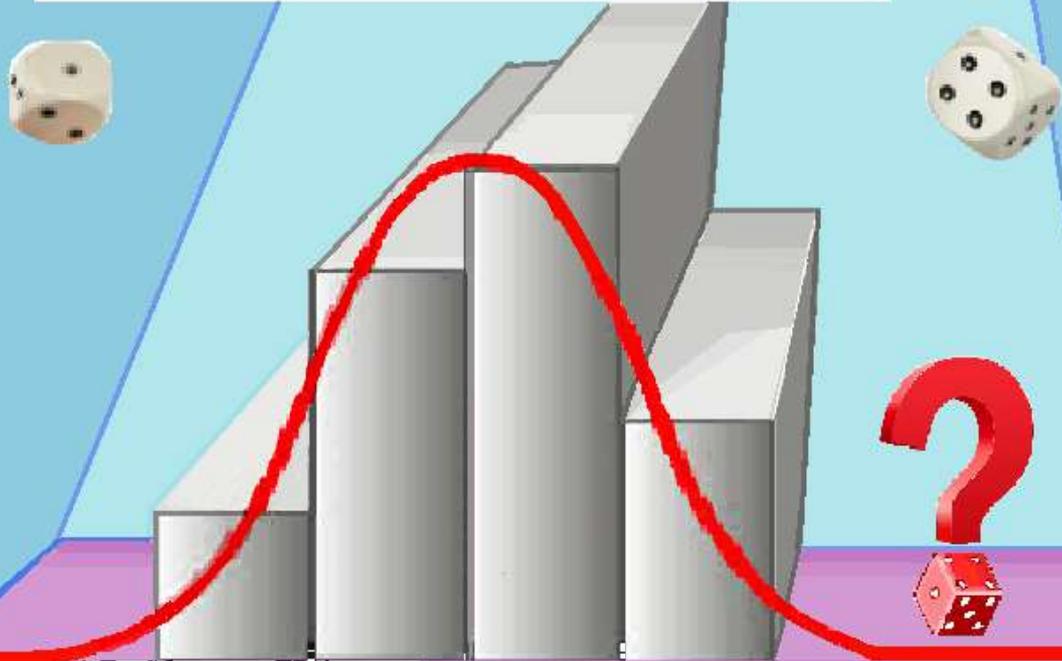

$P(AB)=P(A) \cdot P(B/A)=P(B) \cdot P(A/B)$  $P_n=n!$  $C=A-B$

$P(A)=\frac{k}{n}$  $P(A)+P(\bar{A})=1$  $DX=M(X-MX)^2$  $P(B/A)=\frac{P(AB)}{P(A)}$  $P(A+B)=P(A)+P(B)$  $C_n^m=\frac{n!}{m!(n-m)!}$

Л. Л. Чінарова
І. Л. Андронов

# ЕЛЕМЕНТИ ТЕОРІЇ ЙМОВІРНОСТЕЙ
## та
# МАТЕМАТИЧНОЇ СТАТИСТИКИ

**Рецензенти:**
**І. Б. Вавілова** – член-кор. НАНУ, доктор. фіз.-мат. наук, професор
**О. О. Панько** – доктор. фіз.-мат. наук, професор

**Чінарова Л. Л., Андронов І.Л.**
Ч93 Елементи теорії ймовірностей та математичної статистики. Навчальний посібник з курсу «Теорія ймовірностей та математична статистика» для спеціальностей 122 «Комп'ютерні науки», 124 «Системний аналіз», 125 «Комп'ютерна безпека» (для денної та заочної форм навчання) / Чінарова Л. Л., Андронов І.Л. – Львів: ННВК «АТБ», 2024. – 152 с.
ISBN 978-966-2042-81-8

Навчальний посібник розроблений на базі викладання освітньої компоненти «Теорія ймовірностей та математична статистика» (ТЙМС) викладачами кафедри «Математика, фізика та астрономія» Одеського національного морського університету. Приведений лекційний матеріал по основним аксіомам, теоремам та формулам статистичних розподілів та характеристик, який ілюстрований багатьма прикладами рішення конкретних задач. Обчислення можна проводити із використанням калькулятора або мов програмування чи електронних таблиць.

Навчальний посібник може бути використаний як для комп'ютерних спеціальностей 122, 124, 125, так і для інших технічних та економічних спеціальностей, а також, як додатковий навчальний матеріал, для гуманітаріїв.

# ЗМІСТ

# ВСТУП

Теорія ймовірностей та математична статистика (ТЙМС) є важливими компонентами освітніх програм усіх спеціальностей, які читають на загальних математичних кафедрах. Теорія ймовірностей допомагає моделювати та прогнозувати випадкові події та довірчі інтервали, що є критичним для аналізу та оптимізації алгоритмів у сфері комп'ютерних наук. Вона є продовженням курсу вищої математики, та знаходить подальше застосування в областях криптографії, штучного інтелекту, аналізу даних та інших важливих аспектів інформаційних технологій. Статистичні розподіли або ймовірнісні розподіли є математичними функціями, які визначають ймовірності різних значень, які може приймати випадкова величина. Вони є важливими компонентами теорії ймовірностей та математичної статистиці та знаходять застосування в різних галузях, включаючи фундаментальні, природничі та технічні науки, економіку, менеджмент та інші.

При повністю відомих законах еволюції (руху) та умовах, як початкових, так і майбутніх, система детермінована (класична фізика – точні рішення диференціальних рівнянь із точними законами та точними коефіцієнтами).

Випадковість «виникає» там, де є додаткові невраховані або змінні фактори, та, замість точних рішень, прогнози представляються з довірчими інтервалами. Довірчий інтервал - це інтервал значень, в якому з певною ймовірністю може знаходитися реальне значення величини. Довірчі інтервали є важливим інструментом у статистичному аналізі, оскільки вони враховують невизначеність та варіабельність (змінність) в зразках та допомагають зробити висновки про параметри популяції на основі доступних даних. Для визначення меж інтервалів, необхідне знання статистичного розподілу досліджуваної величини.

У теорії ймовірностей та статистиці є два основних типи випадкових величин – дискретні та неперервні. Різниця між ними полягає в тому, як вони приймають допустимі значення.

Дискретні випадкові величини приймають обмежену, окрему множину значень, які можуть бути перераховані одне за одним, максимально, від мінус нескінченності до плюс нескінченності.. Наприклад, це кількість «успіхів» (подій) при випробуваннях, кількість об'єктів даного класу у вибірці, послідовність вибраних об'єктів.

Неперервні випадкові величини: можуть приймати будь-яке значення в певному інтервалі, тому мають безліч можливих значень і не можуть бути перераховані окремо.

Наприклад, ріст людини, вага предмету, час, який знадобиться для виконання завдання - це приклади неперервних випадкових величин.

Різниця між ними може бути також у вигляді функцій ймовірностей. У випадку дискретних випадкових величин ймовірність, що випадає в конкретне значення, може бути визначена як маса ймовірностей. У випадку неперервних випадкових величин ймовірність попадання в конкретний діапазон визначається густині ймовірностей.

Для комп'ютерних та інших спеціальностей, особливо важливе вміння не лише застосовувати теоретичні знання, а й застосовувати їх на практиці, «доводячи до числа» на калькуляторі чи у програмних середовищах.

Навчальний посібник розроблено у стилі викладання математичних основ та їх безпосередньо застосування для рішення практичних завдань, зокрема, тих, які входять до домашнього завдання – типового розрахунку по першій частині курсу – теорії ймовірностей.

При використанні сучасних електронних таблиць (Microsoft Excel, Google DOCS, Libre/Open Office Calc, GNUmeric та ін.), приведені формули можуть бути використані для генерації нових випадкових послідовностей при кожному редагуванні, що надає «живий ефект анімації» порівняно із «фіксованими» методичними матеріалами по статистиці.

Для кожного статистичного розподілу, запропоновані додаткові завдання із генерації випадкових чисел по даному закону, порівняння із теоретичними значеннями. Такий підхід дозволяє наочно показати різноманітність задач, які розв'язує теорія ймовірностей, та дослідити статистичну значимість результатів.

# 1. ВИПАДКОВІ ПОДІЇ.

## 1.1. Відносна частота події та ймовірність події. Основні визначення понять.

Все, що відбувається або не відбувається в реальній дійсності, називаються подіями (іноді, явищами).

Теорія ймовірностей – це назва розділу математики, який займається дослідженням закономірностей у масових подій, тобто, вивчає математичні моделі випадкових подій.

Подію називають випадковою по відношенню до певного випробування (досліду), якщо в ході цього випробування вона може відбутися, а може й не відбутися.

Наприклад, якщо випробування полягає в одному киданні грального кубика, то в ході цього випробування можливі наступні події (результати випробування): на верхній грані кубика випаде число 1, число 2, ..., число 6. Рух та обертання кубика визначаються фізичними законами та відповідними диференціальними рівняннями. Якби усі початкові умови були ті самі при кожному кидку (положення та орієнтація кубика, однакова швидкість та траєкторія, швидкість вітру, густина повітря), то результат був би однаковий. Те, що результати можуть бути різними, залежить від різниці умов кидка. З іншої сторони, визначеність (= детермінованість) результату не означає того, що його можна заздалегідь передбачити, не знаючи усіх обставин (умов). Випадковість виникає внаслідок неповноти інформації.

Аналогічно, в кожного є номер телефону, паспорту, банківської картки, ідентифікаційний номер та багато інших чисел. Деяка частина з них може бути фіксована (код країни, провайдера, банку), а інша різна у різних користувачів. Це вже зафіксовано для конкретної особистості. Але інша людина може спробувати «вгадати». Вгадає, чи ні, теж може розглядатися, як подія, що відбулася, чи ні.

Контролер у транспорті дає квитки із номерами, що збільшуються. Можна передбачити, в кого буде квиток із наступним номером, якщо знати попередній. Але вгадати номер (принаймні, останні цифри) квитків тих, хто їх отримав, коли Ви ще не купили свій, це випадкова подія. Якщо 10 квитків поспіль придбали, то усі (останні) цифри будуть по одному разу. Але це не обов'язково буде, якщо взяти останні цифри випадковим чином у 10 пасажирів із більшій вибірки.

Випадкові події зазвичай позначаються початковими літерами латинського алфавіту: *A, B, C, D* та ін.

Будь-який результат випадкового експерименту називають випадковою подією. Внаслідок такого експерименту ця подія може чи відбутися, чи ні.

Під випадковим експериментом, розуміють різні експерименти, досліди, спостереження, вимірювання, випробування, результати яких залежать від випадку і які можливо повторювати багато разів в однакових умовах.

Подію U ($\Omega$) називають **достовірною** (вірогідною) по відношенню до певного випробування, якщо в ході цього випробування подія обов'язково відбудеться. $\Omega$ означає «генеральну сукупність» для даного випробування, множину усіх можливих наслідків.

Наприклад, достовірною подією буде поява одного з шести чисел (1,2,3,4,5,6) при одному киданні грального кубика. Тобто, якщо наслідків може бути лише 2, то звичайно поділяють результат на (A, $\bar{A}$), хоча можна писати із індексами, напр., $A_1 = A$, $A_2 = \bar{A}$. Для кубиків, можна взяти 6 індексів, в залежності від числа очок, для цифр – індекси 0, 1,…9 і т.д.

Подію $\emptyset$ (пуста множина) називають **неможливою** по відношенню до певного випробування, якщо в ході цього випробування подія не відбудеться.

Наприклад, неможливою подією є випадання числа 7 при киданні звичайного грального кубика. Чи отримати кількість подій більшу, ніж було випробувань.

Події називають **несумісними**, якщо появу однієї з них виключає поява інших подій в одному й тому ж випробуванні.

У результаті певного випробування обов'язково відбувається одна із взаємовиключних подій, причому кожна з них не поділяється на більш прості.

Такі події називаються елементарними подіями (або елементарними вихідними випробуваннями).

**Визначення:** «Сукупність подій $A_1, A_2, ..., A_n$ утворюють повну групу подій, якщо одна і тільки одна із цих подій в результаті експерименту обов'язково настає: $A_1 \cup A_2 \cup ... \cup A_n = \Omega, A_i \cap A_j = \emptyset$, якщо $i \neq j$. Події, що утворюють повну групу подій, називають

елементарними.»

Сукупність всіх можливих елементарних подій випробування називають **простором елементарних подій**. Нагадуємо позначення ∪ — об'єднання (кон'юнкція), ∩ — перетин (диз'юнкція). Ці позначення частіше використовують у теорії множин, в той час, як у теорії ймовірностей частіше використовують позначення арифметичних операцій, наприклад, $A + B$ та $A \cdot B$ (чи $AB$).

**1.2. Залежні та незалежні події.**

**Означення.** Подія $A$ називається незалежною (independent) від події $B$, якщо ймовірність появи події $A$ не залежить від того, відбулась подія $B$ чи не відбулась.

**Означення.** Подія $A$ називається залежною (dependent) від події $B$, якщо ймовірність події $A$ змінюється в залежності від того, відбулась подія $B$ чи не відбулась.

В якості прикладів розглянемо завдання:

• **Завдання 1**. Визначить, чи є задані події незалежними.

Нехай кинуто три монети. Подія $A$ — монети 1 і 2 впали однією й тією самою стороною, а 3 - іншою. Подія $B$ — монети 2 і 3 впали однією й тією самою стороною.

Відповідь: Ці події слід вважати залежними тому, що при здійсненні події $A$, подія $B$ вже не можлива і навпаки.

• **Завдання 2.** Визначить відповідь, у якій перераховано всі елементарні події, які можуть відбутися в результаті наступного випробування:

1) на підлогу кидається монета й визначається видима сторона.

Пропоновані варіанти відповіді:
- випадання однієї зі сторін - герба чи цифри;
- випадання герба;
- випадання цифри.

Відповідь: Щоб дати правильну відповідь згадаємо визначення - «Події, що утворюють повну групу подій, називають **елементарними**», тобто, це події які можуть відбуватися і зовсім не обов'язково відбудуться, таким чином відповідь - випадання герба чи цифри.

2) Визначить, чи є перераховані елементарні події рівно-можливими:

Відповідь: так.

Якщо при незмінних умовах випадковий експеримент проведено *n* разів і в *n(A)* випадків відбулася подія *A*, то число *n(A)* називається частотою події A.

Відносною частотою (або вибірковою ймовірністю) випадкової події називають відношення числа появ цієї події до загального числа проведених експериментів, тобто, відношення:

$$\frac{n(A)}{n}$$

### 1.3. Комбінації подій. Протилежні події.

Сумою (об'єднанням) подій $A$ і $B$ називається подія, яка полягає в тому, що відбувається хоча б одна з даних подій ($A$ або (or) $B$). Суму подій $A$ і $B$ позначають $A + B$ (або $A \cup B$). На рисунку за допомогою «кіл Ейлера» проілюстровано поняття суми подій $A$ і $B$. При цьому, якщо відбуваються обидві події, то зараховують лише один раз.

Це є історичною назвою, хоча це не обов'язково будуть «кола», а можуть бути еліпси, ромби, прямокутники чи фігури будь-якої форми.

Велике «коло» зображує всі елементарні події, які можуть відбутися в розглянутому випробуванні. Ліве «коло» зображує подію $A$, праве — подію $B$, а зафарбована область — подію $A + B$, (або $A \cup B$) (на рисунку зображена сума а) — для сумісних подій і б) — для несумісних подій).

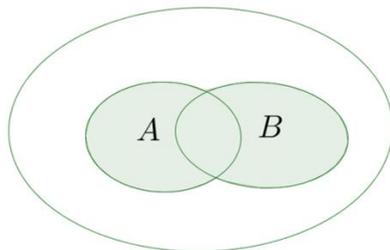

а)

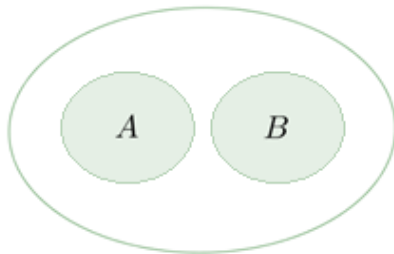

б)

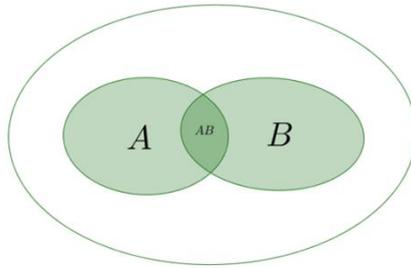

в)

Добутком (перетином) подій $A$ і $B$ називається подія, яка полягає в тому, що відбуваються обидві ці події – $A$ та (and) $B$. Добуток подій $A$ і $B$ позначають $AB$ (або $A \cap B$). Рисунок в) ілюструє за допомогою кіл Ейлера добуток подій $A$ і $B$. Область, зафарбована темніше (спільна частина кіл $A$ і $B$), ілюструє подію $AB$.

Події $A$ і $B$ називають рівними (рівносильними, тотожними, еквівалентними) й записують $A = B$, якщо подія $A$ відбувається тоді й тільки тоді, коли відбувається подія $B$. Наприклад, якщо у випробуванні з одним киданням грального кубика (подія $A$) випало число 6, а подія $B$ — випало найбільше з можливих чисел, то $A = B$.

Подію $\bar{A}$ називають **протилежною події $A$**, якщо подія $\bar{A}$ відбувається тоді й тільки тоді, коли не відбувається подія $A$. Читають, як «не $A$» (англійською "not $A$"). На рисунку проілюстровано взаємозв'язок подій $A$ і $\bar{A}$ на множині всіх елементарних наслідків випробування (подія $\bar{A}$ зображена як зафарбована область рис. г)).

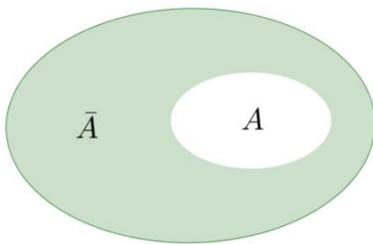    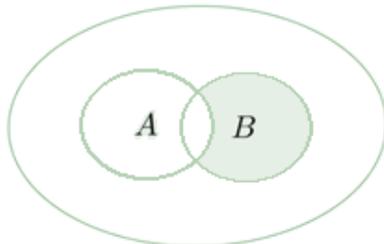

г)                                                                                           д)

Імовірність протилежної події: $P(\bar{A}) = 1 - P(A)$. Тобто, $P(A) + P(\bar{A}) = 1$, події $A, \bar{A}$ несумісні ($A\bar{A} = \emptyset$), отже, створюють «повну групу подій» ($A + \bar{A} = \Omega$) для цього експерименту.

Зауважимо, що, крім позначення $\bar{A}$, яке зручно набирати у рукописному конспекті чи на дошці, використовують «позначення у рядок», які можна набирати у тексті, чаті і т.і. : \A (часто); як у мові програмування С, !A ; іноді (дуже рідко) _A чи #A. Позначення «у рядок» використовують звичайно у коментарях до комп'ютерних програм чи електронних таблиць (Excel, Calc, GNUmeric та ін.), де не можна писати індекси та діакритичні знаки. Тобто, можна використовувати різні позначення тієї самої величини, відповідно до місця ($\bar{A}$ −підручник, конспект, рукописне рішення).

Різницею двох подій $A$ і $B$, тобто, різниця подій $B$ за винятком $A$ є така випадкова подія, для якої сприятливими є всі елементарні події, що є сприятливими для випадкової події $B$, за винятком тих, які є сприятливими для випадкової події $A$. Отже, маємо подію $B$ за винятком події $A$ (рис. д)):
$$B - A = B \backslash A = B\bar{A}$$

За принципом різниці множин, визначаються протилежні події. Так, $\bar{A} = \Omega \backslash A = \Omega\bar{A}$, тобто, сприятливими для події $\bar{A}$ є всі елементарні події простору $\Omega$ за винятком тих, що входять до складу події $A$.

Ще один із варіантів взаємозв'язку подій на прикладі кіл Ейлера (рис. е) ілюструє таке часткове об'єднання подій $A$ і $B$, при якому відбувається подія, яка полягає в тому, що здійснюється лише одна з даних подій при відсутності об'єднання цих подій (тобто, за видаленням події коли відбуваються обидві події $A$ та $B$). У програмуванні, такий результат позначається оператором **xor** (eXclusive OR), тобто $A \text{ xor } B = (A \neq B) = (A - B) + (B - A) = A\bar{B} + \bar{A}B$. Втім, **xor** звичайно використовується ширше - не до окремих подій, а для послідовностей 0 та 1 у представленні кодів символів у двійковій системі.

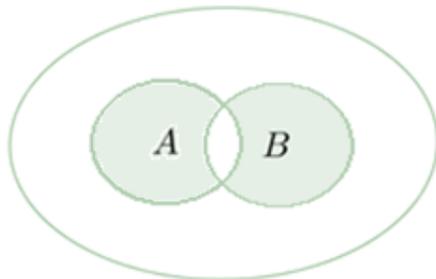

е)

На прикладі двох подій, приведемо таблицю різних комбінацій. Для зручності, можна пронумерувати події та ставити знак 1 (у мовах програмування, «умова виконується», значення логічного виразу є «істина=true», якщо відбувається, та 0 (хиба=false), якщо не відбувається Або «відбувається 0 (або 1) подій.

$\overline{A}\,\overline{B}$ (00) – обидві події не відбуваються (=відбувається 0 подій)

$\overline{A}\,B$ (01) – не відбувається перша, відбувається друга (одна подія)

$A\,\overline{B}$ (10) – відбувається перша, не відбувається друга (одна подія)

$A\,B$ (11) – відбувається як перша, так і друга (дві події)

$\overline{A}\,B + A\,\overline{B}$ (01 or 10) – відбувається рівно одна подія, або перша, або друга

$A = A\,B + A\,\overline{B} = A(B + \overline{B}) -$ відбувається перша подія незалежно від результату другої

$A + B = \overline{A}\,B + A\,\overline{B} + AB = (A + \overline{A})(B + \overline{B}) - \overline{A}\,\overline{B} = \overline{\overline{A}\,\overline{B}} -$ відбувається принаймні (хоча б) одна подія, тобто, відбувається подія, протилежна «відбувається 0 подій» (= не відбувається жодної події).

Для самостійного завдання, знайдіть для перших чотирьох рядків, де приведено «логічне множення» (та=and, кон'юнкція, перетин), значення операторів «або=or» (диз'юнкція, об'єднання), «виключна диз'юнкція = xor» (0, коли обидва аргументи рівні, 1, якщо не рівні), еквівалентність. Більше логічних функцій розглядається у курсі «дискретна математика».

В якості прикладів розглянемо завдання:

• **Приклад 1.** Двадцять карточок пронумеровані числами від 1 до 20. Із них

довільно вибирається одна карточка. Нехай подія $A$ — на карточці записано число, кратне 4; подія $B$ — на карточці записано число, кратне 6.

З'ясуйте, у чому полягає подія $A + B$.

Відповідь: У якості відповіді розглянемо два варіанта - вибрана карточка з числом 12 і вибрана карточка з одним із чисел: 4,6,8,12,16,18,20. Згадаємо визначення події A+B, «сумою (об'єднанням) подій $A$ і $B$ називається подія, яка полягає в тому, що відбувається хоча б одна з даних подій», тобто, будь-яке число кратне чи 4, чи 6. Цьому відповідає відповідь – «вибрана карточка з одним із чисел: 4,6,8,12,16,18,20».

• **Приклад 2**. Випробування складається з двох пострілів по мішені. Подія *A* — влучення по мішені при першому пострілі, подія *B* — при другому пострілі.

Визначте, у чому полягає подія *AB*.

Відповідь: Пропонуються наступні дві відповіді - влучення по мішені при обох пострілах і влучення по мішені хоча б при одному з двох пострілів. По визначенню «Добутком (перетином) подій *A* і *B* називається подія, яка полягає в тому, що відбуваються обидві ці події», згідно з визначенням відповіддю є влучення по мішені при обох пострілах.

• **Приклад 3.** Нехай *A* і *B* — довільні події. Необхідно обрати зі списку

введених позначень - $A\bar{B}, A\bar{B} + \bar{A}B, AB, A + B, \bar{A}\,\bar{B}$, наступну подію:

«відбулася принаймні одна з подій *A* і *B*»

Щоб вказати правильну відповідь проаналізуємо надані записи: перше позначення в цьому запису $A\bar{B}$, де – $\bar{B}$ це протилежна подія для події *B*, але у події, що розглядаємо ничого не говориться про наявність протилежної події, тому із списку можна виключити всі позначення де є позначення протилежної події і залишаться лише два позначення: *AB* – означає, що відбуваються обидві події та *A+B* - означає, що відбувається одна із подій.

Відповідь: $A + B$.

• **Приклад 4.** Зазначте подію, що є протилежною події: «хоча б на одному з двох кинутих гральних тетраедрів (=кубиків =кісток) з'явилося число 1». У якості відповіді пропонується розглянути два варіанта:
   - на обох гральних тетраедрах випали числа, відмінні від 1.
   - на обох гральних тетраедрах випали двійки.

Відповідь: «Подію $\bar{A}$ називають **протилежною події *A***, якщо подія $\bar{A}$ відбувається тоді й тільки тоді, коли не відбувається подія *A*.» Таким чином, згідно з визначенням, протилежною подією являється будь-яке число від двох до шістки на двох гральних тетраедрах, тобто, перша відповідь. Друга відповідь це одна із багатьох ймовірних комбінацій, яка відповідає умові «відмінні від 1».

• **Приклад 5.** Випробування полягає в киданні грального кубика. Подія $A_i$ (i=1,2,3,4,5,6) — випадання $i$ точок; подія $A$ — випадання парної кількості точок; подія $D$ — випадання кількості почок, більшої від 4. Виразіть події $A$ і $D$ через $A_i$ (i=1,2,3,4,5,6).

Відповідь: $A = A_2 + A_4 + A_6;\quad D = A_5 + A_6$

• **Приклад 6.** Дослід: на числову вісь кидається випадковим чином точка.

Події: $A$ — координата точки потрапляє в інтервал −2<*x*<8;

$B$ — точка потрапляє в інтервал 1<*x*<11.

Що становить собою подія *A*+*B*?

Відповідь: Згадаємо визначення події *A*+*B* – «Сумою (об'єднанням) подій

*A* і *B* називається подія, яка полягає в тому, що відбувається хоча б одна з даних подій», тобто, це об'єднання двох інтервалів, яке утворює інтервал $-2 < x < 11$.

• **Приклад 7**. Події *A* і *B* зображені за допомогою кіл Ейлера. Великим колом зображено всі елементарні результати випробування, з яким пов'язані події

*A* і *B*. Треба визначить подію, яка полягає в тому, що «відбулися обидві події *A* і *B*», тобто, вказати номер малюнку, що відповідає цієї події.

| 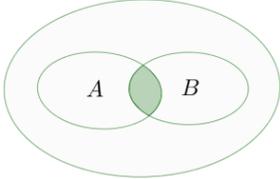 | 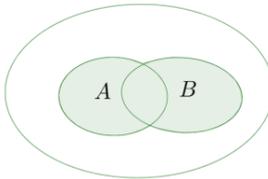 |
|---|---|
| *Мал.* 1 | *Мал.* 2 |
| 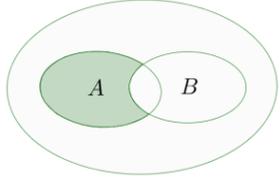 | 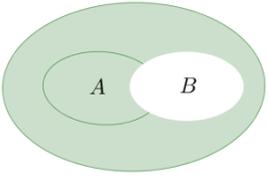 |
| *Мал.* 3 | *Мал.* 4 |

Відповідь: По визначенню «Добутком (перетином) подій *A* і *B* називається подія, яка полягає в тому, що відбуваються обидві ці події. Добуток подій *A* і *B* Позначають *AB* (або *A*∩*B*). На малюнку відповідає область, зафарбована темніше (спільна частина кіл *A* і *B*)», тобто, малюнок №1.

• **Приклад 8**. У коробці лежать червоні й сині фломастери. Випадковим чином вибирають два з них. Подія *A* полягає в тому, що обидва фломастери виявилися червоними. Подія *B* — один із них червоний, а один — синій.

Необхідно визначити, що означають події: $\bar{B}$, $A \cdot B$, $A \cdot \bar{B}$.

1. $\bar{B}$? Для такої задачі надаються наступні відповіді серед яких треба вибрати правильну:

○ *B*

○ *A*

○ ∅

⦿ взяті фломастери одного кольору: обидва червоні або обидва сині

○ хоча б один із вибраних фломастерів — червоний

○ хоча б один із двох фломастерів — синій

Відповідь: «взяті фломастери одного кольору: обидва червоні або обидва сині». Таку відповідь можна пояснити тим, що в даній задачі маємо тільки три ймовірні події – два фломастери різнокольорові чи вони мають один колір, але при цьому це можуть бути два красних, що відповідає події *A* чи два синіх. Згідно з визначенням «Подію $\bar{B}$ називають **протилежною події *B***, якщо подія $\bar{B}$ відбувається тоді й тільки тоді, коли не відбувається подія *B*», тобто, подія *B* — один із фломастерів червоний, а один — синій, не відбудеться якщо обидва фломастери будуть однокольорові і згідно з визначенням це і буде $\bar{B}$.

2. *A*·*B*? Пропоновані відповіді:

○ хоча б один із вибраних фломастерів — червоний

⦿ ∅

○ взяті фломастери одного кольору: обидва червоні або

обидва сині

   ○ хоча б один із двох фломастерів — синій

   ○ B

   ○ A

Відповідь: A·B=∅, якщо згадаємо, що AB це подія, яка полягає в тому, що відбуваються обидві події *A* і *B*, а вони майже протилежні події, тобто, якщо здійснюється *A*, не можлива *B* і навпаки.

3. $A \cdot \bar{B}$? Пропоновані відповіді:

   ○ взяті фломастери одного кольору: обидва червоні або обидва сині

   ○ хоча б один із двох фломастерів — синій

   ○ ∅

   ○ хоча б один із вибраних фломастерів — червоний

   ○ *B*

   ◉ *A*

Відповідь: $A \cdot \bar{B} = A$. Ми вже з'ясували, що $\bar{B}$ = «взяті фломастери одного кольору: обидва червоні або обидва сині», подія *A* по умовам задачі полягає в тому, що обидва фломастери виявилися червоними. а $A \cdot \bar{B}$ – це подія, коли відбуваються обидві події, тобто, їх перетин.

### 1.4. Ймовірність події. Поняття теорії ймовірностей.

Теорія ймовірностей — це розділ математики, що вивчає випадкові події та загальні властивості подій, послідовностей подій, процесів.

Випадковість є наслідком неповноти інформації. Наприклад, в кожного є конкретні цифри у телефонному номері – це вибірка, послідовність, яка вже відбулась, та комусь відома. Але, коли цифри називає людина, яка не знає цього номера, то може вгадати, або ні. Аналогічно, падіння кубика можна промоделювати, розв'язуючи диференціальні рівняння руху із врахуванням точних законів руху та початкових умов. Якщо умови будуть зберігатися, то кубик (чи монета) буде падати кожного разу тією ж самою стороною. В колоді карт (чи квитків на транспорт, бланків і т.д), послідовність детермінована. Але, якщо перемішати, чи розглядати квитки у пасажирів, які входять, переміщуються по транспорту, виходять, то важливою є випадковість.

Розділяють теоретичні характеристики (апріорні, до випробування), та вибіркові (наслідок конкретної реалізації). В останньому випадку, є кількість подій $m$ при $n$ випробуваннях. $m$ також називають частотою (чи абсолютною частотою), а відношення $m/n$ – відносною частотою, чи вибірковою ймовірністю.

Якщо при проведенні великої кількості випадкових експериментів, у кожному з яких може відбутися або не відбутися подія $A$, значення відносної частоти $m/n$ близькі до деякого певного числа, то це число називається імовірністю випадкової події $A$ і позначається $P(A)$, або Prob($A$), $P_A$, від англійського терміну Probability (ймовірність): $0 \leq P(A) \leq 1$.

Теоретична ймовірність визначається, як
$$P(A) = p = \lim_{n \to \infty} \frac{m}{n} .$$

Очевидно, що нескінчену кількість випробувань провести неможливо, можна лише провести обмежену кількість випробувань, але, із збільшенням $n$, точність статистично покращується пропорційно $\frac{1}{\sqrt{n}}$ згідно до «закону великих чисел».

У теорії ймовірностей, експерименти називаються дослідами, або випробуваннями, а можливі результати — наслідками. Усі можливі результати разом створюють множину наслідків. У випадку, коли якась подія відбувається, то часто кажуть «успіх» (математичний синонім довшого виразу «подія відбулася»), якщо не відбувається, то «неуспіх».

Для визначення теоретичної ймовірності, часто використовують моделювання замість випробувань. Наприклад, результат падіння кубика чи монети визначається із врахуванням симетрії (чи неважливих відхилень від неї). Математичне очікування кількості подій (=теоретична абсолютна частота) $m = pn$ не обов'язково має бути цілим. Наприклад, монета може випасти (вибірка) однією, чи іншою стороною (0 чи 1), але математичне очікування $m = np = 1 \cdot 0.5 = 0.5$. Очевидно, що випасти 0.5 герба чи 0.5 цифри випасти не може. Аналогічно, для непарного числа випробувань, математичне очікування буде пів-цілим (напр., $m = 14.5$ при $n = 29$. Навіть, при парному числі, вибіркова частота може співпадати із теоретичною, а може і не співпадати. Напр., при 2 випробуваннях, кількість подій може дорівнювати 0 ($P(0) = 1/4$), 1 $P(1) = 2/4$, або 2 ($P(2) =$ 1/4). Як буде розглянуто нижче, вибіркова частота $m$ сама є випадковим числом, яке має біноміальний розподіл із параметрами $n, p$. На рисунку, показані ймовірності кількості подій $m$ для різних значень кількості випробувань $n$ (для $p = 0.5$).

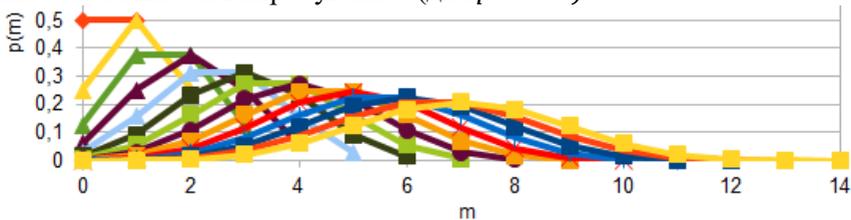

• **Приклад.** Монету кидають тричі й записують, якою стороною вгору вона падає — гербом (*г*) чи цифрою (*ц*). Знайди множину наслідків. Для студентів комп'ютерних спеціальностей, можна зробити позначення у двійковій системі числення (герб у 1,2,3 кидку оцінюється, як відповідно, $2^2=4$, $2^1=2$, $2^0=1$):

| 1-й кидок | 2-й кидок | 3-й кидок | Наслідки | Дв.система |
|---|---|---|---|---|
| цифра | цифра | цифра | *ццц* | 000 =0 |
| цифра | цифра | герб | *ццг* | 001 =1 |
| цифра | герб | цифра | *цгц* | 010 =2 |
| цифра | герб | герб | *цгг* | 011 =3 |
| герб | цифра | цифра | *гцц* | 100 =4 |
| герб | цифра | герб | *гцг* | 101 =5 |
| герб | герб | цифра | *ггц* | 110 =6 |
| герб | герб | герб | *ггг* | 111 =7 |

Множина результатів цього досліду складається з 8 однаково можливих результатів, із позначеннями у двійковій системі числення, це числа від 0 до 7.

У задачах, де 4 події, запис може бути від 0000 до 1111 (тобто, від 0 до 15, це $2^4=16$ наслідків.

Розглянемо поняття події в експерименті, тобто, подія це будь-яке твердження про результат досліду, правильність якого можливо перевірити. Наприклад, події — це *«випаде цифра»* і *«випаде герб»*, які можна позначити, як подія $A$ — *«випаде цифра»* та подія $B$ — *«випаде герб»*.

Подія, яка не може відбутися, називається неможливою. Наприклад, кинувши звичайний гральний кубик, випаде 14 пунктів.

**Визначення.** Неможлива подія (∅) — це подія, що у даному експерименті відбутися не може. P (∅) = 0. Якщо неможлива подія ∅, не відбудеться жодного разу при проведенні n експериментів, тоді її відносна частота буде дорівнювати: $\frac{n(\emptyset)}{n} = \frac{0}{n} = 0$.

Подія, яка відбувається завжди, називається достовірною чи вірогідною. Наприклад, кидаючи монету, випаде герб або цифра (інших можливостей немає). Якщо вірогідна подія U відбувається в кожному з n експериментів, то відносна частота її появи дорівнює:
$$\frac{n(U)}{n} = \frac{n}{n} = 1.$$

**Визначення:** Вірогідна подія — це подія *U*, яка обов'язково відбувається при кожному повторенні експерименту: $P(U) = 1$.

Рівно-можливі (рівно-імовірні) події — це такі події, кожна з яких не має ніяких переваг у з'явленні частіше інших при багаторазових експериментах, що проводяться в приблизно однакових умовах. Імовірності рівно-можливих подій однакові.

Події *A* і *B* називають несумісними, якщо вони не можуть відбутися одночасно в даному експерименті.

Якщо всі результати досліду однаково можливі, то ймовірність $P(A)$ будь-якої події A можна обчислити за формулою:
$$P(A) = \frac{\text{кількість наслідків, сприятливих події A}}{\text{кількість усіх можливих наслідків}}$$

**Визначення**: Ймовірністю випадкової події *A*, називають відношення числа рівно-можливих випадків, що сприяють події *A*, до числа усіх вірогідних випадків: $P(A) = \frac{m}{n}$ ,

де $n -$ загальне число рівноможливих випадків; $m -$ число випадків, що сприяють події A.

Ймовірність протилежної події можна обчислити за формулою: P($\overline{A}$)=1−P(A).

• **Приклад 1.** Кидається гральний кубик. Подія А — випаде цифра 2. Раніше вже було з'ясовано, що P(A)=$\frac{1}{6}$. Протилежна подія $\overline{A}$ не випаде цифра 2 (тобто, випаде 1,3,4,5 або 6).
$$P(\overline{A})=1-P(A)=1-\frac{1}{6}=\frac{5}{6}$$

Цю формулу особливо зручно використовувати, якщо в досліді багато наслідків,

• **Приклад 2.** У кошику лежать 100 пронумерованих куль. Яка ймовірність, що не витягнуть кульку під номером 6? Подія А — виймуть кульку під номером 6.

Подія $\overline{A}$ — вийнята кулька не буде під номером 6.
$$P(\overline{A})=1-P(A)=1-\frac{1}{100}=\frac{99}{100}$$

• **Приклад 3**.

У задачі 5000 варіантів. П'ять з них можна розв'язати набагато легше, ніж інші. Яка ймовірність того, що учневі дістанеться легкий варіант задачі?

Відповідь: Щоб відповісти на питання згадаємо визначення ймовірності – «Ймовірністю випадкової події А, називають відношення числа рівно-можливих випадків, що сприяють події А, до числа усіх вірогідних випадків». В даній задачі під чинниками, що сприяють події ми розуміємо п'ять легких варіантів, а загальна кількість випадків це є кількість усіх варіантів, тобто, маємо: $P(A)=\frac{5}{5000}=0{,}001$

• **Приклад 4.** У мисці 10 чорних кульок, 4 — білих та 5 кульок синього кольору. Навмання витягли одну кульку. Треба обчисли ймовірності:

а) Р (кулька білого кольору)=?

За умовою задачі маємо 4 білих кульки – це факт, що сприяє вірогідності витягнути саме білу кульку, а загальна кількість куль це сума кульок усіх кольорів, тобто, 10+4+5=19.

Відповідь: Р (кулька білого кольору) $=\frac{4}{19}$

Аналогічно обчислюємо наступні завдання:

b) Р (кулька чорного кольору) $=\frac{10}{19}$

c) Р (кулька синього кольору) $=\frac{5}{19}$

d) Р (кулька не синього кольору) $=\frac{14}{19}$, тобто, сума чорних та білих кульок, їх14.

е) Р (кулька чорного або синього кольору)$=\frac{10}{19}+\frac{5}{19}=\frac{15}{19}$, тобто, на відміну від попереднього випадку, маємо обчислити суму ймовірностей – ймовірність випадання чорної кульки плюс ймовірність випадання синьої кульки.

• **Приклад 5**. Сергій забув першу цифру п'ятизначного коду свого мобільного телефону.

1) Яка ймовірність того, що Сергій вгадає цю цифру?

Відповідь: Перша цифра може бути будь-якою цифрою від 0 до 9, тобто, маємо 10 варіантів першої цифри, додаткових умов чи факторів, що сприяють події – вгадати 1 цифру нема, тому $\frac{1}{10}$, тобто:

P=$\frac{1}{10} = 0{,}1$

2) Яка ймовірність того, що ця цифра виявиться парною?

Відповідь: Згадаємо визначення парного числа – «Парним називається число, яке можна поділити на 2 без залишку», тобто, це цифри 0, 2, 4, 6 та 8. Таким чином, всього парних цифр п'ять і це є сприятливий фактор, який замість 1 цифри надає 5 можливих, тому:

P=$\frac{5}{10} = 0{,}5$

3) Яка ймовірність того, що ця цифра не є нулем і ділиться на 8?

Відповідь: Спочатку розглянемо подію А того, що цифра є 0, тоді подія, яка полягає в тому, що цифра не є 0, являється протилежною події А і ми маємо використати формулу:

«Імовірність протилежної події:   P($\overline{A}$)=1− P(A)»

Знайдемо ймовірність події А:

P(A)= $\frac{1}{10}$ =0,1   тоді  P($\overline{A}$)=1−0,1=0,9.

Ймовірність P(B) − того, що цифра ділиться на 8 така ж, як і та, що цифра є нуль, тому що на вісім ділиться тільки цифра вісім, а це лише одна із десяти вірогідних цифр. При відповіді на питання «Яка ймовірність того, що ця цифра не є нулем і ділиться на 8?» ми якби шукаємо перетин цих двох подій. Подія, що цифра ділиться на вісім, тобто, що цифра є 8 − входить у подію, що цифра не є 0. Таким чином, у відповіді буде найменша ймовірність:

P(B)= $\frac{1}{10}$ =0,1

• **Приклад 6.** У комплекті гральних карт (колоді) 52 карти. Навмання витягується одна карта.

1) Скільки всього наслідків у даного експерименту?

Відповідь: Мається на увазі, що можемо витягнути будь-яку карту із 52, що надані, таким чином, наслідків буде 52.

2) Скільки сприятливих наслідків у події, якщо витягнута карта:

а) є картою пікової масті — Відповідь: 13. Обчислимо, скільки всього карт пікової масті у колоді. В колоді 4 різних мастей, тому $\frac{52}{4} = 13$. Кількість карт цієї масті і є кількістю сприятливих наслідків у події.

b) не є числом — Відповідь: 4. В колоді із 52 карт – карти з числами від 2 до 10, тобто, 9 шт. 13-9=4 це і є кількість сприятливих наслідків цієї події.

c) є валетом — Відповідь: 4. В колоді 4 масті і в кожній по одному валету.

d) не є дев'яткою — Відповідь: 48. Дев'яток у колоді 4 шт. по одній кожної масті, тому із загальної кількості віднімаємо кількість дев'яток і отримаємо кількість карт не дев'яток 52-4=48.

e) є валетом черви — Відповідь: 1 всього одна сприятлива подія, тому що картка єдина в колоді.

• **Приклад 7.** Обчисли, яка ймовірність того, що випадково назване двозначне число ділиться на 29.

Відповідь: Число 29 відноситься до простих чисел, тобто, це число яке ділиться на 1 та на само себе. Таким чином, серед двозначних чисел маємо лише одно, що підходить до умов задачі – само число 29. Усього двозначних чисел 90. Звідси маємо:

$$P = \frac{1}{90}$$

Якщо треба визначити протилежну ймовірність, тобто, ймовірність того, що випадково назване двозначне число не ділиться на 29, то маємо:

$$\bar{P} = 1 - \frac{1}{90} = \frac{89}{90}$$

• **Приклад 8.** Кинуто 2 кубики.
1. Скільки всього можливих наслідків у даного експерименту?

Відповідь: Стандартний кубик має шість граней, які пронумеровані від 1 до 6. Таким чином, ми маємо від кидання одного кубика 6 наслідків. З урахуванням ще шістьох наслідків другого кубика мали би вдвічі більше, але при киданні двох кубиків одночасно треба врахувати усі можливі комбінації випадання цифр, тобто, по шість комбінацій, при зміні лише одного кубика:

1-1; 1-2; 1-3; 1-4; 1-5; 1-6
2-1; 2-2…………..; 2-6

……………………………

6-1; 6-2……………..; 6-6.

Тобто, ми маємо шість строк по шість комбінацій, що відповідає 6х6=36 різним комбінаціям випадання кубика, тобто, маємо 36 наслідків.

2. Яка ймовірність наступних подій?

1) А — сума пунктів, що випали, дорівнює 7. Р(А)=?

Відповідь: Спочатку треба з'ясувати кількість комбінацій цифр, що складають суму 7. На гранях кубика маємо цифрі: 1, 2, 3, 4, 5, 6. В сумі сім дають наступні комбінації цифр: 1+6; 2+5; 3+4, крім того цифри, що випали на першому кубику та на другому можуть помінятися і отримаємо ще комбінації:

6+1; 5+2; 4+3. Таким чином маємо шість комбінацій цифр, що задовольняють умовам події А, тобто, маємо 6 наслідків, сприятливих події А. Кількість можливих наслідків ми з'ясували вище =36. Формула для підрахунку ймовірності події:

$$P(A) = \frac{\text{кількість наслідків, сприятливих події A}}{\text{кількість усіх можливих наслідків}} = \frac{6}{36} = \frac{1}{6}$$

2) В — сума пунктів, що, дорівнює 2. Р(В)=?

Відповідь: Сума, що дорівнює 2, можлива тільки при одній комбінації 1+1, тому Р(В)=$\frac{1}{36}$.

3) С — сума очок, що випали, більше, ніж 8. Р(С)=?

Відповідь: Більше 8 означає, що підходять тільки комбінації:

3+6

4+5; 4+6;

5+4, 5+5; 5+6;

6+3, 6+4, 6+5, 6+6, тобто, усього 1+2+3+4=10 комбінацій, що означає 10 наслідків, сприятливих події С. Таким чином, Р(С)=$\frac{10}{36} = \frac{5}{18}$

• **Приклад 9.** Кидаючи монетку 300 разів, 180 разів випала цифра.

1) Обчисли відносну частоту випадання цифри.

Відповідь: Згадаємо визначення: «Відносною частотою випадкової події називають відношення числа появ цієї події до загального числа проведених експериментів, тобто, відношення: $\frac{n(A)}{n}$», тобто, маємо $\frac{180}{300} = \frac{3}{5}$.

2) Знайди відносну частоту випадання герба.

Відповідь: Якщо усього монета кидалася 300 разів і 180 разів випадала цифра, то герб випадав решту – 120 разів. Тоді маємо $\frac{n(\text{В})}{n} = \frac{120}{300} = \frac{2}{5}$

• **Приклад 10.** Валерій, виїхавши в село, сім днів поспіль рахував і записував кількість машин, що проїжджали. У результаті вийшла таблиця:

| День тижня | Понеділок | Вівторок | Середа | Четвер | П'ятниця | Субота | Неділя |
|---|---|---|---|---|---|---|---|
| Легкові машини | 7 | 11 | 9 | 9 | 7 | 17 | 16 |
| Автобуси | 3 | 6 | 3 | 5 | 2 | 7 | 8 |

Яка ймовірність того, що в неділю першою проїжджала легкова машина?

Відповідь: Нас цікавить тільки неділя і за цей день всього було 16+8=24 машини. Підобчислюємо ймовірність події А – того, що проїжджала легкова машина P(A)=$\frac{16}{24} = \frac{2}{3}$. Ймовірність події В того, що проїжджали автобуси P(В)=$\frac{8}{24} = \frac{1}{3}$. Таким чином, ймовірність того, що першою проїжджала легкова машина у два рази більша $P = \frac{P(A)}{P(B)} = 2$. Але відповідь на поставлене питання це P(A)=$\frac{2}{3}$.

• **Приклад 11.** Приклад вибірки із $p = 0.75$ – зміни кількості подій (точки). Нахилені лінії показують максимальну ($p = 1$) та задану ($p = 0.75$) теоретичну кількість подій

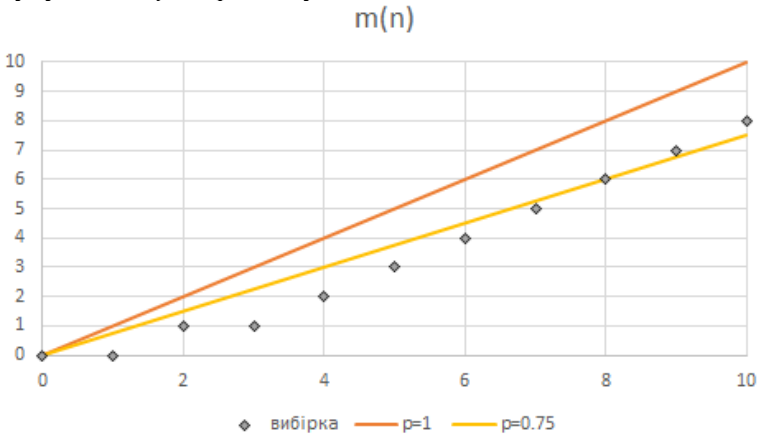

На наступних рисунках показані, для 10 різних вибірок, залежності від кількості випробувань кількості подій ($m(n)$, зліва) та відповідних відносних частот (вибіркових ймовірностей $p(n)$). Статистично, відхилення вибіркових значень від теоретичних пропорційні $\sqrt{n}$ для $m(n)$ та $1/\sqrt{n}$ для $p(n)$. Отже, для покращення точності $p(n)$ у $\beta$ разів, потрібно збільшити $n$ у $\beta^2$ разів.

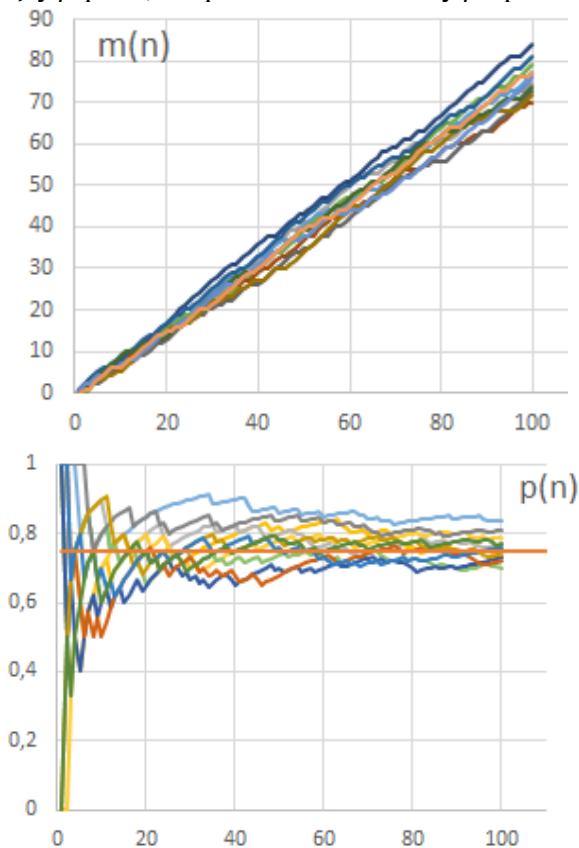

### 1.5. Несумісні події. Додавання ймовірностей.

Дві випадкові події A і B називаються несумісними, якщо їх добуток є неможливою подією, тобто, A·B = ∅ ( за інших позначень A∩B= ∅).

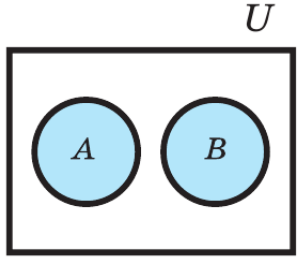

$A \cdot B = \varnothing$

Наприклад, при підкиданні грального кубика розглядають події: *A* — випало парне число очок, *B* — випало 1 очко, *C* — випало число очок, кратне 3. Тоді події *A* і *B* та події *B* і *C* — несумісні (не можуть відбутися одночасно). Події *A* і *C* — сумісні (можуть відбутися одночасно, якщо випаде 6 очок, тобто, A·C $\neq \varnothing$).

Ймовірність суми двох несумісних подій дорівнює сумі ймовірностей цих подій, тобто: $P(A + B) = P(A) + P(B)$.

У загальному випадку сумісних подій, $P(A + B) = P(A) + P(B) - P(AB)$.

Події є неспільними або несумісними, якщо поява однієї з них виключає появу іншої.

**Наслідок.** Сума ймовірностей протилежних подій дорівнює одиниці, тобто: $P(A) + P(\bar{A}) = 1$.

• **Приклад:** У ящику лежать 9 куль, із яких 2 — білі, 3 — червоні та 4 — зелені Навмання береться одна куля. Яка ймовірність того, що ця куля є кольоровою (не білою)? Розглянемо два способи вирішення цієї задачі.

<u>1 спосіб</u>. Нехай подія *A* — поява червоної кулі, подія *B* — поява зеленої кулі, тоді подія *A* + *B* — поява кольорової кулі. Очевидно, що: $P(A)=\frac{3}{9}=\frac{1}{3}$ та $P(B)=\frac{4}{9}$. Оскільки події *A* і *B* несумісні, до них можна застосувати теорему додавання ймовірностей: $P(A+B)=P(A)+P(B)=\frac{1}{3}+\frac{4}{9}=\frac{7}{9}$.

<u>2 спосіб</u>. Нехай подія C — поява білої кулі, тоді протилежна їй подія $\bar{C}$ — поява не білої (кольорової) кулі. Очевидно, що $P(C)=\frac{2}{9}$, а згідно з наслідком із теореми, маємо: $P(\bar{C})=1-P(C)=1-\frac{2}{9}=\frac{7}{9}$,

**Зауваження**

1. Теорема, аналогічна першій теоремі, правильна для будь-якої конкретної кількості подій, тобто, P($A_1$+$A_2$+...+$A_n$)=P($A_1$)+P($A_2$)+...+P($A_n$), де $A_1$;$A_2$;...;$A_n$ — попарно несумісні події.

2. Якщо $A_1$;$A_2$;...;$A_n$ — усі елементарні події деякого випробування, то їхня сукупність називається **полем подій**. Очевидно, що ці події попарно несумісні й $A_1$+$A_2$+...+$A_n$=U, де U — достовірна подія.

P(U)=P($A_1$+$A_2$+...+$A_n$) = P($A_1$)+P($A_2$)+...+P($A_n$)=1

• **Приклад 1.** В урні 6 білих кульок, 8 — чорних, 14 — у смужку та 10 — у клітинку. Знайди ймовірність того, що з урни буде вилучена одноколірна кулька.

Відповідь: Спочатку обчислимо загальну кількість кульок в урні:

6+8+14+10 =38. Тоді ймовірність того, що випаде біла кулька P(Б)=$\frac{6}{38}$, а ймовірність вилучення чорної кульки P(Ч)=$\frac{8}{38}$. Ми знаємо що «Ймовірність суми двох несумісних подій дорівнює сумі ймовірностей цих подій, тобто: P(A+B)=P(A)+P(B)», таким чином, ймовірність того, що вилучать однорідну кульку, тобто, білу чи чорну буде:

$$P(Б+Ч)=\frac{6}{38}+\frac{8}{38}=\frac{14}{38}=\frac{7}{19}$$

• **Приклад 2.** Є 100 лотерейних квитків. Відомо, що на 8 квитків трапляється виграш по 20 *грн*, на 11 — по 15 *грн*, на 13 — по 10 *грн*, на 26 — по 2 *грн*, на решту — нічого. Знайди ймовірність того, що на куплений квиток буде отримано виграш не менш ніж 10 *грн*.

Відповідь: Нас цікавить тільки виграш не менше ніж 10 грн, ця фраза означає, що це може бути 20 грн, 15 грн, 10 грн., тобто, маємо 8+11+13=32 наслідку, які сприяють події, що на куплений квиток буде отримано виграш не менш ніж 10 *грн*. Таким чином,

$$P=\frac{32}{100}=\frac{16}{50}=0{,}32.$$

Також цю задачу можна вирішити трохи інакше – знайти ймовірності окремо, а потім їх скласти:

$$P = \frac{8}{100}+\frac{11}{100}+\frac{13}{100}=\frac{32}{100}=0{,}32$$

• **Приклад 3.** Ймовірність того, що студент складе іспит на відмінно дорівнює 0,2; на добре — 0,3; на задовільно — 0,1; на незадовільно — 0,4. Знайди ймовірність того, що студент складе іспит.

Відповідь: При підрахунку ймовірності треба не враховувати тільки не задовільно – 0,4, тобто, усі інші ймовірності маємо скласти (або відняти від одиниці ймовірність «незадовільно») і отримаємо відповідь:
$$P = 0{,}2 + 0{,}3 + 0{,}1 = 1 - 0{,}4 = 0{,}6$$

• **Приклад 4.** У коробці 250 лампочок, із них 61 — по 100 *Вт*, 33 — по 60 *Вт*, 70 — по 25 *Вт*, 86 — по 15 *Вт*. Обчисли ймовірність того, що потужність будь-якої взятої навмання лампочки не перевищить 60 *Вт*.

Відповідь: Фраза «Потужність будь-якої взятої навмання лампочки не перевищить 60 *Вт*.» означає, що треба врахувати лампочки по 60 Вт.; 25 Вт. і 15 Вт., їх кількість це кількість наслідків, сприятливих події, ймовірність якої треба обчислити, тобто, маємо 33+70+85=188 лампочки.

Ймовірність того, що потужність будь-якої взятої навмання лампочки не перевищить 60 *Вт*.:
$$P = \frac{188}{250} = \frac{94}{125} = 0{,}752.$$

• **Приклад 5.** Дві фабрики випускають однакове скло для автомобільних фар. Перша фабрика випускає 36% цього скла, друга — 64%. Перша фабрика випускає 5% бракованого скла, а друга — 3%. Обчисли ймовірність того, що випадково куплене в магазині скло виявиться бракованим.

Відповідь: Спочатку обчислимо скільки бракованого скла випускає кожна фабрика окремо. Так перша фабрика випускає брак 5% від 36%, а друга 3% від 64%. Можемо скласти пропорції і порахувати процент браку від всієї продукції для кожної фабрики, тобто, ми приведемо кількість браку до одиниці продукції, щоб мати змогу порівняти. Маємо:

| 36 | 100% | та | 64 | 100% |
| 5 | $x_1$ | | 3 | $x_2$ |

звідки знаходимо:

$$x_1 = \frac{500}{36} \approx 13{,}9\% \quad \text{та} \quad x_2 = \frac{300}{64} \approx 4{,}7\%.$$

Таким чином, ми з'ясували процент бракованої продукції на обох фабриках, це 13,9%+4,7%=18,6% від всієї продукції на обох фабриках, тобто, від 100%. Ймовірність того, що випадково куплене скло виявиться бракованим:

$$P = \frac{18{,}6}{100} = 0{,}186$$

**1.6. Залежні та незалежні події. Множення ймовірностей.**

**Означення 1.** Подія *A* називається **незалежною** від події *B*, якщо ймовірність появи події *A* не залежить від того, відбулась подія *B* чи не відбулась.

**Означення 2.** Подія *A* називається *залежною* від події *B*, якщо ймовірність події *A* змінюється в залежності від того, відбулась подія *B* чи не відбулась.

Тобто, події *A* і *B* називаються незалежними, якщо поява однієї з них не змінює ймовірності появи іншої. Часто про незалежність подій вдається зробити висновок на підставі того, як організований дослід, у якому вони відбуваються. Незалежні події виникають тоді, коли дослід складається з декількох незалежних випробувань (наприклад, як у досліді з киданням двох гральних кубиків).

**Теорема.** Якщо випадкові події *A* і *B* незалежні, то імовірність суміщення подій *A* і *B* дорівнює добутку ймовірностей появи цих подій.

$$P(AB) = P(A) \cdot P(B).$$

• **Приклад 1.**

Кидають дві гральні кості (кубики). Визначити ймовірність того, що:

а) сума числа очок не перевищує *N*;

б) добуток числа очок не перевищує *N*;

в) добуток числа очок ділиться на *N*.

Відповідь:

**а)** Для *N*=3. Якщо сума очок не перевищує 3, то маємо три комбінації випадання очок на двох гральних картах. Це 1+1, 1+2, 2+1.

Позначимо подію випадання будь-якого фіксованого числа на першому кубику, як «А», тоді випадання (можливо, іншого) фіксованого числа на другому кубику – це подія В.

На гранях гральних костей, маємо від 1 до 6 очок, тоді ймовірність випадання 1 на одному кубику буде:

$$P(A) = \frac{\text{кількість наслідків, сприятливих події A}}{\text{кількість усіх можливих наслідків}}$$, чи $P(A) = \frac{M_t}{N_t}$

$P(A) = \frac{1}{6}$ та $P(B) = \frac{1}{6}$

Але, по умовам задачі, маємо дві гральні кості. Тобто, повинні відбутися обидві події, тому застосуємо формулу ймовірності добутку незалежних подій

$$P(AB) = P(BA) = \frac{1}{6} \cdot \frac{1}{6} = \frac{1}{36}$$

Це ймовірність будь-якої із 36 різних можливих послідовностей чисел (напр., (1,2) та (2,1) – це різні послідовності.

Ймовірність того, що сума чисел не перевищить числа 3 (має три сприятливих фактора) до числа усіх вірогідних можливостей =36, тому:

$$P(N = 3) = \frac{3}{36} = \frac{1}{12}$$

**б)** добуток числа очок не перевищує 3, тобто, $N \leq 3$.

В цьому випадку можливі такі комбінації очок на двох гральних к.:

$1 \cdot 1; 1 \cdot 2; 2 \cdot 1; 3 \cdot 1; 1 \cdot 3$. Тобто, маємо 5 варіантів, що є сприятливими до події, при загальній кількості комбінацій =36 тому:

$$P(AB_{N \leq 3}) = \frac{5}{36}$$

**в)** добуток числа очок ділиться на три (N=3).

Добуток очок, що ділиться на 3 складають 20 пар:

$1 \cdot 3; \ 3 \cdot 1; \ 2 \cdot 3; \ 3 \cdot 2; \ 3 \cdot 3; \ 4 \cdot 3; \ 3 \cdot 4; \ 5 \cdot 3; \ 3 \cdot 5;$
$1 \cdot 6; \ 6 \cdot 1; \ 2 \cdot 6; \ 6 \cdot 2; \ 3 \cdot 6; \ 6 \cdot 3; 4 \cdot 6; \ 6 \cdot 4; \ 5 \cdot 6; \ 6 \cdot 5; \ 6 \cdot 6$

Тобто, маємо двадцять сприятливих наслідків до події, тоді ймовірність цього:

$$P(AB_{:3}) = \frac{20}{36} = \frac{5}{9}$$

• **Приклад 2.** У компанії 12 акціонерів. Із них троє мають привілейовані акції. На збори акціонерів з'явилося 7 осіб. Обчисли ймовірність P(A) того, що серед акціонерів, які з'явилися, усі троє акціонерів з привілейованими акціями відсутні.

Відповідь: Спочатку треба з'ясувати яка ймовірність того, що один із сьомі присутніх акціонерів не є з привілейованими акціями. Усього акціонерів 12 та 3 з привілейованими акціями, тобто, 12-3=9 акціонерів, що нас цікавлять. Ймовірність, що один із цих 9 акціонерів з'явиться P=$\frac{9}{12}$ . Крім того, нам треба з'ясувати ймовірність присутності наступного акціонера, це буде P=$\frac{8}{11}$. Так підраховуємо ймовірності присутності кожного наступного акціонера до 7, тому що їх всього прийшло семеро, тобто, при кожному наступному підрахунку ми зменшуємо кількість наслідків, сприятливих події, тобто, зменшуємо кількість акціонерів без привілейованих акцій та загальну кількість акціонерів. Таким чином, ймовірність присутності останнього сьомого акціонера – P=$\frac{3}{6}$ . Це ми з'ясували ймовірності присутності кожного із сьомі акціонерів окремо, а в завданні треба з'ясувати ймовірність присутності усіх сімох акціонерів одночасно, тобто, всі сім подій повинні відбутися одночасно.

Згадаємо визначення «**Я**кщо випадкові події *A* і *B* незалежні, то імовірність суміщення подій *A* і *B* дорівнює добутку ймовірностей появи цих подій P(AB)=P(A)·P(B)», тобто, ми маємо перемножити всі події.

$$P(A) = \frac{9}{12} \cdot \frac{8}{11} \cdot \frac{7}{10} \cdot \frac{6}{9} \cdot \frac{5}{8} \cdot \frac{4}{7} \cdot \frac{3}{6} = \frac{1}{22}$$

• **Приклад 3.** У коробці лежить 25 краваток, причому 10 із них — червоні, а інші — білі. Обчисли ймовірність того, що з 4 навмання витягнених краваток всі вони будуть одного й того самого кольору.

Відповідь: Ймовірність того, що всі 4 витягнуті краватки будуть тільки червоні буде меншою від ймовірності того, що всі витягнуті краватки будуть білі, тому що білих в коробці більше – 15 шт. Нехай ймовірність події А це витягнута краватка червоного кольору. Тоді:

P($A_1$)=$\frac{10}{25}$; P($A_2$) = $\frac{9}{24}$; P($A_3$) = $\frac{8}{23}$; P($A_4$) = $\frac{7}{22}$,

а ймовірність того, що будуть витягнути всі чотири красні краватки це добуток цих чотирьох ймовірностей:

$$P(A) = \frac{10}{25} \cdot \frac{9}{24} \cdot \frac{8}{23} \cdot \frac{7}{22} = \frac{21}{1265} = 0{,}0166.$$

Тепер аналогічно пообчислюємо ймовірність події В того, що буде витягнута біла краватка. Маємо:
$$P(B_1)=\frac{15}{25};\ P(B_2)=\frac{14}{24};\ P(B_3)=\frac{13}{23};\ P(B_4)=\frac{12}{22} \quad \text{та}$$
$$P(B)=\frac{15}{25}\cdot\frac{14}{24}\cdot\frac{13}{23}\cdot\frac{12}{22}=\frac{273}{2530}\approx 0{,}108$$

### 1.7. Статистичне та класичне означення ймовірності.

Статистичною (або вибірковою) ймовірністю події $A$ називається відношення кількості $m$ випробувань, в яких подія $A$ відбулась, до загальної кількості виконаних випробувань $n$: $W(A)=m/n$

Теоретичною ймовірністю випадкової події $A$ називається відношення кількості елементарних подій $m$, які сприяють появі цієї події (становлять множину її елементарних подій), до загальної кількості $n$ рівно-можливих елементарних подій, що утворюють простір елементарних подій $W$ (чи $\Omega$): $P(A)=m/n$.

Якщо $m=0$, то така подія є неможливою, оскільки нема жодної елементарної події, яка б була сприятливою для цієї події, і за формулою її ймовірність дорівнює нулю. З точки зору термінології, побутове «не відбулось жодної події» = математичне «відбулось 0 подій».

Щоб обчислити ймовірність події $A$ за цією формулою, потрібно знайти кількість елементарних подій у просторі $W$, а також кількість їх у множині, яка відповідає події $A$. Для цього застосовують формули комбінаторики.

### 1.8. Геометричне визначення ймовірності.

У випробуваннях з незчисленною кількістю результатів для підрахунку ймовірностей вводять поняття геометричної ймовірності.

Нехай в результаті випробування точка з однаковою ймовірністю обов'язково може потрапити в будь-яке місце деякої області. Вважаємо, що подія відбулася, якщо точка при цьому потрапляє в область. Геометричною ймовірністю події називають ймовірність, обчислену за формулою: ,

$P(A)=\frac{\text{міра } A}{\text{міра } \Omega}=\frac{\mu(A)}{\mu(\Omega)}$, де $\mu(A)$ – це міра області $A$, яку звичайно позначають грецькою літерою «мю» ($\mu$).

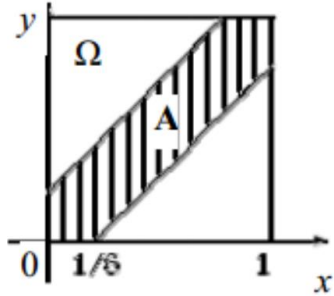

Тобто, під мірою розуміють довжину, площу, об'єм в одно -, двох - і тривимірному випадках відповідно. Очевидно, це можна узагальнити й на багатовимірний простір (чисельний чи, за ознаками у базі даних із багатьма параметрами.

Таким чином, можемо записати наступні рівняння:

1. Якщо відрізок $l$ складає частину відрізка L і на відрізок L навмання поставлена точка, то при допущенні, що ймовірність потрапляння точки на відрізок $l$ пропорційна довжині цього відрізку і не залежить від його розміщення відносно відрізку L маємо записати ймовірність потрапляння точки на відрізок $l$:

$$P = \frac{\text{Довжина } l}{\text{Довжина } L}.$$

2. Якщо плоска фігура $g$ складає частину плоскої фігури G і на плоску фігуру G навмання поставлена точка, то при допущенні, що ймовірність потрапляння точки на фігуру $g$ пропорційна площі цієї фігури і не залежить ні від її розміщення відносно фігури G, ні від форми $g$ маємо записати ймовірність потрапляння точки в фігуру $g$:

$$P = \frac{\text{Площа } g}{\text{Площа G}}.$$

3. Аналогічно визначається ймовірність потрапляння точки в об'ємну фігуру $v$, яка складає частину об'ємної фігури $V$:

$$P = \frac{\text{Об'єм } v}{\text{Об'єм } V}$$

• **Приклад 1.** На відрізку L довжиною 20 см розміщений менший відрізок l довжиною 10 см. Знайти ймовірність того, що точка, яка довільно розміщена на великому відрізку, попаде також і на менший відрізок.

Відповідь:

*l=10 см*

*L=20 см*

$$P = \frac{\text{Довжина } l}{\text{Довжина } L} = \frac{10}{20} = \frac{1}{2} = 0{,}5 \text{ см}$$

• **Приклад 2.** У відрізку одиничної довжини навмання з'являється точка. Визначити ймовірність того, що точка буде у одному з двох діапазонів [0.1,0.3] або [0.2, 0.4].

Відповідь:

Маємо задачу на геометричне визначення ймовірності.

$P(A) = \frac{\text{міра } A}{\text{міра } \Omega} = \frac{\mu(A)}{\mu(\Omega)}$, де міра $\Omega$ – це одинична довжина, а подія $A = A_1 + A_2$, де елементарні події – попадання в перший, або другий інтервал. Міра $\mu(A) = \mu(A_1) + \mu(A_2) - \mu(A_1 A_2) = \mu(0.1 \leq x \leq 0.3) + \mu(0.2 \leq x \leq 0.4) - \mu(0.2 \leq x \leq 0.3) = (0.3 - 0.1) + (0.4 - 0.2) - (0.3 - 0.2) = 0.2 + 0.2 - 0.1 = 0.3 = \mu(0.1 \leq x \leq 0.4)$.

Отже, відповідь $P(A) = \frac{\mu(A)}{\mu(\Omega)} = \frac{0.3}{1} = 0.3$.

*Зауваження.* Слід відмітити, що такий результат справедливий для неперервних координат $x$. При дискретизації, слід бути уважним, чи входять у множину межі інтервалів. Наприклад, якщо округляємо до десятих, то в перший інтервал входять лише числа 0.1, 0.2, 0.3, отже, 3 числа. Множина $\Omega$ теж може бути під питанням – це від 0.0 до 1.0, тобто, 11 чисел? Тоді $P(A) = \frac{\mu(A)}{\mu(\Omega)} = \frac{3}{11} = 0.2727 \ldots$ . Якщо числа до сотих, то $P(A) = \frac{\mu(A)}{\mu(\Omega)} = \frac{31}{101} = 0{,}3069 \ldots$ ближче до отриманого нами правильного значення 0.3. Чим буде більше знаків, тим буде ближче до правильного значення, яке є межею відношення для кількості знаків, що прямує до нескінченності. Для виправлення цієї неточності, часто роблять інтервал «піввідкритим», напр. (0.1, 0.399999), але 0.4 в нього не входить. На побутовому рівні, за хвилину після 23:59, буде не 24:00, а 00:00 наступної доби. Багато цін ставлять, як на «дрібничку» менше, ніж основна ціна, напр. 99.99, а не 100. Інше зауваження в тому, що «дискретна» «міра» насправді може співпадати із «неперервною. Напр. до 0.0 будуть округлені числа $0.0 \leq x \leq 0.04999 \ldots$, до 0.1 – числа $0.05 \leq x \leq 0.14999 \ldots$, Тобто, крайні інтервали вдвічі менші за інтервали в середині. Якщо зберігати кількість дискретних чисел рівною кількості інтервалів, то не потрібно включати значення 1.0. Тобто, «обрізати» (trunc) числа до почату інтервалу, або брати середину інтервалу (0.05, 0.15, …, 0.95). Ці нюанси потрібно враховувати при програмуванні. Наприклад, у популярній мові Python, range(n) дає вектор із компонентами 0,1,..(n-1), що аналогічно linspace(0,0.9,10) – обидва не дають 1.

• **Приклад 3.** Нехай в нас є поле розміром 50 на 100 (м). Перед нами стоїть ціль, яка ймовірність того, що парашутист попадав в коло на цьому полі радіусом 10(м), якщо він може приземлитися будь-де у цьому прямокутніку?

Відповідь: $\mu(A)$ – площа кола, тобто, $\mu(A) = S_{\text{коло}} = \pi R^2 = \pi \cdot 10^2 = \pi \cdot 100$

$\mu(\Omega)$ – площа прямокутника, тобто, $\mu(\Omega) = S_{\text{пр.}} = 50 \cdot 100$; Тоді:

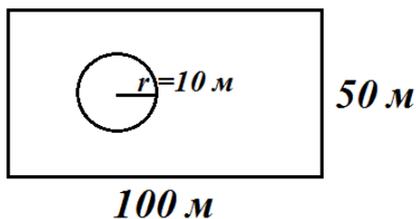

$$P(A) = \frac{\mu(A)}{\mu(\Omega)} = \frac{\pi \cdot 10^2}{50 \cdot 100} = \frac{\pi}{50} = 0.0628$$

• **Приклад 4.** На аудіокасеті записані концерти трьох співаків: першого -

протягом 40 хв. Звучання, другого – протягом 30 хв., третього – протягом 20 хв. Запис перемотується і навмання включається. Яка ймовірність, що звучить пісня у виконанні другого співака?

Відповідь: Час звучання всього запису $T(\Omega) = 90$ хв., час звучання другого співака $T(A) = 30$ хв. За формулою $P(A) = \frac{\text{міра } A}{\text{міра } \Omega}$ маємо:

$$P(A) = \frac{T(A)}{T(\Omega)} = \frac{30}{90} = \frac{1}{3}.$$

• **Приклад 5.**

<u>1 варіант задачі.</u> У коло радіусу R розміщено менше коло радіусу r. Знайти ймовірність того, що точка, яка навмання кинута в велике коло, попаде також і в менше коло.

Відповідь:

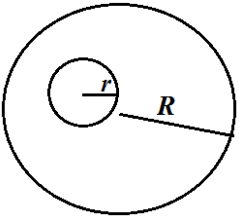

За формулою $P = \dfrac{\text{Площа } g}{\text{Площа } G} = \dfrac{\pi r^2}{\pi R^2} = \dfrac{r^2}{R^2}$

<u>2 варіант задачі.</u>

На площині накреслені два концентричних кола з радіусами 5 та 10 см відповідно. Знайти ймовірність того, що точка, яка кинута навмання у велике коло, потрапить також і в кільце, що утворилося побудованими колами.

Відповідь:

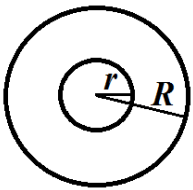 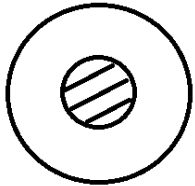 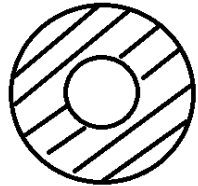

$S_1 = \pi r^2$
$S_2 = \pi R^2$

Спочатку обчислимо площу утвореного кільця (заштриховане):
$$S_\text{к} = S_2 - S_1 = \pi R^2 - \pi r^2 = \pi(R^2 - r^2)$$

Підставимо знайдений вираз площі кільця у формулу $P = \dfrac{\text{Площа } g}{\text{Площа } G}$, де $g$ площа, що в задачі відповідає фігурі в яку потрапляє точка, тобто:
$$S_\text{к} = \pi(R^2 - r^2),$$

а площа $G$ це $S_2 = \pi R^2$, тобто, маємо ймовірність потрапляння в кільце:
$$P(k) = \frac{\pi(R^2 - r^2)}{\pi R^2} = \frac{R^2 - r^2}{R^2} = 1 - \frac{r^2}{R^2} = 1 - \frac{25}{100} = 1 - 0{,}25 = 0{,}75$$

• **Приклад 6.** На площину з нанесеною сіткою квадратів зі

стороною $a$ навмання кинута монета радіусом $r < \frac{a}{2}$. Знайти ймовірність того, що монета не перетне ні однієї із сторін квадрата. Передбачається, що ймовірність потрапляння точки на фігуру пропорційна площі цієї фігури і не залежить від її розміщення.

Відповідь: Щоб краще зрозуміти умови завдання спочатку побудуємо рисунок так, що $a > 2r$ (згідно умови $r < \frac{a}{2}$).

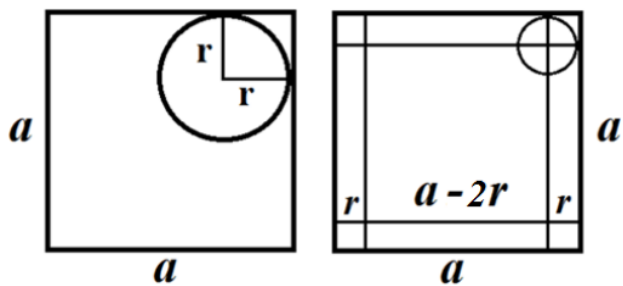

Згідно з другою умовою задачі коло не повинно торкатися ні однієї зі сторін квадрата. Практично це означає, що коло повинно бути на відстані свого радіусу від сторони даного квадрата. Щоб визначити площу фігури, в якій може знаходитися центр кола, поглянемо на другий рисунок.

Ми отримали смугу шириною $r$ вздовж периметру даного квадрата, в яку не може потрапити центр кола і менший квадрат по центру даного зі сприятливими умовами потрапляння, тобто, площа цього меншого квадрату і є $\mu(A)$. Сторону меншого квадрата маємо записати як: $(a - 2r)$, тоді площа фігури буде $(a - 2r)^2$. За умовою задачі $\mu(\Omega) = S_{\text{кв.}} = a^2$.

Таким чином отримаємо:
$$P(A) = \frac{\mu(A)}{\mu(\Omega)} = \frac{(a-2r)^2}{a^2}.$$

• **Приклад 7.** Довільно взято два додатних числа X та Y, кожне з яких не перевищує одиниці. Знайти ймовірність того, що сума X + Y не перевищить одиниці (тобто, менше або рівно 1), а добуток X·Y не менше 0,09 (тобто, більше або рівно 0,09).

Відповідь. По умовам задачі маємо інтервали $0 < X \leq 1$; $0 < Y \leq 1$ (тобто, X та Y два додатних числа, які не перевищують одиниці). Крім того X+Y≤ 1 та X·Y≥ 0,09.

Для розв'язання цієї задачі застосуємо декартову систему координат і на неї побудуємо графіки, отриманих умов. Для наочності запишемо умови наступним чином:

замість X·Y≥ 0,09 запишемо $Y = \frac{0,09}{X}$;

замість X+Y≤ 1 запишемо $Y = 1 - X$.

Згадаємо, що графіком функції $Y = \frac{1}{X}$ це є гіпербола, а приймаючи до уваги додаткові обмеження: $0 < X \leq 1$; $0 < Y \leq 1$, маємо розглядати тільки додатну частину графіка.

Графік функції $Y = 1 - X$ це є пряма, яку треба побудувати на цьому ж графіку.

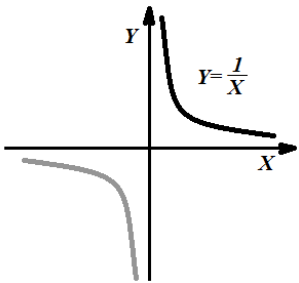 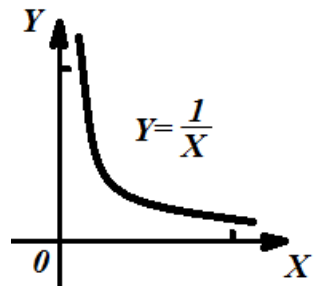

Таким чином, врахувавши усі умови задачі, ми отримали графічно сектор (заштрихована частина), площа якого, є сприятливим фактором при знаходженні ймовірності, де сума X + Y не перевищує одиниці а добуток X·Y не менше 0,09. Тепер необхідно знайти чисельно цю площу – позначену як $S_1$.

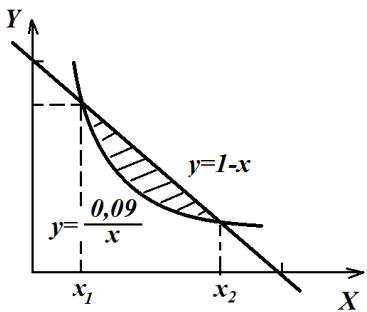 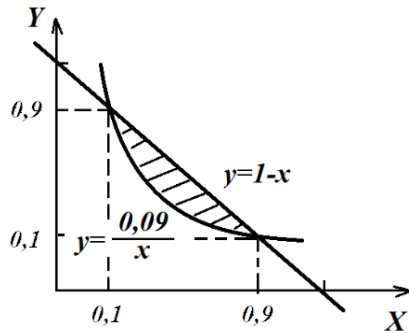

Застосуємо визначений інтеграл для знаходження площі фігури обмеженої лініями, в нашому випадку кривою $Y = \frac{0,09}{X}$ і прямою $Y = 1 - X$, тобто, знайдемо площу за формулою:

$$S = \int\limits_{x_1}^{x_2} \bigl(f(x) - g(x)\bigr)dx$$

де функція $f(x)$ це $Y = 1 - X$, а $g(x)$ це $Y = \frac{0{,}09}{X}$ (віднімати треба із тієї функції, яка графічно розташована вище).

Щоб знайти межі інтегрування треба визначити точки перетину кривої $Y = \frac{0{,}09}{X}$ з прямою $Y = 1 - X$, для чого маємо вирішити систему рівнянь:

$$\begin{cases} Y = \dfrac{0{,}09}{X} \\ Y = 1 - X \end{cases}$$

$$\frac{0{,}09}{X} = 1 - X$$

$$0{,}09 = X(1 - X) = X - X^2$$

$$X^2 - X + 0{,}09 = 0$$

Отримали квадратне рівняння, знайдемо корні цього рівняння за формулами:

$$ax^2 + bx + c = 0 \quad \text{тоді} \quad = \frac{-b \pm \sqrt{b^2 - 4ac}}{2a}$$

$X^2 - X + 0{,}09 = 0.$

$$x_{1,2} = \frac{-(-1) \pm \sqrt{(-1)^2 - 4 \cdot 1 \cdot 0{,}09}}{2 \cdot 1} = \frac{1 \pm \sqrt{1 - 0{,}36}}{2} = \frac{1 \pm 0{,}8}{2}$$

Отримали:

$$x_1 = \frac{1 - 0{,}8}{2} = \frac{0{,}2}{2} = 0{,}1$$

$$x_2 = \frac{1 + 0{,}8}{2} = \frac{1{,}8}{2} = 0{,}9$$

Таким чином, визначений інтеграл має межі від 0,1 до 0,9, тобто, ми отримали додаткові значення координат, в межах яких є «сприятливі значення»: $0{,}1 \leq x \leq 0{,}9$ та $0{,}1 \leq y \leq 0{,}9$ за формулою:

$$S = \int\limits_{x_1}^{x_2} \bigl(f(x) - g(x)\bigr)dx \quad \text{маємо записати:}$$

$$S_1 = \int\limits_{0,1}^{0,9} \left((1-x) - \left(\frac{0,09}{x}\right)\right)dx = \int\limits_{0,1}^{0,9}(1-x)dx - \int\limits_{0,1}^{0,9}\frac{0,09}{x}dx =$$

$$= \int\limits_{0,1}^{0,9} 1 dx - \int\limits_{0,1}^{0,9} x dx - 0,09 \cdot \int\limits_{0,1}^{0,9}\frac{1}{x}dx$$

$$= x\Big|_{0,1}^{0,9} - \frac{x^2}{2}\Big|_{0,1}^{0,9} - 0,09\ln x\Big|_{0,1}^{0,9} =$$

$$= 0,8 - 0,4 - 0,198 = 0,202$$

За умовою задачі «Довільно взято два додатних числа X та Y, кожне з яких не перевищує одиниці», тобто, мається на увазі, що це є квадрат зі сторонами 1x1, тобто, площа $S_2 = 1 \cdot 1 = 1$. Зараз ми вже можемо з'ясувати ймовірність:

$$P = \frac{S_1}{S_2} = \frac{0,202}{1} = 0,202$$

**Приклад 8.** Випадкова точка попадає у квадрат із одиничною стороною. Яка ймовірність того, що її відстань від лівого нижнього кута не перевищує 1?\

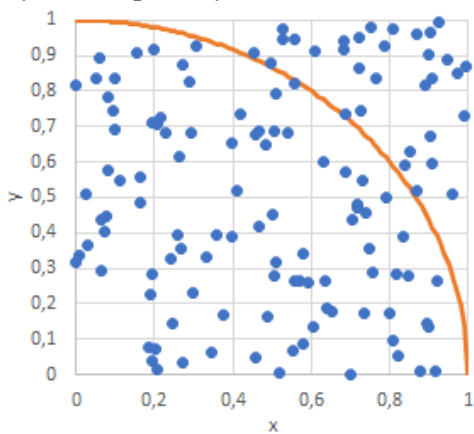

Введемо систему координат, як показано на рис. Площа квадрата $S_n = 1$, площа чверті круга $S_m = \pi/4$. Теоретична геометрична ймовірність $p = S_m/S_n = \pi/4 \approx 0{,}78539816$. З $n = 200$ точок із згенерованими випадковими координатами, $m = 155$ точок задовольняють умові $x^2 + y^2 \leq 1$, отже, вибіркова ймовірність $p_{виб} = m/n = 155/200 = 0.775$. Математичне очікування (теоретичне значення) кількості подій $M = pn \approx 157.08$ не є цілим, хоча вибіркові значення, які можуть відрізнятися від теоретичного та між собою, є цілими, й є випадковими числами, які підкоряються біноміальному розподілу, який буде розглянуто окремо. Хоча точки мають випадкові координати, вони, для кожної вибірки, створюють окремі групи (кластери), аналогічно «сузір'ям» на небі. Це відрізняється від детермінованого рівномірного розподілу, де точки б «вишукувались» у рівномірні рядки та стовпчики. Дослідження властивостей згенерованих точок називають «методом Монте-Карло».

## 2. ЕЛЕМЕНТИ КОМБІНАТОРИКИ.

Комбінаторика — розділ математики, в якому вивчаються способи вибору і розміщення елементів деякої скінченної множини на основі якихось умов. Вибрані (або вибрані і розміщені) групи елементів називають сполученнями або вибіркою.

Якщо всі елементи заданої впорядкованої множини різні — отримаємо перестановки без повторень, а якщо в заданій упорядкованій множині елементи повторюються, то дістанемо перестановки з повтореннями.

### 2.1. Перестановки.

Перестановкою з n елементів називається будь-яка впорядкована множина з n елементів. Інакше кажучи, це така множина, для якої указано, який елемент знаходиться на першому місці, який — на другому, ..., який — на n-му).

Формула числа перестановок ($P_n$):

$(P_n) = n!$

де n! = 1·2·3·... n (знак ! це факторіал).

Однією із найважливіших властивостей факторіалу є

$$n! = n(n-1)! = n(n-1)(n-2)! = \cdots$$
$$= n(n-1)\ldots(n-k+1)(n-k)!$$

для будь-якого цілого $0 \leq k \leq n$. Прийнято, що 0!=1. Адже, коли добуток у цій формулі помножується на останній вираз $(n-k)!$ при $n = k$, то це мусить бути одиниця.

• **Приклад:** Кількість різних шестизначних чисел, які можна скласти з цифр 1, 2, 3, 4, 5, 6, не повторюючи ці цифри в одному числі, дорівнює

$P_6$ = 6! = 1·2·3·4·5·6 = 720.

_Перестановки з повтореннями_. Число різних перестановок, які можна утворити з $n$ елементів, серед яких є $k_1$ елементів першого типу, $k_2$ елементів другого типу, ...., $k_m$ елементів $m$-того типу, дорівнює:

$$P_n(k_1, k_2, \ldots, k_m) = \frac{n!}{k_1! \cdot k_2! \cdot \ldots \cdot k_m!}$$

• **Приклад:** Скільки різних слів (перестановок) можна утворити перестановкою літер у слові "математика"?

Відповідь: Слово "математика" містить n = 10 літер. Літер одного типу:

"м" - дві, "а" - три, "т" - дві, "е" - одна, "и" - одна, "к" - одна. За формулою:

$$P_n(k_1, k_2, \ldots, k_m) = \frac{n!}{k_1! \cdot k_2! \cdot \ldots \cdot k_m!}$$

одержимо:

$$P_{10}(2,3,2,1,1,1) = \frac{10!}{2! \cdot 3! \cdot 2! \cdot 1! \cdot 1! \cdot 1!!} = 151200.$$

Оскільки серед цих слів, лише одне є «математика», то ймовірність випадкового набору цього слова буде дорівнювати 1/151200.

### 2.2. Розміщення.

Розміщенням (перестановкою, permutation) з $n$ елементів по $k$ називається будь-яка впорядкована множина з $k$ елементів, складена з елементів $n$-елементної множини.

Формула числа розміщень $\left(A_n^k\right)$ – (читається: «А з $n$ по $k$»):

$$A_n^k = n(n-1)\ldots(n-k+1) = \frac{n!}{(n-k)!}$$

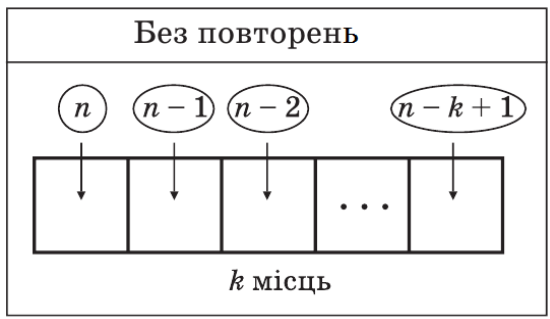

Пояснення: Розглянемо, скільки всього можна скласти розміщень з $n$ елементів по $k$ (без повторень). Складання розміщення уявимо собі як послідовне заповнення $k$ місць, які ми будемо зображати у вигляді клітинок (див. рис.). На перше місце ми можемо вибрати один з n елементів заданої множини (тобто, елемент для першої клітинки можна вибрати n способами).

Якщо елементи не можна повторювати, то на друге місце можна вибрати тільки один елемент із тих, що залишилися, тобто, з (n – 1). Тепер вже два елементи використані і на третє місце можна вибрати тільки один з (n – 2) елементів і т. д. На k-те місце можна вибрати тільки один з n – (k –1) = n – k +1 елементів.

Оскільки нам потрібно вибрати елементи і на перше місце, і на друге, ..., і на k-те, то використовуємо правило добутку і одержуємо таку формулу числа розміщень з n елементів по k:

$$A_n^k = \underbrace{n(n-1)(n-2)\ldots(n-k+1)}_{k \text{ множників}}$$

За допомогою факторіалів, цю формулу для числа розміщень без повторень можна записати в іншому вигляді. Для цього помножимо і поділимо вираз у цієї формулі на добуток $(n-k)\cdot(n-k-1)\cdot\ldots\cdot 2\cdot 1 = (n-k)!$. Одержуємо:

$$A_n^k = n(n-1)(n-2)\ldots(n-k+1) =$$
$$= \frac{n\cdot(n-1)\cdot(n-2)\cdot\ldots\cdot(n-k+1)\cdot(n-k)\cdot(n-k-1)\cdot\ldots\cdot 3\cdot 2\cdot 1}{(n-k)\cdot(n-k-1)\cdot\ldots\cdot 3\cdot 2\cdot 1} = \frac{n!}{(n-k)!}.$$

Отже, формула для числа розміщень без повторень з n елементів по k може бути записана так:

$$A_n^k = \frac{n!}{(n-k)!}$$

Для того щоб цією формулою можна було користуватися при всіх значеннях $k$ і, зокрема при $k = n - 1$ та при $k = n$, домовилися вважати, що

$1! = 1$ і $0! = 1$.

Отже,
$$A_n^n = n!, \quad A_n^0 = 1, \ A_n^k = (n-k+1)A_n^{k-1}$$

• **Приклад**: Кількість різних тризначних чисел, які можна скласти з цифр 1, 2, 3, 4, 5, 6, якщо цифри не можуть повторюватися, дорівнює:

$$A_6^3 = \frac{6!}{(6-3)!} = \frac{6!}{3!} = \frac{1\cdot 2\cdot 3\cdot 4\cdot 5\cdot 6}{1\cdot 2\cdot 3} = 4\cdot 5\cdot 6 = 120.$$

### Розміщення з повтореннями.

Нехай $M$ множина, що містить n елементів. Розміщенням з повтореннями з $n$ елементів по $k$ називається довільна впорядкована підмножина з $k$ елементів з множини $M$ (елементи не обов'язково різні). Число розміщень з повтореннями з $n$ елементів по $k$ знаходиться за формулою:

$$\widetilde{A_n^k} = n^k$$

Тобто, $\widetilde{A_n^k}$ – це число способів, якими можна розкласти k різних предметів по n ящиках.

• **Приклад 1.** Є 4 події ABCD, які можуть відбуватися чи ні. Класична форма запису для «не A» є із рискою зверху ($\bar{A}$). Для

запису «у рядок» (без формульного редактора), записують \A. Втім, такі впорядковані послідовності зручно записувати через 0 (0 подій=подія не відбулася) чи 1 (відбулася 1 подія). Така нумерація зручна, щоб не пропустити якусь послідовність, адже можна розглядати послідовність 0 та 1, як запис числа у двійковій системі числення.

| 16-річні | 10-річні | ABCD | Класично, Риска | у рядок | Кількість подій |
|---|---|---|---|---|---|
| 0 | 0 | 0000 | $\bar{A}\,\bar{B}\,\bar{C}\,\bar{D}$ | \A\B\C\D | 0 |
| 1 | 1 | 0001 | $\bar{A}\,\bar{B}\,\bar{C}\,D$ | \A\B\C D | 1 |
| 2 | 2 | 0010 | $\bar{A}\,\bar{B}\,C\,\bar{D}$ | \A\B C\D | 1 |
| 3 | 3 | 0011 | $\bar{A}\,\bar{B}\,C\,D$ | \A\B C D | 2 |
| 4 | 4 | 0100 | $\bar{A}\,B\,\bar{C}\,\bar{D}$ | \A B\C\D | 1 |
| 5 | 5 | 0101 | $\bar{A}\,B\,\bar{C}\,D$ | \A B\C D | 2 |
| 6 | 6 | 0110 | $\bar{A}\,B\,C\,\bar{D}$ | \A B C\D | 2 |
| 7 | 7 | 0111 | $\bar{A}\,B\,C\,D$ | \A B C D | 3 |
| 8 | 8 | 1000 | $A\,\bar{B}\,\bar{C}\,\bar{D}$ | A\B\C\D | 1 |
| 9 | 9 | 1001 | $A\,\bar{B}\,\bar{C}\,D$ | A\B\C D | 2 |
| A | 10 | 1010 | $A\,\bar{B}\,C\,\bar{D}$ | A\B C\D | 2 |
| B | 11 | 1011 | $A\,\bar{B}\,C\,D$ | A\B C D | 3 |
| C | 12 | 1100 | $A\,B\,\bar{C}\,\bar{D}$ | A B\C\D | 2 |
| D | 13 | 1101 | $A\,B\,\bar{C}\,D$ | A B\C D | 3 |
| E | 14 | 1110 | $A\,B\,C\,\bar{D}$ | A B C\D | 3 |
| F | 15 | 1111 | $A\,B\,C\,D$ | A B C D | 4 |

• **Приклад 2.** У ліфт 12-поверхового будинку зайшло на першому поверсі 10 чоловік. Скількома способами вони можуть вийти з ліфта?

Відповідь: Задача зводиться до розкладання 10 предметів по 11 ящиках.

Число таких способів дорівнює $\widetilde{A}_{11}^{10} = 11^{10}$.

• **Приклад 3.** В ліфт 6-поверхового будинку сіли 4 пасажира. Кожен, незалежно від інших з однаковою ймовірністю, може вийти на будь-якому (починаючи з другого) поверсі. Визначити

ймовірність того, що:

а) всі вийшли на різних поверхах;

б) принаймні двоє зійшли на одному поверсі.

Відповідь: Існує декілька варіантів виходу пасажира (у кожного з n=4). Визначимо кількість варіантів виходу на різних поверхах: в першого пасажира є (k-1) де -1 означає вихід (починаючи з другого поверху), а тому що *n* пасажирів, то маємо *n* разів (кожний повинен вийти на іншому поверху), таким чином, можливість в другого (*k*-2), в третього (*k*-3), в *n*-го (*k*-n), тобто, загалом

$$A_{k-1}^n = \frac{(k-1)!}{(k-1-n)!} = (k-1)(k-2)\ldots(k-n).$$

Загальна ж кількість можливостей $(k-1)^n$, тому ймовірність

$$p_n = \frac{(k-1)\cdot(k-2)\cdot\ldots\cdot(k-n)}{(k-1)\cdot(k-1)\cdot\ldots\cdot(k-1)}.$$

Підставимо наші значення

$$p_4 = \frac{(6-1)(6-2)(6-3)(6-4)}{(6-1)(6-1)(6-1)(6-1)} = \frac{5\cdot 4\cdot 3\cdot 2}{5^4} = \frac{120}{625} = 0{,}192$$

Ми з'ясували ймовірність того, що четверо пасажирів вийдуть на різних поверхах (починаючи з другого) шести-поверхового будинку. Протилежною подією *q* до цієї, є вихід декількох пасажирів на одному поверсі (принаймні, двох пасажирів принаймні на одному поверсі) тому відповідь на б): ймовірність протилежної події:

$$q = 1 - p = 1 - 0{,}192 = 0{,}808$$

**2.3. Комбінації (combinations).**

Комбінацією без повторень з *n* елементів по *k* називається будь-яка

*k*-елементна підмножина n-елементної множини.

Формула числа комбінацій $\left(C_n^k\right)$

$$C_n^k = \frac{n!}{k!\,(n-k)!}$$

За означенням, вважається, що $C_n^0 = 1$.

Пояснення: Розглянемо скільки всього можна скласти комбінацій без повторень з *n* елементів по *k*. Для цього використаємо відомі нам формули числа розміщень і перестановок.

Складання розміщення без повторень з *n* елементів по *k* проведемо в два етапи. Спочатку виберемо *k* різних елементів із заданої *n*-елементної множини, не враховуючи порядок вибору цих елементів (тобто, виберемо $k$-елементну підмножину з $n$-елементної множини — комбінацію без повторень з $n$-елементів по $k$). За нашим позначенням, це можна зробити $C_n^k$ способами. Після цього, одержану множину з $k$ різних елементів впорядкуємо. Її можна впорядкувати $P_k = k!$ способами. Отже, кількість розміщень без повторень з $n$ елементів по $k$ в $k!$ разів більша за число комбінацій без повторень з $n$ елементів по $k$. Тобто, $A_n^k = C_n^k \cdot k!$. Звідси:

$$C_n^k = \frac{A_n^k}{k!} = \frac{n!}{(n-k)!\, k!} = \frac{n!}{k!\,(n-k)!}.$$

Якщо у цієї формули скоротити чисельник і знаменник на $(n-k)!$, то дістанемо формулу, за якою зручно обчислювати $C_n^k$ при малих значеннях $k$:

$$C_n^k = \frac{\overbrace{n(n-1)(n-2)\ldots(n-k+1)}^{k \text{ множників}}}{k!} = \frac{\overbrace{n(n-1)(n-2)\ldots(n-k+1)}^{k \text{ множників}}}{\underbrace{1\cdot 2\cdot\ldots\cdot k}_{k \text{ множників}}}.$$

Наприклад, $C_{25}^2 = \dfrac{\overbrace{25\cdot 24}^{2 \text{ множники}}}{1\cdot 2} = 25\cdot 12 = 300$, $C_8^3 = \dfrac{\overbrace{8\cdot 7\cdot 6}^{3 \text{ множники}}}{1\cdot 2\cdot 3} = 8\cdot 7 = 56$.

Слід зауважити властивість симетрії $C_n^k = C_n^{n-k}$. Отже, перед обчисленням, можна рекомендувати вибрати у якості верхнього індексу менше з чисел $n$ та $n-k$. Наприклад,

$$C_{10}^8 = \frac{A_{10}^8}{8!} = \frac{10\cdot 9\cdot 8\cdot 7\cdot 6\cdot 5\cdot 4\cdot 3}{1\cdot 2\cdot 3\cdot 4\cdot 5\cdot 6\cdot 7\cdot 8} = \frac{10\cdot 9}{1\cdot 2} = 45$$

Значно швидше обчислити, змінив $k$ та $n-k$ та навпаки:

$$C_{10}^8 = C_{10}^{10-8} = C_{10}^2 = \frac{A_{10}^2}{2!} = \frac{10\cdot 9}{1\cdot 2} = 45$$

• **Приклад 1:** Із класу, що складається з 25 учнів, можна виділити 5 учнів для чергування по школі $C_{25}^5$ способами, тобто

$$C_{25}^5 = \frac{25!}{5!(25-5)!} = \frac{25!}{5!\cdot 20!} = \frac{21\cdot 22\cdot 23\cdot 24\cdot 25}{1\cdot 2\cdot 3\cdot 4\cdot 5} = 53\,130$$

способами.

• **Приклад 2**. Серед $n=10$ лотерейних квитків, $k=6$ виграшних. Навмання взяли $m=4$ квитків. Визначити ймовірність того, що серед них $l=2$ виграшних.

Відповідь:

Для вирішення задачі використаємо формулу для числа комбінацій. Запишемо ще раз отримані дані і комбінації: $n = 10$ – загальна кількість квитків із них виграшних («першого ґатунку») $k = 6$. Взято $m = 4$ і із них виграшних $l = 2$. Тобто, число способів взяти 4 квитка із 10 буде комбінація $C_{10}^4$, але, крім 4 виграшних з 6 ($C_6^2$ комбінацій), є два невиграшних з 4 ($C_4^2$ комбінацій). Загальна кількість сприятливих комбінацій $C_6^2 C_4^2$.

Ймовірність дорівнює відношенню кількості сприятливих комбінацій до загальної.

Загальна формула для ймовірності того, що буде взято $l$ об'єктів з $k$:

$$p_l = \frac{C_k^l \cdot C_{n-k}^{m-l}}{C_n^m}.$$

Такий статистичний розподіл називають «гіпергеометричним» розподілом. По формулі числа комбінацій $\left(C_n^k\right)$ пообчислюємо комбінації $C_k^l, C_{n-k}^{m-l}$ та $C_n^m$:

$$C_n^k = \frac{n!}{k!(n-k)!}$$

$C_4^2 = \frac{4!}{2!(4-2)!} = \frac{1 \cdot 2 \cdot 3 \cdot 4}{2 \cdot 2} = 6$

$C_6^2 = \frac{6!}{2!(6-2)!} = \frac{4! \cdot 5 \cdot 6}{2! \cdot 4!} = 15$

$C_{10}^4 = \frac{10!}{4!(10-4)!} = \frac{6! \cdot 7 \cdot 8 \cdot 9 \cdot 10}{4! \cdot 6!} = 210$

Ймовірність дорівнює:

$p_2 = \frac{C_6^2 \cdot C_4^2}{C_{10}^4} = \frac{15 \cdot 6}{210} = \frac{3}{7} = 0{,}429$

• **Приклад 3.** Маємо вироби чотирьох ґатунків, причому число виробів $i$-го ґатунку дорівнює $n_i$, i=1,2,3,4. Для контролю, навмання беруться $m$ виробів. Визначити ймовірність того, що серед них $m_1$ першого ґатунку, $m_2, m_3$ і $m_4$ другого, третього і четвертого ґатунків відповідно.

| $n_1$ | $n_2$ | $n_3$ | $n_4$ | $m_1$ | $m_2$ | $m_3$ | $m_4$ |
|---|---|---|---|---|---|---|---|
| 1 | 2 | 3 | 4 | 1 | 1 | 2 | 3 |

Відповідь:

Для вирішення задачі, скористаємось формулою для числа комбінацій

$$C_n^k = \frac{n!}{k!\,(n-k)!}$$

і запишемо число комбінацій для виробів кожного ґатунку:

Маємо виробів першого ґатунку $m_1 = 1$ при $n_1 = 1$, тобто, $C_1^1 = 1$, із властивостей числа комбінацій $C_n^n = 1$.

Виробів другого ґатунку $m_2 = 1$, при $n_2 = 2$, маємо записати

$$C_2^1 = \frac{2!}{1!\,(2-1)!} = \frac{2!}{1!} = \frac{1 \cdot 2}{1} = 2$$

Виробів третього ґатунку $m_3 = 2$, при $n_2 = 3$, маємо записати

$$C_3^2 = \frac{3!}{2!\,(3-2)!} = \frac{3!}{2!} = \frac{1 \cdot 2 \cdot 3}{1 \cdot 2} = 3$$

Виробів четвертого ґатунку $m_4 = 3$, при $n_2 = 4$, маємо записати

$$C_4^3 = \frac{4!}{3!\,(4-3)!} = \frac{4!}{3!} = \frac{1 \cdot 2 \cdot 3 \cdot 4}{1 \cdot 2 \cdot 3} = 4$$

Таким чином з'ясували число комбінацій без повторень для виробів кожного ґатунку.

Для контролю навмання беруться $m$ виробів і треба визначити ймовірність того, що серед них першого ґатунку 1 деталь, другого 1, третього 2 і четвертого 3 виробів, тобто, усього взято 7. Маємо врахувати всі комбінації одночасно, для чого їх треба перемножити. Щоб з'ясувати ймовірність такої події, тобто, що буде взяті саме ці деталі із загальної кількості $n$, треба перемножені комбінації поділити на комбінацію взятих до загальної кількості. Щоб визначити загальну кількість виробів, треба скласти усі $n$.

$$\sum n_i = 1 + 2 + 3 + 4 = 10$$

$$C_n^m = C_{10}^7 = \frac{10!}{7!\,(10-7)!} = \frac{7! \cdot 8 \cdot 9 \cdot 10}{7! \cdot 3!} = \frac{8 \cdot 9 \cdot 10}{1 \cdot 2 \cdot 3} = 4 \cdot 3 \cdot 10 = 120$$

Отримаємо ймовірність цієї події:

$$P = \frac{C_1^1 \cdot C_2^1 \cdot C_3^2 \cdot C_4^3}{C_{10}^7} = \frac{1 \cdot 2 \cdot 3 \cdot 4}{120} = \frac{24}{120} = 0{,}2$$

Крім розглянутих вище властивостей числа комбінацій без повторень (які також називають біноміальними коефіцієнтами, адже вони використовуються також у Біномі Ньютона):

$C_n^k = C_n^{n-k}$
$C_n^n = C_n^{n-n} = C_n^0 = 1,$

Важливими є наступні:
$$C_n^k + C_n^{k+1} = C_{n+1}^{k+1}$$

(використовується звичайно для обчислення не окремих значень, а таблиць). Якщо потрібно створити таблицю для конкретного $n$, то, щоб не робити повторень тих самих обчислень, можна рекомендувати рекурентною формулою, у якій кожен наступний елемент легко може бути обчислений із використанням значення попереднього:

$$C_n^k = C_n^{k-1} \frac{n-k+1}{k} = C_n^{k-1} \frac{n-(k-1)}{k}$$

Із врахуванням $C_n^0 = 1$, $C_n^1 = C_n^{1-1} \frac{n-1+1}{1} = 1 \frac{n}{1} = \frac{n}{1} = n,$
$C_n^2 = C_n^1 \frac{n-1}{2} = \frac{n}{1} \frac{n-1}{2},$
$C_n^3 = C_n^2 \frac{n-2}{3} = \frac{n}{1} \frac{n-1}{2} \frac{n-2}{3}$ і т.д.

Звичайно, можна таким чином обчислити лише першу половину коефіцієнтів для $k \leq n/2$, а для решти, скористатися симетрією: $C_n^{n-k} = C_n^k$.

**Комбінації з повтореннями.** Комбінацією з повтореннями з $n$ елементів по $k$ називається довільна підмножина з $k$ елементів з множини $M$ (елементи не обов'язково різні).

Число комбінацій з повтореннями з $n$ елементів по $k$ знаходиться за формулою:

$$\widetilde{C_n^k} = C_{n+k-1}^k = \frac{(n+k-1)!}{k!\,(n-1)!}$$

Тобто, $\widetilde{C_n^k}$ – це число способів, якими можна розкласти $k$ однакових предметів по $n$ урнах.

• **Приклад 1**. Якими та скількома способами можна покласти 3 однакових кулі у 2 урни?

Відповідь: Можна покласти в першу+другу урну 0+3,1+2,2+1 або 3+0 кульки, отже, 4 способа. За формулою, це число способів дорівнює $\widetilde{C_2^3} = \frac{(2+3-1)!}{3!(2-1)!} = 4$.

• **Приклад 2**. Якими та скількома способами можна покласти 2 однакових кулі у 3 урни?

Відповідь: за формулою, це число способів дорівнює $\widetilde{C_2^3} =$

$\frac{(3+2-1)!}{2!(3-1)!} = 6$. Це 0+0+2. 0+1+1, 0+2+0, 1+0+1, 1+1+0, 2+0+0.

• **Приклад 3.** В кондитерській є 6 різних сортів тістечок. Скільки є

способів купити 8 тістечок?

Відповідь:.Шукане число дорівнює

$$\widetilde{C_6^8} = C_{13}^8 = C_{13}^5 = \frac{13 \cdot 12 \cdot 11 \cdot 10 \cdot 9}{1 \cdot 2 \cdot 3 \cdot 4 \cdot 5} = 1287$$

• **Приклад 4**. Скількома способами можна покласти 15 однакових куль у 5 урн?

Відповідь: Це число способів дорівнює

$$\widetilde{C_{15}^5} = \frac{19 \cdot 18 \cdot 17 \cdot 16}{1 \cdot 2 \cdot 3 \cdot 4} = 3876.$$

**Задача 3**.

Серед $n$ лотерейних квитків, $k$ виграшних. Навмання взяли $m$ квитків. Визначити ймовірність того, що серед них $l$ виграшних.

Розглянемо для варіанту: $n = 10$, $l = 2$, $m = 4$, $k = 6$.

Відповідь:

Для вирішення задачі, використаємо формулу для числа комбінацій. Запишемо ще раз отримані дані і комбінації: $n$=10 загальна кількість квитків із них виграшних («першого ґатунку») $k$=6. Взято $m$=4 і із них виграшних $l$=2. Тобто, число способів взяти 4 квитка із 10 буде комбінація $C_{10}^4$, але, крім 4 виграшних з 6 ($C_6^2$ комбінацій), є два невиграшних з 4 ($C_4^2$ комбінацій). Загальна кількість сприятливих комбінацій $C_6^2 C_4^2$.

Ймовірність дорівнює відношенню кількості сприятливих комбінацій до загальної.

Загальна формула для ймовірності того, що буде взято $l$ об'єктів з $k$:

$$p_l = \frac{C_k^l \cdot C_{n-k}^{m-l}}{C_n^m}.$$

Такий статистичний розподіл називають «гіпергеометричним» розподілом. По загальній формулі, число комбінацій

$$C_n^k = \frac{n!}{k!(n-k)!}$$

при $C_n^k$ ($0 \le k \le \min(k, m)$), інакше $C_n^k = 0$.

Обчислимо комбінації $C_k^l, C_{n-k}^{m-l}$ та $C_n^m$, перевіряючи параметри на виконання природніх умов (кількість успіхів не може бути від'ємною та перевищувати а ні число спроб, ні число об'єктів у ресурсі) :

$$C_4^2 = \frac{4!}{2!\,(4-2)!} = \frac{1 \cdot 2 \cdot 3 \cdot 4}{2 \cdot 2} = 6$$
$$C_6^2 = \frac{6!}{2!\,(6-2)!} = \frac{4! \cdot 5 \cdot 6}{2! \cdot 4!} = 15$$
$$C_{10}^4 = \frac{10!}{4!\,(10-4)!} = \frac{6! \cdot 7 \cdot 8 \cdot 9 \cdot 10}{4! \cdot 6!} = 210$$

Ймовірність дорівнює:
$$p_2 = \frac{C_6^2 \cdot C_4^2}{C_{10}^4} = \frac{15 \cdot 6}{210} = \frac{3}{7} = 0{,}429$$

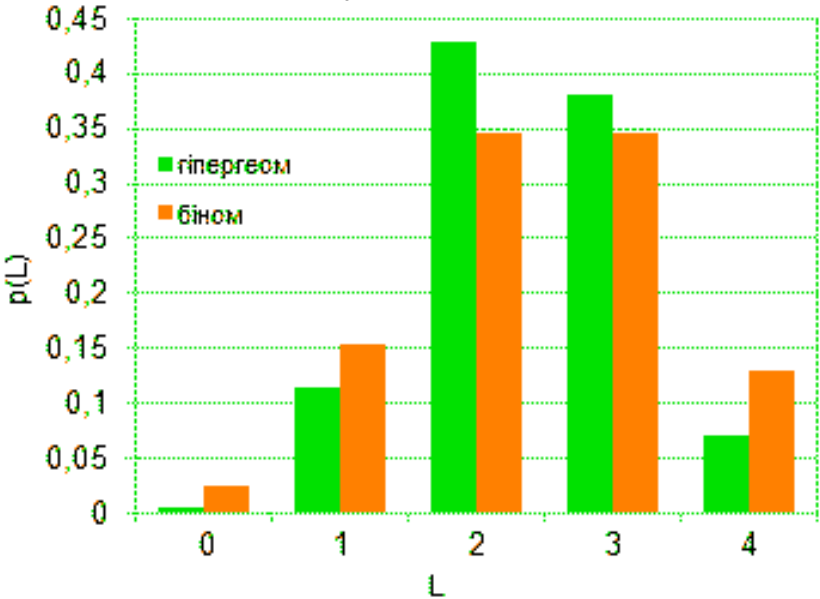

Для ілюстрації, на рис. показана гістограма розподілу для усіх можливих значень $l = 0,1,2,3,4$ для інших початкових даних задачі (стовпчики зліва), а не лише для $l = 2$. Це є «вибірка без повернення», тобто ймовірність наступного результату залежить від того, який об'єкт вибрано. У випадку «вибірки з поверненням», ймовірність чергового результату постійна, й дорівнює $p = k/n$. Ймовірності для кількості успіхів в цьому випадку можуть бути обчислені за законом біноміального розподілу й показані правими стовпчиками. При пропорційному збільшенні кількості елементів у ресурсі, гіпергеометричний розподіл буде асимптотично наближатись до біноміального $p_l = C_m^l \, p^l \, (1-p)^{m-l}$.

В основі розв'язування багатьох комбінаторних задач лежать два основних правила — правило суми і правило добутку.

**Правило суми**: Якщо елемент $A$ можна вибрати $m$ способами, а інший елемент $B$ — $n$ способами, то $A$ або $B$ можна вибрати $(m + n)$ способами.

Наприклад: Якщо на тарілці лежить 5 груш і 4 яблука, то вибрати один фрукт (тобто, грушу або яблуко) можна 9 способами (5 + 4 = 9).

**Правило добутку**: Якщо елемент A можна вибрати $m$ способами, а після цього елемент B — $n$ способами, то A та B можна вибрати $(m \cdot n)$ способами.

Це твердження означає, що, оскільки для кожного з $m$ елементів $A$, можна взяти в пару будь-який з $n$ елементів $B$, то кількість пар дорівнює добутку $m \cdot n$.

Наприклад: Якщо в кіоску продають ручки 5 видів і зошити 4 видів, то вибрати набір з ручки і зошита (тобто, пару — ручка і зошит) можна 5·4 = 20 способами (оскільки до кожної з 5 ручок, можна взяти будь-який із 4 зошитів).

Отже, якщо доводиться вибирати або перший елемент, або другий, або третій і т. д. елемент, кількості способів вибору кожного елементу додають, а коли доводиться вибирати набір, у який входить і перший, і другий, і третій, і т.д. елементи, кількості способів вибору кожного елементу перемножують.

Розглянуті вище правила, аксіоми та властивості ймовірностей лежать в основі основних теорем теорії ймовірностей.

# 3. ОСНОВНІ ТЕОРЕМИ ТА ФОРМУЛИ ТЕОРІЇ ЙМОВІРНОСТІ

## 3.1. Теореми додавання та добутку ймовірностей.

**Теорема 1.** (додавання ймовірностей несумісних подій).

Ймовірність об'єднання попарно несумісних подій дорівнює сумі ймовірностей цих подій. Так, для двох несумісних подій:
$$P(A \cup B) = P(A + B) = P(A) + P(B) \qquad (1)$$

Висновок із теореми. Якщо події $A_1, A_2, \ldots, A_n$ утворюють повну групу несумісних подій, то сума їх ймовірностей дорівнює 1:
$$\sum_{i=1}^{n} P(A_i) = 1$$

• **Приклад.** В ящику 10 червоних і 6 синіх ґудзиків. Навмання виймають два ґудзики. Яка ймовірність того, що ґудзики будуть одного кольору?

Відповідь:. Введемо позначення: подія А – ґудзики одного кольору, подія В – ґудзики червоні, подія С– ґудзики сині. Очевидно, $A = B \cup C$ і події В і С несумісні. За формулою (1) P(A)=P(B)+P(C). Знайдемо P(B) і P(C). Число способів взяти 2 ґудзики з 16 дорівнює $C_{16}^2$. Число наслідків, сприятливих для події В, дорівнює $C_{10}^2$, а для події С це $C_6^2$.

Одержимо:
$$P(B) = \frac{C_{10}^2}{C_{16}^2} = \frac{10 \cdot 9}{1 \cdot 2} \cdot \frac{1 \cdot 2}{16 \cdot 15} = \frac{3}{8}$$
$$P(C) = \frac{C_6^2}{C_{16}^2} = \frac{6 \cdot 5}{1 \cdot 2} \cdot \frac{1 \cdot 2}{16 \cdot 15} = \frac{1}{8}$$

Отже:
$$P(A) = \frac{3}{8} + \frac{1}{8} = \frac{1}{2}$$

**Теорема 2.** (додавання будь-яких двох подій).

Ймовірність появи об'єднання будь яких двох подій дорівнює сумі їх ймовірностей без ймовірності перетину (сумісної появи) $P(A \cap B) = P(AB)$:
$$P(A \cup B) = P(A + B) = P(A) + P(B) - P(AB) \qquad (2)$$

Висновок із теореми. Ймовірність суми трьох сумісних подій

$$P(A \cup B \cup C) = P(A) + P(B) + P(C) -$$
$$-P(A \cap B) - P(A \cap C) - P(B \cap C) + P(A \cap B \cap C).$$

• **Приклад.** Знайти ймовірність того, що навмання вибране двозначне число є кратним 2, або 5, або тому й іншому одночасно.

Відповідь: Введемо позначення: подія $A$ – вибране число кратне 2, або 5, або тому й іншому одночасно, подія $B$ – число кратне 2, подія $C$ – число кратне 5. Очевидно, що $A = B \cup C$. За формулою (2): P(A)=P(B)+P(C) –P(BC). Всього є $n$=90 двозначних чисел (10,12,…98), з них, 45 кратних 2. Кількість двозначних чисел, кратних 5, за правилом добутку дорівнює n=9·2 (маємо 9 десятків і в кожному десятку два числа кратні 5. тобто, (10,15,…,95), з них 9 чисел, кратних і 2 і 5 (тому що кожне двозначне число, що закінчується на нуль кратне 2 та 5, а таких чисел 9 (10,20,…,90)). Маємо P(B)=$\frac{45}{90} = \frac{1}{2}$, P(C)=$\frac{18}{90} = \frac{1}{5}$, P(BC)=$\frac{9}{90} = \frac{1}{10}$. Отже,

P(A)= $\frac{1}{2} + \frac{1}{5} - \frac{1}{10} = \frac{5+2-1}{10} = \frac{6}{10} = \frac{3}{5} = 0{,}6$

**Умовні ймовірності та незалежні події.**

Умовною ймовірністю події А за умови, що відбулася подія В, називається величина $P(A/B) = P(A|B)$:
$$P(A/B) = \frac{P(AB)}{P(B)}, P(B) > 0. \tag{3}$$

Подія А називається незалежною від події В, коли ймовірність події А не залежить від того, відбулася подія В, чи не відбулася.

Подія А називається залежною від події В, коли ймовірність події А змінюється залежно від того, відбулася подія В, чи не відбулася.

**Теорема добутку ймовірностей**.

Ймовірність сумісної появи двох подій дорівнює ймовірності однієї з них, помножену на умовну ймовірність іншої:
$$P(A \times B) = P(A) \times P(B/A) = P(B) \times P(A/B) \tag{4}$$

Висновок: Якщо подія А не залежить від події В, то і подія В не залежить від події А. У разі $n$ незалежних подій:
$$P\left(\prod_{i=1}^{n} A_i\right) = \prod_{i=1}^{n} P(A_i).$$

**Теорема 3.** (добутку ймовірностей для незалежних подій).

Ймовірність сумісної появи двох незалежних подій дорівнює добутку їх ймовірностей:
$$P(A \times B) = P(A) \times P(B).$$

• **Приклад 1**. Кинуто три монети. Визначити, залежні чи незалежні події:

А – випав герб на першій монеті, В – випала хоча б одна цифра.
Відповідь:

Простір елементарних подій: $\Omega =$ {ггг, ггц, гцг, цгг, гцц, цгц, ццг, ццц},

подія А={ггг, ггц, гцг, гцц},

подія В = {ггц, гцг, цгг, гцц, цгц, ццг, ццц},

подія АВ ={ггц, гцг, гцц}.

Очевидно, що:
$$P(A)=\frac{1}{2}, \; P(B)=\frac{7}{8}, \; P(AB)=\frac{3}{8},$$
але, якщо порахувати по формулі, то
$$P(A) \cdot P(B) = \frac{1}{2} \cdot \frac{7}{8} = \frac{7}{16},$$
тобто, формула для незалежних подій не виконується, що свідчить о залежності подій. Події А і В є залежні.

• **Приклад 2.** В урні M=12 занумерованих куль з номерами від 1 до М. Кулі виймаються по одній без повернення. Розглянемо такі події: А – номери куль в порядку вилучення утворюють послідовність 1,2,…,M; В – хоч один раз збігається номер кулі і порядковий номер вилучення; С – нема жодного збігу номера кулі та порядкового номера вилучення. Визначити ймовірність подій А, В, С.

Відповідь:

Ця задача про вибірку без повернення, оскільки кульки мають власні номери (позначки) і не повертаються. За умовами задачі, «сприятливим» є лише одне сполучення $1,2,...,m$. Загальна кількість перестановок $m$ елементів з $m$ дорівнює $N_m = m!$

Ймовірність того, що першою буде кулька із першим номером, дорівнює $P(A_1) = \frac{1}{m}$, для другої, залишається на одну можливість менше, отже, $P(A_2|A_1) = \frac{1}{m-1}$, (ймовірність умовна, бо потрібно, щоб попередні кульки були із «правильними» номерами). Згадаємо, що «умовною ймовірністю події А за умови, що відбулася подія В, називається величина:

$$P(A|B) = P(A \vee B) = \frac{P(AB)}{P(B)}, P(B) > 0$$

В нашому випадку замість В маємо події $A_1$, і т.д. Таким чином, для третьої кульки маємо записати ймовірність як $P(A_3|A_2 A_1) = \frac{1}{m-2}$, і т.д. Ймовірність складної події (усі кульки витягнуті із «правильними» номерами), дорівнює

$$p_m = P(A_m A_{m-1} \ldots A_2 A_1) = P(A_m|A_{m-1} \ldots A_2 A_1) P(A_{m-1}| \ldots A_2 A_1) =$$
$$= \frac{1}{1} \frac{1}{2} \ldots \frac{1}{m-1} \cdot \frac{1}{m} = \frac{1}{m!}$$

Підставимо $m = 12$ і отримаємо:

$$p_m = \frac{1}{m!} = \frac{1}{12!} = \frac{1}{1 \cdot 2 \cdot 3 \cdot \ldots \cdot 12} = \frac{1}{479001600} = 2 \cdot 10^{-9}$$

*Примітка: Якщо якийсь елемент повторюється, наприклад, (1,2,2), то ймовірність дорівнює нулю, бо така послідовність неможлива у вибірці «без повернення».*

*Можливий варіант задачі, коли потрібно, щоб збігалися не усі $m$ чисел, а лише перші $k$. У цій задачі, нас «не цікавлять» послідовності останніх $(m-k)$ чисел, а таких буде $(m-k)!$ замість єдиної. Отже, маємо:*

$$p_{k,m} = \frac{(m-k)!}{m!} = \frac{1}{A_m^k}.$$

Ймовірності події В і С відповідають ймовірностям протилежних подій («є хоча б один збіг номеру кульки із порядковим номером у послідовності» та «є нуль (тобто, немає жодного) збігів номеру кульки із порядковим номером у послідовності»). Отже, можна скористатись властивістю, що сума двох протилежних подій дорівнює одиниці, і спершу обчислити щось одне. Візьмемо, наприклад, подію С «жодного збігу». Обчислимо кількість «сприятливих» послідовностей для перших значень $m$, і позначимо її $M_m$:

$M_1 = 0$, адже у послідовності (1) є збіг;

$M_2 = 1$, адже без збігу лише послідовність (2,1);

$M_3 = 2$, адже без збігу лише послідовності (2,3,1), (3,1,2)

Можна скористатись рекурентною формулою, обчислюючи

послідовно:
$M_1 = 0; M_2 = 2M_1 + 1; M_3 = 3M_2 - 1; …; M_m = mM_{m-1} + (-1)^m;$

Ця формула нагадує вираз для факторіала, де $m! = N_m = mN_{m-1}$, але додається $+1$, якщо $m-$ парне число, або $-1$, якщо непарне. Пам'ятаємо, що $N_m = m!$ загальна кількість перестановок $m$ елементів з $m$.

Для числа послідовностей «хоча б один збіг», $L_m = N_m - M_m$,
$L_1 = 1; L_2 = 2L_1 - 1; L_3 = 3L_2 + 1; …; L_m = mL_{m-1} - (-1)^m;$

Відповідні ймовірності:
$$P(C) = p_{mc} = \frac{M_m}{N_m} = 1 - \frac{1}{1!} + \frac{1}{2!} - \frac{1}{3!} + \cdots + \frac{(-1)^m}{m!}.$$

Розрахунок за формулою: $1 - \frac{1}{1!} + \frac{1}{2!} - \frac{1}{3!} + \cdots + \frac{(-1)^{12}}{12!}$

надає більш точний результат, але можна скористатися і формулою де в чисельнику кількість «сприятливих» послідовностей, а у знаменнику загальна кількість перестановок:

$$P(C) = p_{mc} = \frac{M_m}{N_m} = \frac{mM_{m-1} + (-1)^m}{mN_{m-1}} = \frac{12M_{11} + (-1)^{12}}{12 \cdot N_{11}}$$
$$= \frac{176214841}{12 \cdot 11!} =$$
$$= \frac{176214841}{479001600} = 0{,}368$$

$$P(B) = p_{mв} = \frac{L_m}{N_m} = 1 - p_{mc} = \frac{+1}{1!} - \frac{1}{2!} + \frac{1}{3!} + \cdots - \frac{(-1)^m}{m!},$$

$P(B) = p_{mв} = 1 - p_{mc} = 1 - 0{,}368 = 0{,}632$

*Слід зауважити, що сума для $P(C)$ є обмеженою сумою ряду Тейлора (у модифікації Мак-Лорена), тому ймовірності асимптотично наближуються відповідно до $(1 - e^{-1})$ та $e^{-1} = \frac{1}{e} = \frac{1}{2{.}17281728…} = 0{.}36787944…$*

Якщо розглядати вибірку із поверненням (наприклад, телефонні номери), то $\tilde{P}(A) = m^{-m}$, тобто $1{.}12 \cdot 10^{-13}$, що значно менше значення для вибірки без повернення $\frac{1}{12!} = 2{.}09 \cdot 10^{-9}$. Втім, $\tilde{P}(C) = (1 - \frac{1}{m})^m$, тобто $0{.}351995628$ – відносно значно ближче до $P(C)$.

**Теорема 4.** (правило знаходження ймовірності протилежної події).

Ймовірність настання принаймні однієї з подій $A_1, A_2, \ldots, A_n$, незалежних в сукупності, знаходиться за формулою:
$$P(A) = 1 - q_1 q_2 \ldots q_n \quad (5)$$

де $q_n$ — ймовірність протилежної події, тобто, $P(\overline{A_n})$.

Наслідок: Якщо події $A_1, A_2, \ldots, A_n$ мають однакову ймовірність $p$, то ймовірність настання принаймні однієї з цих подій дорівнює:
$$P(A) = 1 - q^n \quad (6)$$

• **Приклад 1.** У кожному з трьох ящиків лежить по 10 деталей; у першому ящику 2 деталі браковані, у другому – 3, у третьому – 1. З кожного ящика беруть по одній деталі. Знайти ймовірність того, що серед них є принаймні одна стандартна.

Відповідь:

Позначимо події: A – серед трьох деталей принаймні одна стандартна, $A_i (i = 1,2,3)$ – деталь, взята з $i$-го ящика, стандартна. Маємо:

у першому ящику 8 деталей не браковані, тому $P(A_1) = \frac{8}{10}$,

у другому ящику 7 деталей не браковані, тому $P(A_2) = \frac{7}{10}$,

у третьому ящику 9 деталей не браковані, тому $P(A_3) = \frac{9}{10}$,

звідси маємо, що:
$$q_1 = P(\overline{A_1}) = \frac{2}{10}, q_2 = P(\overline{A_2}) = \frac{3}{10}, q_3 = P(\overline{A_3}) = \frac{1}{10}.$$

За формулою (5) знаходимо ймовірність того, що серед трьох деталей принаймні одна стандартна:
$$P(A) = 1 - 0{,}2 \cdot 0{,}3 \cdot 0{,}1 = 1 - 0{,}006 = 0{,}994$$

• **Приклад 2.** Ймовірність хоча б одного влучення в ціль при чотирьох пострілах дорівнює 0,9984. Знайти ймовірність влучення в ціль при одному пострілі.

Відповідь:

Нехай подія A – хоча б одне влучення в ціль при чотирьох пострілах,

p – ймовірність влучення при одному пострілі, тоді:

$q = 1 - p$ – ймовірність промаху при одному пострілі і $P(A) = 1 - q^4$.

За умовою, що $1 - q^4 = 0{,}9984$, звідки $q^4 = 0{,}0016$, маємо, що:
$$q = 0{,}2, \quad p = 1 - 0{,}2 = 0{,}8.$$

• **Приклад 3.** Ймовірність того, що у ціль влучив при одному пострілі перший стрілок $p_1 = 0{,}61$, другий – $p_2 = 0{,}55$. Перший здійснив $n_1 = 2$, а другий – $n_2 = 3$ пострілів. Визначити ймовірність того, що в ціль не влучили.

Відповідь.:

Ця задача відноситься до типу задач, в яких одно і теж випробування повторюється багато разів. Як результат з'являється, або не з'являється якась подія. В задачі головне загальне число появ якоїсь події в серії випробувань.

Позначимо через А – подію, що перший стрілок влучив, а через В, що другий стрілок влучив. За умовами задачі, маємо ймовірність попадання в ціль при одному пострілі для першого стрілка $p_1 = 0{,}61$ та $p_2 = 0{,}55$ для другого. Позначимо ймовірності протилежних подій, тобто, того, що не влучив перший, $q_1 = 1-p_1 = 1 - 0{,}61 = 0{,}39$, не влучив другий, $q_2 = 1-p_2 = 1 - 0{,}55 = 0{,}45$. Ймовірність того, що не влучив перший при $n_1$ пострілах, дорівнює:
$$P(\overline{A}) = q_1 q_1 \ldots q_1 = q_1^{n_1}.$$

Ймовірність того, що не влучив другий при $n_2$ пострілах, дорівнює:
$$P(\overline{B}) = q_2 q_2 \ldots q_2 = q_2^{n_2}.$$

Оскільки кількість та результати пострілів незалежні, то ймовірність того, що обидва стрілка не влучили в ціль при $n_1 = 2$ пострілах і $n_2 = 3$ пострілах:
$$P(\overline{AB}) = P(\overline{A})P(\overline{B}) = q_1^{n_1} q_2^{n_2} = 0{,}39^2 \cdot 0{,}45^3 = 0{,}15 \cdot 0{,}09 = 0{,}0135.$$

У випадку ненульової кількості влучень, потрібно розглядати складові можливості. Наприклад, 2 влучання можуть бути розподілені між стрільцями, як 0+2, 1+1, 2+0, отже, потрібно використовувати формулу повної ймовірності.

**3.2. Формула повної ймовірності.**

Формула повної ймовірності – висновки теорем додавання і множення ймовірностей. Нехай деяка подія А може відбутися тільки разом з однією з подій $H_1, H_2, \ldots, H_n$. Події $H_1, H_2, \ldots, H_n$ складають повну групу несумісних подій і називаються гіпотезами.

**Теорема 5.**

Ймовірність події A, яка може відбутись лише за умови появи однієї з несумісних подій $H_1, H_2, \ldots, H_n$, що утворюють повну групу, дорівнює сумі добутків ймовірностей кожної з цих подій на відповідну умовну ймовірність події A:

$$P(A) = \sum_{i=1}^{n} P(H_i) \cdot P(A / H_i). \quad (7)$$

Ця формула (7) називається **формулою повної ймовірності**.

• **Приклад 1**. У трьох урнах лежать білі і чорні кулі. У першій – 3 білі і 1 чорна, у другій – 6 білих і 4 чорних, у третій – 9 білих і 1 чорна. З навмання взятої урни виймають одну кулю. Знайти ймовірність того, що вона біла.

Відповідь: Введемо позначення: події $H_i (i = 1,2,3)$ – вибрана $i$-та урна, подія A – взята куля біла.

Очевидно, що:

$$P(H_1) = P(H_2) = P(H_3) = \frac{1}{3},$$

Іноді, замість $H_1$ (Hypothesis, гіпотеза), використовують $B_i$ (літера абетки, наступна після $A$), наприклад, $B$, як у формулі умовної ймовірності:

$$P(A|B) = \frac{P(AB)}{P(B)}$$

По умові задачі треба знайти P(A), яке можна порахувати по формулі:

P(A) = P(H$_1$) P(A|H$_1$) + P(H$_2$) P(A|H$_2$) + P(H$_3$) P(A|H$_3$)

Треба ще знайти $P(A|H_1), P(A|H_2), P(A|H_3)$, тобто, ймовірність взяття білої кулі із першої урни, із другої та третьої. Для цього, ділимо кількість білих кульок на загальну кількість:

$$P(A|H_1) = \frac{3}{4}$$
$$P(A|H_2) = \frac{6}{10}$$
$$P(A|H_2) = \frac{9}{10}.$$

Підставимо знайдені значення у формулу і знайдемо P(A)

$$P(A) = \frac{1}{3} \cdot \frac{3}{4} + \frac{1}{3} \cdot \frac{6}{10} + \frac{1}{3} \cdot \frac{9}{10} = \frac{1}{4} + \frac{1}{5} + \frac{3}{10} = \frac{5+4+6}{20} = \frac{15}{20} = \frac{3}{4}$$

Часто виникає потреба визначити умовну ймовірність, напр.

$$P(H_1|A) = \frac{P(A \cdot H_1)}{P(A)} = \frac{P(A|H_1) \cdot P(H_1)}{P(A)} = \frac{3/4 \cdot 1/3}{3/4} = \frac{1}{3}$$

• **Приклад 2.**

В двох партіях $k_1 = 71\%$ і $k_2 = 47\%$. якісних виробів. Навмання вибирають по одному виробу з кожної партії. Яка ймовірність виявити серед них:

а) хоч один неякісний;

б) два неякісних;

в) один якісний і один неякісний?

Відповідь:

Ця задача на формулу «повної ймовірності» у різних умовах. Позначимо через А подію потрапляння якісного виробу із першої партії і через В подію потрапляння якісного виробу із другої партії, тоді ймовірності цих подій:

$P(A) = \frac{k_1\%}{100\%}, P(B) = \frac{k_2\%}{100\%}$, підставимо наші значення і отримаємо

$P(A) = 0.71$ і $P(B) = 0{,}47$. Ймовірність протилежних подій:

$P(\overline{A}) = 1 - P(A) = 1 - 0{,}71 = 0{,}29$ та $P(\overline{B}) = 1 - P(B) = 1 - 0{,}47 = 0{,}53$.

Тоді можливі комбінації незалежних подій

$P(\overline{A} \cdot \overline{B}) = P(\overline{A}) \cdot P(\overline{B}) = 0{,}29 \cdot 0{,}53 = 0{,}1537$

$P(\overline{A} \cdot B) = P(\overline{A}) \cdot P(B) = 0{,}29 \cdot 0{,}47 = 0{,}1363$

$P(A \cdot \overline{B}) = P(A) \cdot P(\overline{B}) = 0{,}71 \cdot 0{,}53 = 0{,}3763$

$P(A \cdot B) = P(A) \cdot P(B) = 0{,}71 \cdot 0{,}47 = 0{,}3337$

Питання: «Яка ймовірність виявити серед них: хоч один неякісний?» має на увазі, що це може бути в тому числі і два неякісних вироби, тому можливі три комбінації подій: $P(\overline{A} \cdot \overline{B}); P(\overline{A} \cdot B); P(A \cdot \overline{B})$. Повинна відбутися будь-яка комбінація із трьох можливих, тому, щоб визначити ймовірність цієї події маємо знайти об'єднання комбінацій, тобто, їх суму.

Таким чином, відповідь на а) Р=0,1537+0,1363+0,3763=0,6663

б) Яка ймовірність виявити серед них два неякісних? – для цієї події підходить тільки одна комбінація, тому відповідь:

Р=$P(\overline{A} \cdot \overline{B}) = P(\overline{A}) \cdot P(\overline{B}) = 0{,}29 \cdot 0{,}53 = 0{,}1537$

в) Яка ймовірність виявити серед один якісний і один неякісний? – для цієї події підходять дві комбінації $P(\overline{A} \cdot B)$ та $P(A \cdot$

$\overline{B}$), тому:

P=0,1363+0,3763=0,5126

• **Приклад 3.** З 1000 ламп $n_i$ належать $i$–ій партії $i$=1,2,3,
$$\sum_{i=1}^{3} n_i = 1000.$$

В першій ($n_1 = 100$) партії 6%, у другій ($n_2 = 250$) – 5%, у третій - 4% неякісних ламп. Навмання вибирається одна лампа. Визначити ймовірність того, що вона неякісна.

Відповідь: По умовам задачі маємо 100 ламп з першої партії та 250 ламп з другої партії при загальній кількості 1000 ламп. Знайдемо кількість ламп з третьої партії $n_3 = 1000 - (100 + 250) = 650$.

Задача вирішується за формулою повної ймовірності, згадаємо теорему і формулу: «Ймовірність події А, яка може відбутись лише за умови появи однієї з несумісних подій $H_1, H_2, ..., H_n$, що утворюють повну групу, дорівнює сумі добутків ймовірностей кожної з цих подій на відповідну умовну ймовірність події А»:

$$P(A) = \sum_{i=1}^{n} P(H_i) \cdot P(A / H_i).$$

Введемо позначення: подія $A$– взята лампа неякісна, події $H_i = (H_1, H_2, H_3)$ – що лампа належить до однієї із трьох партій. Визначимо умовні ймовірності того, що взята неякісна лампа належить до відповідної групи за формулою:

$$P(A|B) = \frac{P(AB)}{P(B)}$$

отримаємо $P(A|H_1) = 6\% = 0{,}06, P(A|H_2) = 5\% = 0{,}05, P(A|H_3) = 4\% = 0{,}04$.

Формула повної ймовірності:
$P(A) = P(H_1) \cdot P(A|H_1) + P(H_2) \cdot P(A|H_2) + P(H_3) \cdot P(A|H_3)$

Визначимо $P(H_1), P(H_2), P(H_3)$, тобто, ймовірності того, що вибрана лампочка належить до першої, другої та третьої партії відповідно за формулою:

$$P(H_i) = \frac{n_i}{n}$$

Значення $n_i$ надані в умовах задачі: $n_1 = 100, n_2 = 250, n_3 = 650$ при

$n = 1000$.

$P(H_1) = \dfrac{100}{1000} = 0{,}1; P(H_2) = \dfrac{250}{1000} = 0{,}25; P(H_3) = \dfrac{650}{1000} = 0{,}65;$

Підставимо отримані значення в формулу повної ймовірності і отримаємо:
$P(A) = 0,1 \cdot 0,06 + 0,25 \cdot 0,05 + 0,65 \cdot 0,04 = 0,006 + 0,0125 + 0,026 = 0,0445$

Враховуючи те, що неякісна лампа може належати до однієї із трьох партій, ймовірність, що лампочка неякісна для цієї задачі дорівнюється 0,0445.

• **Приклад 4.**

В першій урні $N_1 = 4$ білих та $M_1 = 1$ чорних куль, в другій $N_2 = 2$ білих та $M_2 = 5$ чорних. З першій в другу переклали $k = 3$ куль, потім з другої вилучили одну кулю. Визначити ймовірність того, що вибрана з другої урни куля – біла.

Відповідь: Ця задача також на формулу повної ймовірності. Рішення задачі треба почати з розгляду кульок, що були перекладені з першої урни до другої, тобто, показник k=3. Умовами задачі не визначено скільки куль білого та чорного кольору було серед перекладених до другої урни тому потрібно розглянути всі можливі комбінації щодо наявності білої кульки серед перекладених. За умовами цієї задачі перекладено три кульки. Поглянемо на кількість білих куль у першій урні - $N_1 = 4$, тобто, можливий варіант, що всі 3 перекладені кульки білі. Наступний варіант – 2 білі і 1 чорна, далі 1 біла і 2 чорні і останній варіант – 0 білих і 3 чорних. Таким чином, кількість білих кульок $i$ може приймати значення у максимальному інтервалі $0 \leq i \leq k$, але, не може бути витягнуто білих більше, ніж їх є у ресурсі.

Введемо позначення: подія $A$ – витягнута біла куля і P($A$) – ймовірність того, що з другої урни буде витягнуто білу кульку, подія $B_i$ – полягає в тому, що перекладено $i-$ кульок білого кольору із $N_1 = 4$, P($B_i$) – ймовірність події $B_i$. Відповідно, маємо $(k-i)$ кульок чорного кольору з $M_1$ і тут, для багатьох варіантів, діапазон може бути звужений. Так, наприклад, в нашому випадку у першій урні була тільки одна кулька чорного кольору, тому варіанти коли взято 1 біла і 2 чорні та 0 білих і 3 чорних треба зняти з розгляду при $.M_1 = 1$.

Хоч, при розгляді всіх варіантів без винятку, ми отримаємо 0 в тих випадках, що є неможливі, в кінцевому розрахунку маємо правильну відповідь.

Умовна ймовірність витягнути білу кулю, якщо у другий кошик переклали $i$ білих куль та $k-i$ чорних, дорівнює відношенню кількості білих куль на загальну кількість білих + чорних:
$$P(A|B_i) = \frac{N_2 + i}{N_2 + M_2 + k}$$

Ймовірність того, що випадково буде взято при перекладенні $i$ білих куль, визначається гіпергеометричним розподілом (який ми розглядали у задачі 3). Із позначеннями задачі 10, можна записати
$$P(B_i) = \frac{C_{N_1}^i C_{M_1}^{k-i}}{C_{N_1+M_1}^k}$$

Таким чином, повна ймовірність події $A$:
$$P(A) = \sum_{i=0}^{k} P(A|B_i) \cdot P(B_i) = \sum_{i=0}^{k} \frac{N_2 + i}{N_2 + M_2 + k} \cdot \frac{C_{N_1}^i C_{M_1}^{k-i}}{C_{N_1+M_1}^k}$$

Розташуємо всі дані в вигляді таблиці:

| $i$ – число білих серед взятих | $k-i$ – число чорних серед взятих | $N_2 + i$ – число білих у 2 урні після перекладення | $C_{N_1}^i$ – число комбінацій, якими можна взяти $i$ білих з $N_1$ | $C_{M_1}^{k-i}$ – число комбінацій, якими можна взяти $k-i$ чорних з $M_1$ | $(N_2+i)C_{N_1}^i C_{M_1}^{k-i}$ – чисельник остаточної відповіді |
|---|---|---|---|---|---|
| 0 | 3 | 2 | $C_4^0 = 1$ | $C_1^3 = 0$ | 0 |
| 1 | 2 | 2+1=3 | $C_4^1 = 4$ | $C_1^2 = 0$ | 0 |
| 2 | 1 | 2+2=4 | $C_4^2 = 6$ | $C_1^1 = 1$ | 24 |
| 3 | 0 | 2+3=5 | $C_4^3 = 4$ | $C_1^0 = 1$ | 20 |

загальна кількість куль у другій урні після перекладення:
$N_2 + M_2 + K = 2 + 5 + 3 = 10$

Комбінацію для загальної кількості кульок у першій урні (як і інші) розобчислюємо за формулою:
$$C_n^k = \frac{n!}{k!(n-k)!}$$

Маємо:
$$C_{N_1+M_1}^K = C_{4+1}^3 = \frac{5!}{3!(5-3)!} = \frac{4 \cdot 5}{2} = 10$$

Підставимо отримані дані у формулу повної ймовірності:

$$P(A) = \sum_{i=0}^{k} \frac{N_2 + i}{N_2 + M_2 + K} \cdot \frac{C_{N_1}^i C_{M_1}^{k-i}}{C_{N_1+M_1}^k} = \frac{4}{10} \cdot \frac{6 \cdot 1}{10} + \frac{5}{10} \cdot \frac{4 \cdot 1}{10}$$
$$= 0{,}24 + 0{,}2 = 0{,}44$$

Це повна ймовірність того, що з другого кошика буде витягнута біла кулька (відповідно, 1-0.44=0.56 – ймовірність того, що буде витягнута чорна кулька).

### 3.3. Формули Байєса.

Нехай тепер деяка подія A відбулася разом з однією із подій $H_1, H_2, \ldots, H_n$, що складають повну групу несумісних подій. Визначимо ймовірність гіпотез $H_1, H_2, \ldots, H_n$ за умови, що подія A вже відбулася, тобто, обчислимо умовні ймовірності $P(H_1|A), \ldots, P(H_n|A)$.

За теоремою множення $P(A \cdot H_1) = P(A) \cdot P(H_1|A)$, звідки
$$P(H_1|A) = \frac{P(A \cdot H_1)}{P(A)} = \frac{P(H_1) \cdot P(A|H_1)}{P(A)}$$

Ймовірність $P(A)$ обчислюється по формулі (7). Тоді:
$$P(H_1|A) = \frac{P(H_1) \cdot P(A|H_1)}{\sum_{i=1}^{n} P(H_i) P(A|H_i)}$$

Аналогічно можна обчислити умовну ймовірність $P(H_2|A), \ldots, P(H_n|A)$. Тому загальну формулу треба записати у вигляді:
$$P(H_k|A) = \frac{P(H_k) \cdot P(A|H_k)}{\sum_{i=1}^{n} P(H_i) P(A|H_i)}, \quad k = 1, 2, \ldots, n \qquad (8)$$

Формули Байєса (8) дозволяють перерахувати ймовірність гіпотез після того, як подія A відбулася в результаті випробування.

• **Приклад 1.** На фабриці виготовляють гвинти. Перша машина виготовляє 25%, друга – 35%, третя – 40% усіх гвинтів. Частка браку відповідно 5%, 4% і 2%. Випадково вибраний гвинт виявився бракованим. Яка ймовірність того, що його зроблено першою, другою, третьою машинами?

Відповідь: Позначимо події: $H_i (i = 1,2,3)$ – вибраний гвинт виготовлений і-ю машиною, A – вибраний гвинт бракований. Ймовірність того, що гвинт зроблений першою машиною, тобто, враховуємо до загальної кількості – 100%. $P(H_1) = \frac{25\%}{100\%} = 0{,}25$,

другою - $P(H_2) = \frac{35\%}{100\%} = 0{,}35$

і третьою - $P(H_3) = \frac{40\%}{100\%} = 0{,}2$

Ймовірність події А - того, що гвинт бракований:

для першої машини - $P(A|H_1) = \frac{5\%}{100\%} = 0{,}05$ ,

для другої - $P(A|H_2) = \frac{4\%}{100\%} = 0{,}04$

та третьої - $P(A|H_3) = \frac{2\%}{100\%} = 0{,}02$.

Застосуємо формулу повної ймовірності:
$$P(A) = \sum_{i=1}^{n} P(H_i) \cdot P(A/H_i).$$

Отримаємо:
$$P(A) = 0{,}25 \cdot 0{,}05 + 0{,}35 \cdot 0{,}04 + 0{,}40 \cdot 0{,}02 = 0{,}0345.$$

Тепер, по формулі (8):
$$P(H_k|A) = \frac{P(H_k) \cdot P(A|H_k)}{\sum_{i=1}^{n} P(H_i) P(A|H_i)}, \quad k = 1,2,\ldots,n$$

Сума у знаменнику це $P(A)$, яке ми вже порахували. Підставимо пораховане:

$$P(H_1/A) = \frac{P(H_1)P(A/H_1)}{P(A)} = \frac{0{,}25 \cdot 0{,}05}{0{,}0345} = 0{,}3623$$

$$P(H_2/A) = \frac{P(H_2)P(A/H_2)}{P(A)} = \frac{0{,}35 \cdot 0{,}04}{0{,}0345} = 0{,}4058$$

$$P(H_3/A) = \frac{P(H_3)P(A/H_3)}{P(A)} = \frac{0{,}40 \cdot 0{,}02}{0{,}0345} = 0{,}2319$$

**Приклад 2.**

До магазину надходять однотипні вироби з трьох заводів, причому $i$-й завод постачає $m_i\%$ виробів. Серед виробів $i$-го заводу $n_i\%$ першого ґатунку. Куплено один виріб. Він виявився виробом першого ґатунку. Визначити ймовірність того, що куплений виріб виготовлено $j$-м $(i)$ заводом.

Маємо:.

| $m_1$ | $m_2$ | $m_3$ | $n_1$ | $n_2$ | $n_3$ | $j\ (i)$ |
|---|---|---|---|---|---|---|
| 50 | 30 | 20 | 70 | 80 | 90 | 1 |

Відповідь: Задача на формулу Байєса, для якої потрібно використання формули повної ймовірності. За умовами задачі, маємо три заводи, кожний з яких виготовляє вироби (дані надані в %) так, що $m_1 + m_2 + m_3 = 50 + 30 + 20 = 100\%$, та лише $n_i\%$ виробів першого ґатунку.

Позначимо: через $P(A)$ шукану ймовірність того, що виріб виявився першого ґатунку; через $P(B_i)$ ймовірність події $B_i -$ що виріб був від постачальника (заводу) номер $i$. Запишемо формулу повної ймовірності у загальному вигляді:

$$P(A) = \sum_{i=1}^{n} P(A \cdot B_i) = \sum_{i=1}^{n} P(A|B_i) \cdot P(B_i)$$

Умовна ймовірність того, що виріб від заводу $i$ буде першого ґатунку пообчислюємо за формулою (дані задані в умові задачі):

$$P(A|B_i) = \frac{n_i \%}{100\%}.$$

Однак, ймовірність того, що виріб від заводу $i$, дорівнює $P(B_i) = \frac{m_i\%}{100\%}$

Отже, повна ймовірність з підставленими даними для 1 варіанту
$P(A) = 0.70 * 0.50 + 0.80 * 0.30 + 0.90 * 0.20 = 0.35 + 0.24 + 0.18 = 0.77$

За формулою Байєса, умовна ймовірність події $B_i$ за умови, що відбулась подія $A$, тобто, що виріб був від заводу номер $i$ та, що виріб виявився першого ґатунку :

$$P(B_i|A) = \frac{P(A \cdot B_i)}{P(A)} = \frac{P(A|B_i) \cdot P(B_i)}{P(A)}$$

Ця формула є наслідком формули множення ймовірностей
$P(A \cdot B_i) = P(A|B_i) \cdot P(B_i) = P(B_i|A) \cdot P(A)$

Для нашого варіанту, потрібно узнати ймовірність лише для першого заводу (позначення індексу $i$, $j$, чи іншою літерою, значення не має):

$$P(B_1|A) = \frac{P(AB_1)}{P(A)} = \frac{0.35}{0.35 + 0.24 + 0.18} = \frac{0.35}{0.77} = \frac{35}{77} = \frac{5}{11} = 0.4545 \ldots$$

Аналогічно, можна обчислити так інші
$$P(B_2|A) = \frac{P(A \cdot B_2)}{P(A)} = \frac{0.24}{0.35 + 0.24 + 0.18} = \frac{0.24}{0.77} = \frac{24}{77} = 0.3117$$
$$P(B_3|A) = \frac{P(A \cdot B_3)}{P(A)} = \frac{0.18}{0.35 + 0.24 + 0.18} = \frac{0.18}{0.77} = \frac{18}{77} = 0.2338$$

### 3.4. Формула Бернуллі.

Часто зустрічаються задачі, в яких один і той же дослід (випробування) повторюються неодноразово. В результаті кожного випробування може з'явитися або не з'явитися подія $A$. Цікавить не результат окремого досліду, а загальне число появ події $A$ в серії дослідів.

Наведемо умови для випробувань, що повторюються:

1) кількість випробувань – число обмежене, позначаємо $n$;

2) ймовірність появи події $A$ в кожному випробуванні не залежить від результатів інших випробувань, такі випробування називаються незалежними щодо події $A$;

3) ймовірність появи події $A$ в кожному випробуванні – величина постійна, позначимо $P(A) = p$.

Повторення випробувань може бути обумовлене повторенням у часі випробування одного і того ж об'єкта, як згадане підкидання монети однією людиною. Проте, повторення випробувань може відбуватися незалежно від часу, коли випробуванню піддається декілька однакових об'єктів.

Позначимо символом $q$ ймовірність не появи події $A$: $P(\bar{A}) = 1 - p = q$.

Нехай потрібно визначити ймовірність того, що подія $A$ відбулася рівно $m$ раз у $n$ випробуваннях. Така ймовірність позначається $P_{m.n}(A)$. Очевидно, що в загальному випадку $0 \leq m \leq n$. Позначимо $B_n$ подію, яка полягає в тому, що подія А відбулася рівно $m$ раз у $n$ випробуваннях. Нехай $A_i$ поява події в $i$-м випробуванні; $\overline{A_k}$ – не поява події в $k$-м випробуванні. Тоді:

$$B_m = A_1 A_2 \ldots A_m \bar{A}_{m+1} \ldots \bar{A}_n + A_1 A_2 \ldots A_{m-1} \bar{A}_m \ldots \bar{A}_{n-1} A_n + \bar{A}_1 \ldots \bar{A}_{n-m} A_{n-m+1} \ldots A_n$$

У кожен добуток А входить $m$ раз, $\bar{A}$ входить $(n - m)$ раз. Число всіх добутків дорівнює числу комбінацій з $n$ елементів по $m$ елементів: $C_n^m$. Ймовірність кожного добутку визначається за теоремою множення ймовірностей незалежних подій і дорівнює $p^m q^{n-m}$. Добутки між собою несумісні (у кожній парі добуток є протилежні події з однаковим індексом), тому ймовірність суми добутків визначається за теоремою складання несумісних подій:

$$P(B_m) = P(A_1 A_2 \ldots A_m \bar{A}_{m+1} \ldots \bar{A}_n) + \cdots + P(\bar{A}_1 \ldots \bar{A}_{n-m} A_{n-m+1} \ldots A_n) =$$
$$= p^m q^{n-1} + p^m q^{n-1} + \cdots + p^m q^{n-1}$$

Таким чином,
$$P_{m,n}(A) = C_n^m p^m q^{n-m} = \frac{n!}{m!\,(n-m)!} p^m (1-p)^{n-m} \qquad (9)$$

Формула (9) називається **формулою Бернуллі.**
Ймовірність того, що у *n* випробуваннях подія A наступить:
менше *m* разів: $P_{0,n}(A) + P_{1,n}(A) + \cdots + P_{m-1,n}(A)$;
більше *m* разів: $P_{m+1,n}(A) + P_{m+2,n}(A) + \cdots + P_{n,n}(A)$;
не менше *m* разів: $P_{m,n}(A) + P_{m+1,n}(A) + \cdots + P_{n,n}(A)$;
не більше *m* разів: $P_{0,n}(A) + P_{1,n}(A) + \cdots + P_{m,n}(A)$;

На наступному рисунку, показані біноміальні розподіли для різної кількості випробувано та $p = 0.8$. Очевидно, для $n = 0$ спроб, ймовірність 0 успіхів дорівнює 1, тобто одне ненульове значення. Загалом, ненульових значень $n + 1$, тобто, $0..n$. Оскільки сума усіх ймовірностей дорівнює одиниці для усіх $n$, то розподіли стають довшими та нижчими.

Середнє значення
$$\bar{k} = \sum_{k=0}^{n} k\, p_k$$

Для біноміального розподілу, $\bar{k} = np$.

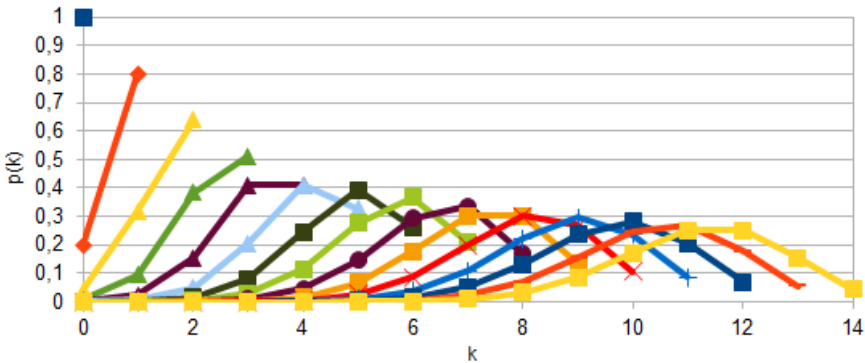

• **Залежність ймовірності $p(k)$ кількості успіхів $k$ при кількості спроб $n = 0, 1, \ldots 14$) та ймовірності успіху із однієї спроби $p = 0.8$.**

**Приклад 1**. Монету кинуто 6 разів. Знайти ймовірність того, що герб випаде

0, 1, 2, 3, 4, 5, 6 раз.

Відповідь: Позначимо, що подія випаде герб це А, ймовірність цієї події $p$ чи її не можливість $q$ однакові, тобто, маємо $n = 6$, $p = \frac{1}{2}$, $q = \frac{1}{2}$.

Тоді по формулі (9):
$$P_{m.n}(A) = C_n^m p^m q^{n-m} = \frac{n!}{m!(n-m)!} p^m (1-p)^{n-m}$$

маємо записати:

$P_6(0) = C_6^0 \left(\frac{1}{2}\right)^0 \left(\frac{1}{2}\right)^6 = \frac{1}{64};$

$P_6(1) = C_6^1 \left(\frac{1}{2}\right)^1 \left(\frac{1}{2}\right)^5 = \frac{6}{1} \cdot \frac{1}{2} \cdot \frac{1}{32} = \frac{6}{64}$

$P_6(2) = C_6^2 \left(\frac{1}{2}\right)^2 \left(\frac{1}{2}\right)^4 = \frac{6 \cdot 5}{1 \cdot 2} \cdot \frac{1}{4} \cdot \frac{1}{16} = \frac{15}{64}$

$P_6(3) = C_6^3 \left(\frac{1}{2}\right)^3 \left(\frac{1}{2}\right)^3 = \frac{6 \cdot 5 \cdot 4}{1 \cdot 2 \cdot 3} \cdot \frac{1}{8} \cdot \frac{1}{8} = \frac{20}{64}$

$P_6(4) = C_6^4 \left(\frac{1}{2}\right)^4 \left(\frac{1}{2}\right)^2 = \frac{15}{64}$

$P_6(5) = C_6^5 \left(\frac{1}{2}\right)^5 \left(\frac{1}{2}\right)^1 = \frac{6}{64}$

$P_6(6) = C_6^6 \left(\frac{1}{2}\right)^6 \left(\frac{1}{2}\right)^0 = \frac{1}{64}$

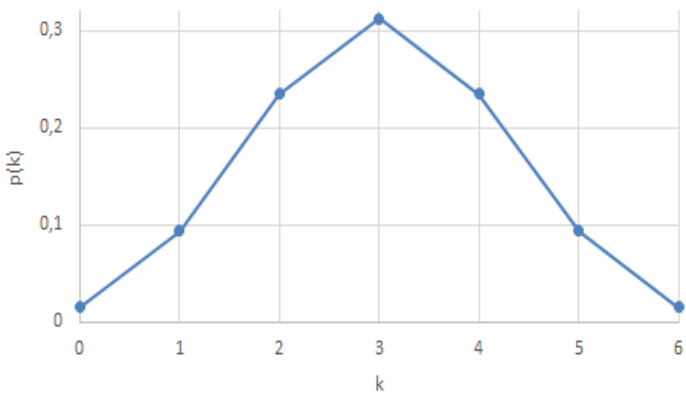

• **Приклад 2**.

Для стрільця, який виконує вправи в тирі, ймовірність влучити в «яблучко» дорівнює $p$=0,23. Спортсмен зробив N=10 пострілів. Знайти ймовірність таких подій:

А – «Точно одне влучення»;

В – «Точно два влучення»;

С – «Хоч одне влучення»;

Д – «Не менше М=4 влучень».

Відповідь:

Розподіл значень ймовірності кількості влучень описується формулою для біноміального розподілу:
$$p_k = C_N^k p^k (1-p)^{N-k} = C_N^k p^k q^{N-k}$$
де $p_k$ – ймовірність того, що при $N$ випробуваннях із постійною ймовірністю «успіху» $p$, відбудеться рівно $k$ «успіхів» (подій). Ймовірність того, що подія не відбудеться, дорівнює $q = 1 - p$. Ймовірність будь-якої послідовності з $k$ успіхів та, відповідно, $N - k$ «неуспіхів», дорівнює $p^k q^{N-k}$, а кількість різних комбінацій дорівнює $C_N^k$.

Для цієї задачі, питання А) та В) відрізняються тільки значенням $k$, тобто, треба визначити $p_1$ та $p_2$ за формулою:
$p_k = C_N^k p^k (1-p)^{N-k}$, де $C_N^k$ можна знайти за формулою:
$$C_n^k = \frac{n!}{k!\,(n-k)!}$$
чи за формулою: $C_N^k = \frac{N}{1} \frac{N-1}{2} \frac{N-2}{3} \dots \frac{N-(k-1)}{k}$ і далі
$$p_k = C_N^k p^k q^{N-k} = \frac{N}{1} \frac{N-1}{2} \frac{N-2}{3} \dots \frac{N-(k-1)}{k} p^k q^{N-k}$$
Підставимо дані умови задачі: $p$=0,23; N=10
$$q = 1 - p = 1 - 0{,}23 = 0{,}77$$
та отримаємо:
$$C_{10}^1 = \frac{10!}{1!\,(10-1)!} = \frac{10!}{9!} = \frac{9! \cdot 10}{9!} = 10$$
$p_1 = C_{10}^1 p^1 q^{10-1} = 10 \cdot 0{,}23^1 \cdot 0{,}77^9 = 2{,}3 \cdot 0{,}095 \approx 0{,}22$
та для $C_{10}^2$ і $p_2$:
$$C_{10}^2 = \frac{10!}{2!\,(10-2)!} = \frac{10!}{2 \cdot 8!} = \frac{8! \cdot 9 \cdot 10}{8! \cdot 2} = \frac{90}{2} = 45$$

$$p_2 = C_{10}^2 p^2 q^{10-2} = 45 \cdot 0{,}23^2 \cdot 0{,}77^8 = 45 \cdot 0{,}053 \cdot 0{,}12 \approx 0{,}29$$

Для випадку С – «Хоч одне влучення» можна спочатку розглянути ймовірність жодного влучення як $p_0$, тоді протилежним цієї позиції буде $q$ – «Хоч одне влучення», яке ми знайдемо з формули: $q = 1 - p_0$. Знайдемо $p_0$ за формулою:
$$p_k = C_N^k p^k (1-p)^{N-k}$$
$$p_0 = C_{10}^0 p^0 (1-p)^{10} =$$

За означенням вважається, що $C_n^0 = 1$, $a^0 = 1$ при $a \neq 0$, тоді
$$= 1 \cdot 1 \cdot 0{,}77^{10} = 0.073$$

Маємо наступний вигляд формули для $k = 0$:
$$p_0 = C_N^0 p^0 q^{N-0} = 1 \cdot 1 \cdot q^N = q^N$$

Таким чином, ймовірність для випадку С – «Хоч одне влучення»
$$q = 1 - p_0 = 1 - 0{,}073 = 0{,}927.$$

Для випадку Д – «Не менше М влучень», фраза «не менше» означає, що кількість влучень може бути більшою. Ми пам'ятаємо, що коли нас цікавить хоч одна подія з декількох за означенням - «Сумою (об'єднанням) подій А і В називається подія, яка полягає в тому, що відбувається хоча б одна з даних подій». Таким чином, ймовірність можна знайти за формулою:
$$p_M + p_{M+1} + \cdots + p_N,$$

тобто, сума всіх ймовірностей з кількістю влучень від М до N. Також можна знайти цю ймовірність і за формулою для протилежної події, що буде коротше (в залежності від М та N.) $q = 1 - p_0 - p_1 - \cdots - p_{M-1}$

Для варіанту М=4 маємо:
$$P(\text{М} \geq 4) = p_4 + p_5 + p_6 + p_7 + p_8 + p_9 + p_{10}$$
$$\approx 0{,}12 + 0{,}04 + 0{,}01 + 0{,}001 \ldots \approx 0{,}171$$

• **Приклад 3.**

Монету кидають, доки герб не випаде $n = 3$ рази. Визначити ймовірність того, що цифра при цьому випаде $m = 2$ рази.

Відповідь: Розподіл числа неуспіхів $m$ до $n$ успіхів при ймовірності успіху $p$ (та ймовірності $q = 1 - p$ «неуспіху» = «невдачі») називають «від'ємним біноміальним» розподілом (Negative Binomial, NB):
$$p_{NBm,n} = C_{m+n-1}^m p^n q^m = \frac{(m+n-1)!}{m!(n-1)!} \cdot p^n q^m = \frac{n}{m+n} p_{Bm,n}$$

При вирішенні такого типу задач, частою помилкою є використання біноміального (Binomial, B) закону замість «від'ємного біноміального», тобто, формули «загалом було $n + m$ випробувань, із них $m$ «неуспіхів»:

$$p_{Bm,n} = C_{m+n}^m p^n q^m = \frac{(m+n)!}{m! \, n!} \cdot p^n q^m = \frac{m+n}{n} p_{NBm,n}$$

Тобто, ці обидві формули розв'язують дві різні постановки задач. «Біноміальний розподіл» фіксує загальну кількість випробувань $N = n + m$, отже, $0 \leq m \leq N$, у той час, як для «від'ємного біноміального розподілу», кількість спроб (до $n$ «успіхів») необмежена: $0 \leq m < \infty$.

Оскільки останнім, за умовою повинен відбутися «успіх» (напр. випасти герб), що фіксовано, то до того, «неуспіх» (цифра) випаде $m$ разів, а «успіх» (герб) – $j = (n-1)$раз. Тоді ймовірність (до завершального успіху) буде підпорядкована «звичайному» біноміальному закону, але наприкінці потрібно ще помножити на ймовірність завершальної події:

$$p_{m,n} = \left(C_{m+j}^j p^j q^m\right)p = (C_{m+n-1}^{n-1} \cdot p^{n-1} q^m)p = C_{m+n-1}^{n-1} \cdot p^n q^m$$
$$= C_{m+n-1}^m \cdot p^n q^m$$

Щоб знайти ймовірність $p_{m,n}$ − що цифра випаде $m$ разів, треба визначити значення ймовірності успіху $p$(випадіння герба) та неуспіху $q$ (випадіння цифри). Такого виду задачі вже розглядалися і визначалось, що при киданні монети, $p = q = \frac{1}{2}$, тобто, маємо:

$$p^n q^m = \left(\frac{1}{2}\right)^n \cdot \left(\frac{1}{2}\right)^m = \left(\frac{1}{2}\right)^{n+m} = \frac{1}{2^{n+m}}$$

Тепер знайдемо $C_{m+n-1}^m$ за формулою

$$C_n^k = \frac{n!}{k! \, (n-k)!}$$

підставивши значення $n = 3; m = 2$ (дані з умов задачі) і отримаємо:

$$C_{m+n-1}^m = C_{2+3-1}^2 = C_4^2 = \frac{4!}{2!(4-2)!} = \frac{4!}{2! \cdot 2!} = \frac{3 \cdot 4}{2} = 6$$

Таким чином, отримаємо:

$$p_{m,n} = \frac{C_{m+n-1}^m}{2^{m+n}} = \frac{6}{2^{3+2}} = \frac{6}{2^5} = \frac{6}{32} \approx 0{,}19$$

**Приклад 4.**

Ймовірність виграшу в лотерею на один квиток дорівнює $p = 0{,}3$. Куплено $n = 10$ квитків. Знайти найімовірніше число виграшних квитків і відповідну ймовірність.

Відповідь: Ця задача на застосування біноміального розподілу. Але потрібно обчислити не декілька значень, а лише одне, яке відповідає «моді» розподілу (найбільшому значенню ймовірності). Позначимо через $m$ найбільш ймовірну кількість подій виграшу.

*Примітка. Щоб визначити значення при якому ймовірність буде максимальною, знайдемо відношення ймовірності поточного значення $p_m$ до ймовірності попереднього. При збільшенні, воно повинно бути не меншим за одиницю, тому маємо записати:*

$$p_m = C_n^m p^m q^{n-m}$$
$$p_{m-1} = C_n^{m-1} p^{m-1} q^{n-(m-1)}$$
$$\frac{p_m}{p_{m-1}} = \frac{C_n^m p^m q^{n-m}}{C_n^{m-1} p^{m-1} q^{n-(m-1)}} =$$

Якщо застосуємо формулу $C_n^k = \frac{n!}{k!(n-k)!}$ та скористаємось загальними властивостями ступенів $\frac{a^n}{a^m} = a^{n-m}$, то отримаємо:

$$\frac{p_m}{p_{m-1}} = \frac{(m-1)!}{m!} \frac{(n-m+1)!}{(n-m)!} \frac{p}{q} = \frac{n-m+1}{m} \cdot \frac{p}{q} \geq 1$$
$$(n-m+1)p \geq mq$$

розкриємо дужки та перегрупуємо

$$(n+1)p - mp \geq mq$$
$$(n+1)p \geq mq + mp = m(q+p) = m \cdot 1 = m \Rightarrow$$
$$m \leq (n+1)p$$

Таким чином, найбільш ймовірна кількість подій виграшу $m$ визначається подвійною нерівністю

$$(n+1)p - 1 \leq m \leq (n+1)p$$

маючи на увазі, що $(q = 1 - p)$ можемо визначити $m$ такою нерівністю:

$$np - q \leq m \leq np + p$$

Відповідна ймовірність обчислюється за звичайною формулою для $k = m$: $p_m = C_n^m p^m q^{n-m}$, де $C_n^m = \frac{n!}{m!(n-m)!}$

Підставимо дані для першого варіанту – $p = 0.3, n = 10$,
$(n + 1)p = (10 + 1) \cdot 0.3 = 3.3$,
Отже, $3.3 - 1 \leq m \leq 3.3$,
$m = 3$,

Підставимо знайдене значення $m$ і дані значення $p = 0.3, n = 10$ в формулу: $p_m = C_n^m p^m q^{n-m}$ і отримаємо:

$$p_3 = C_{10}^3 0.3^3 0.7^7 = \frac{10 \cdot 9 \cdot 8}{1 \cdot 2 \cdot 3} \cdot 0.027 \cdot 0.0823543 = 120 \cdot 0.027 \cdot 0.0823543$$
$$= 0.266827932 \approx 0.2668 \approx 27.7\%$$

**3.5. Локальна теорема Муавра – Лапласа.**
**Теорема.**

Нехай в кожному з $n$ незалежних випробувань імовірність настання події А однакова і дорівнює $p$ *(0<p<1)*, тоді справедлива наближена рівність

$$P_n(k) = \frac{1}{\sqrt{npq}} \varphi\left(\frac{k-np}{\sqrt{npq}}\right) \qquad (10)$$

де $\varphi(x) = \frac{1}{\sqrt{2\pi}} e^{\frac{-x^2}{2}}$.

Функція $\varphi(x)$ парна: $\varphi(-x) = \varphi(x)$. Таблиця функції $\varphi(x)$ наведена в додатку 1. В таблиці по значенню $x$, яке визначається по формулі $x = \frac{k-np}{\sqrt{npq}}$, визначаємо $\varphi(x)$. Формула (10) дає добре наближення, якщо $n$ достатньо велике, $p$ та $q$ не дуже близькі до нуля ( *npq*>9).

• **Приклад.** Ймовірність успіху у кожному випробуванні дорівнює 0,25.

Яка ймовірність того, що при 300 випробуваннях успішними будуть:

а) рівно 75 випробувань;

Відповідь: За умовою $n$=300, $k$ =75, $p$ =0,25, тоді $q = 1 - p = 0{,}75$.

Спочатку обчислимо вираз:
$$\frac{k-np}{\sqrt{npq}} = \frac{75 - 300 \cdot 0{,}25}{\sqrt{300 \cdot 0{,}25 \cdot 0{,}75}} = \frac{0}{75} = 0.$$

Це ми знайшли значення $x = 0$ і щоб визначити $\varphi(x)$ дивимось у таблицю значень функції нормального розподілу Гаусса Лапласа $\varphi(x)$:

| Цілі і десяті частини | Соті частини $x$ | | | | | | | | | |
|---|---|---|---|---|---|---|---|---|---|---|
| $x$ | 0 | 1 | 2 | 3 | 4 | 5 | 6 | 7 | 8 | 9 |
| 0,0 | 0,3989 | 0,3989 | 0,3989 | 0,3988 | 0,3986 | 0,3984 | 0,3982 | 0,3980 | 0,3977 | 0,3973 |
| 0,1 | 0,3970 | 0,3965 | 0,3961 | 0,3956 | 0,3951 | 0,3945 | 0,3939 | 0,3932 | 0,3925 | 0,3918 |

Значенню $x = 0$ відповідає $\varphi(x) = 0{,}3989$, яке підставимо у

формулу (10):
$$P_{300}(75) \approx \frac{1}{7,5}\varphi(0) = \frac{0,3989}{7,5} = 0,0532$$

б) рівно 85 випробувань?

Відповідь: Аналогічно, як і в варіанті а). По умовам задачі маємо:

$n$=300, $k$=75, $p$=0,25, тоді $q = 1 - p = 0,75$, обчислимо вираз:
$$\frac{k - np}{\sqrt{npq}} = \frac{85 - 300 \cdot 0,25}{\sqrt{300 \cdot 0,25 \cdot 0,75}} = \frac{10}{75} = 1,33.$$

Дивимось по таблиці яке значення $\varphi(x)$ відповідає отриманому $x = 1,33$

| Цілі і десяті частини | Соті частини $x$ | | | | | | | | | |
|---|---|---|---|---|---|---|---|---|---|---|
| $x$ | 0 | 1 | 2 | 3 | 4 | 5 | 6 | 7 | 8 | 9 |
| 1,2 | 0,1942 | 0,1919 | 0,1895 | 0,1872 | 0,1849 | 0,1826 | 0,1804 | 0,1781 | 0,1758 | 0,1736 |
| 1,3 | 0,1714 | 0,1691 | 0,1669 | 0,1647 | 0,1626 | 0,1604 | 0,1582 | 0,1561 | 0,1539 | 0,1518 |

Значенню $x = 1,33$ відповідає $\varphi(x) = 0,1647$, яке підставимо у формулу (10):
$$P_{300}(85) \approx \frac{1}{7,5}\varphi(1,33) = \frac{0,1647}{7,5} = 0,0219.$$

### 3.6. Інтегральна теорема Муавра-Лапласа.
**Теорема.**

Ймовірність того, що в $n$ незалежних випробуваннях, в кожному з яких подія А може відбутись з ймовірністю $p$ ($0<p<1$), подія А відбудеться не менше $k_1$ і не більше $k_2$ раз, наближено дорівнює :

$$P_n\{k_1 \leq k \leq k_2\} \approx \Phi\left(\frac{k_2-np}{\sqrt{npq}}\right) - \Phi\left(\frac{k_1-np}{\sqrt{npq}}\right), \quad (11)$$

де $\Phi(x) = \frac{1}{\sqrt{2\pi}}\int_0^x e^{\frac{t^2}{2}}dt$, таблиця функції $\Phi(x)$ наведена в додатку 2.

• **Приклад.** Ймовірність виходу з ладу за час t одного приладу дорівнює 0,1. Визначити ймовірність того, що за час t зі 100 приладів вийде з ладу:

**а)** від 6 до 18 приладів.

Відповідь:

За умовою $n = 100, k_2 = 18, k_1 = 6, p = 0{,}1, q = 1 - p = 0{,}9$. Знаходимо:
$$\frac{k_2 - np}{\sqrt{npq}} = \frac{18 - 100 \cdot 0{,}1}{\sqrt{100 \cdot 0{,}1 \cdot 0{,}9}} = \frac{8}{3} = 2{,}66;$$

$$\frac{k_1 - np}{\sqrt{npq}} = \frac{6 - 100 \cdot 0{,}1}{\sqrt{100 \cdot 0{,}1 \cdot 0{,}9}} = -\frac{4}{3} = -1{,}33;$$

Дивимось по таблиці які значення $\Phi(x)$ відповідають отриманим $x = 2{,}66$ та $x = 1{,}33$.

| 1,3 | 0,4032 | 0,4049 | 0,4066 | 0,4082 | 0,4099 | 0,4115 | 0,4131 | 0,4147 | 0,4162 | 0,4177 |
| 2,6 | 0,4953 | 0,4955 | 0,4956 | 0,4957 | 0,4959 | 0,4960 | 0,4961 | 0,4962 | 0,4963 | 0,4964 |

Підставимо отримані значення у формулу (11):
$$P_{100}\{6 \leq k \leq 18\} \approx \Phi(2{,}66) - \Phi(-1{,}33) = \Phi(2{,}66) + \Phi(1{,}33) =$$
$$= 0{,}49609 + 0{,}40824 = 0{,}90433.$$

**б)** не менше 20.
Відповідь: По умовам задачі маємо:
$$n = 100, k_2 = 100, k_1 = 20, p = 0{,}1, q = 1 - p = 0{,}9.$$
$$\frac{k_2 - np}{\sqrt{npq}} = \frac{100 - 100 \cdot 0{,}1}{\sqrt{100 \cdot 0{,}1 \cdot 0{,}9}} = \frac{90}{3} = 30;$$
$$\frac{k_1 - np}{\sqrt{npq}} = \frac{20 - 100 \cdot 0{,}1}{\sqrt{100 \cdot 0{,}1 \cdot 0{,}9}} = \frac{10}{3} = 3{,}3;$$

Дивимось по таблиці які значення $\Phi(x)$ відповідають отриманим $x = 30$ та $3{,}3$.

| 3,3 | 0,4995 | 0,4995 | 0,4995 | 0,4996 | 0,4996 | 0,4996 | 0,4996 | 0,4996 | 0,4996 | 0,4997 |
| 3,9 | 0,5000 | 0,5000 | 0,5000 | 0,5000 | 0,5000 | 0,5000 | 0,5000 | 0,5000 | 0,5000 | 0,5000 |

Підставимо отримані значення у формулу (11):
$$P_{100}\{20 \leq k \leq 100\} \approx \Phi(30) - \Phi(3{,}3) = 0{,}5 - 0{,}4995 = 0{,}0005.$$

**3.7. Наближена формула Пуассона при повторенні випробувань.**

Формула Бернуллі не може бути застосована, якщо число випробувань дуже велике, оскільки при цьому неприпустимо велика похибка обчислень.

Нехай $n \to \infty$, $p \to 0$, але величина $p \cdot n = \lambda$ зберігає постійне значення. У формулу Бернуллі підставимо $p = \frac{\lambda}{n}$ і визначимо межу при $n \to \infty$.

$$\lim_{n \to \infty} P_{m,n}(A) = \lim_{n \to \infty} \frac{n!}{m!(n-m)!} \frac{\lambda^m}{n^m}\left(1-\frac{\lambda}{n}\right)^{n-m} =$$

$$= \frac{\lambda^m}{m!} \lim_{n \to \infty} \frac{n!}{(n-m)!} \left(\left(1-\frac{\lambda}{n}\right)^{\frac{-n}{\lambda}}\right)^{\frac{\lambda}{n}(n-m)} = \frac{\lambda^m}{m!} \cdot 1 \cdot e^{\lim_{n \to \infty}\frac{-\lambda n + \lambda m}{n}} = \frac{\lambda^m}{m!}e^{-\lambda}$$

Де $\lim_{n \to \infty}\frac{n!}{(n-m)! \cdot n^m} = 1$ як відношення нескінченно великих величин одного порядку зростання. Отримана наближена формула:

$$P_{m,n}(A) \approx \frac{\lambda^m}{m!}e^{-\lambda}, \text{ де } \lambda = p \cdot n \qquad (12)$$

називається формулою Пуассона і застосовується, коли $n$ – велике, а ймовірність $p$ – мала (орієнтування $p < 0{,}1$, $npq < 9$. Формула Пуассона використовується в завданнях, де розглядаються рідкісні події.

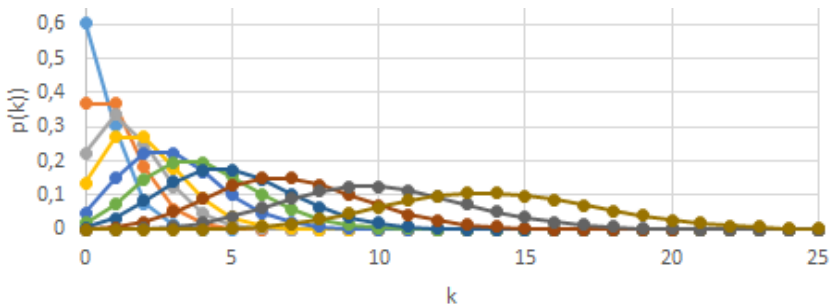

**Розподіл Пуасона для різних значень параметра $\lambda$=0.5, 1, 1.5, 2, 3, 4, 6, 7, 10, 14.**

• **Приклад 1.** Середній брак при виробництві продукції становить 0,1%. Перевіряється партія з 1000 деталей. Яка ймовірність того, що бракованих буде від 2 до 4 деталей?
Відповідь: За умовою $n$=1000, $p$=0,001 отже $\lambda = p \cdot n = 0{,}001 \cdot 1000 = 1$

$$P\{2 \leq m \leq 4\} = \frac{1^2}{2!}e^{-1} + \frac{1^3}{3!}e^{-1} + \frac{1^4}{4!}e^{-1} = 0{,}26058$$

• **Приклад 2.**
Ймовірність «відмови» у роботі телефонної станції при кожному виклику дорівнює $p$= 0,002. Надійшло $n = 1000$ викликів. Визначити ймовірність $m$= 7 «відмов».

Відповідь: Оскільки виконується умова $p \ll 1$, $n \gg 1$, тобто, маємо велике значення $n$ та мале значення ймовірності, то ймовірність кількості подій обчислюють не за біноміальним розподілом, а за його наближенням, яке називають «розподіл Пуасона»:

$$p_m = \frac{\lambda^m}{m!}e^{-\lambda} = \frac{\lambda^m}{m!}\exp(-\lambda)$$

Де параметр $\lambda = np$ — математичне очікування (теоретична середня кількість) подій, $e \approx 2{,}718281828\ldots$ — основа натуральних логарифмів.

Для $p = 0.002, n = 1000, m = 7$, отримаємо:
$\lambda = np = 1000 \cdot 0{,}002 = 2$

$$p_7 = \frac{2^7}{7!}e^{-2} = \frac{128}{5040} \cdot \frac{1}{e^2} \approx 0{,}0254 \cdot 0{,}135 \approx 0{,}00343$$

• **Приклад 3.**

Книги ($n = 7$) випадково ставлять на полицю. Яка кількість комбінацій того, що ($m = 3$) вибрані книги займуть перші місця?

Відповідь: $m!(n-m)!$, адже $m!$ Перестановок для обраних книг, та 4 для тих, що залишилися. Для даних чисел, $3!(7-3)! = 3!4! = 6 * 24 = 144$. Загальна кількість можливих перестановок $n!$. Отже, ймовірність такої випадкової події $P = \frac{m!(n-m)!}{n!} = 1/C_n^m$.

Якщо у завданні не «перші m місць», а «будь-які m місця поруч», то число перестановок (та відповідну ймовірність) потрібно збільшити у $(n - m + 1)$ разів, бо ці $m$ позицій можуть починатися з $1, 2, \ldots (n - m + 1)$.

Якщо «по колу» (наприклад, гості за круглим столом), то початок групи з $m$ місць поспіль може починатися з будь якого номера (1,2,…n), і тому потрібно домножати на $n$.

## 4. ВИПАДКОВІ ВЕЛИЧИНИ.
### 4.1. Дискретні випадкові величини.

**Означення.** Випадковою величиною називають таку величину, яка внаслідок випробування може набути лише одне числове значення з можливих, яке зумовлене результатом експерименту і заздалегідь не можуть бути врахованими.

Отже, випадковою величиною, пов'язаною з певним дослідом, називають величину, яка під час кожного проведенням досліду може набувати того чи другого числового значення, залежно від випадку.

Між випадковими подіями і випадковими величинами є тісний зв'язок. Випадкова подія – це якісна характеристика випадкового результату досліду, а випадкова величина – його кількісна характеристика. Випадкові величини за певними властивостями поділяються на дискретні та неперервні.

**Означення.** Дискретною випадковою величиною (ДВВ) називають таку величину, яка внаслідок випробування може набути відокремлених, ізольованих одне від одного числових значень з відповідними ймовірностями.

Всі можливі значення дискретної випадкової величини можуть бути перенумеровані: $x_1, x_2, ..., x_n, ....$

**Означення.** Неперервною випадковою величиною (НВВ) називають величину, яка може набувати будь-якого числового значення з певного обмеженого інтервалу (a, b) або необмеженого інтервалу $(-\infty; +\infty)$. Наприклад, випадкова величина X – час безвідмовної роботи приладу, вона неперервна, оскільки її можливе значення t > 0.

**Означення.** Співвідношення, яке встановлює зв'язок між можливими значеннями випадкової величини і ймовірностями, з якими приймають ці значення, називають законом розподілу ймовірностей випадкової величини. Для дискретної випадкової величини X закон розподілу може бути заданий у вигляді таблиці або графіку.

У першому випадку, коли закон розподілу задається у вигляді таблиці, його називають рядом розподілу ймовірностей випадкової величини X:

| $X$ | $x_1$ | $x_2$ | ... | $x_n$ | ... |
|---|---|---|---|---|---|
| $P$ | $p_1$ | $p_2$ | ... | $p_n$ | ... |

У першому рядку таблиці записують усі можливі значення випадкової величини, а в другому − відповідні їм ймовірності. Оскільки події {X= $x_1$}, {X=$x_2$}, ..., {X= $x_n$} становлять повну групу несумісних подій, то за теоремою додавання ймовірностей маємо:

$$\sum_{k=1}^{\infty} p_k = 1 \qquad (1)$$

тобто, сума ймовірностей усіх можливих значень випадкової величини дорівнює одиниці.

За рядом розподілу (1) можна побудувати функцію розподілу дискретної випадкової величини X:

$$F(x) = \sum_{x_k < x} p_k. \qquad (2)$$

де підсумування поширюється на ті індекси $k$, для яких $x_k < x$.

Графічна форма закону розподілу називається багатокутником розподілу: по осі абсцис відкладаємо можливі значення $x_k$ випадкової величини $X$, а по осі ординат – імовірності $p_k$ цих значень; точки ($x_k$, $p_k$) послідовно з'єднуємо відрізками прямих.

• **Приклад.** В партії з шести деталей чотири стандартні. Навмання вибрано три деталі. Знайти:

1) ряд розподілу дискретної випадкової величини X – числа стандартних деталей серед відібраних;

2) функцію розподілу $F(x)$;

3) $P\{1{,}5 \le X < 2{,}5\}$.

Відповідь: 1) Знайти ряд розподілу дискретної випадкової величини X – числа стандартних деталей серед відібраних.

По умові задачі випадкова величина $X$ – число стандартних деталей серед відібраних – може приймати такі значення: $x_1 = 1, x_2 = 2, x_3 = 3$, тобто, три стандартні деталі це максимальна можливість, тому що взято тільки три деталі. Значення 0 відсутнє, тому що кількість нестандартних деталей 6-4=2, тому неможливо витягнути усі три нестандартні деталі. Знайдемо ймовірності можливих значень $X$: $P\{X = 1\}, P\{X = 2\}$ та $P\{X = 3\}$.

$P\{X = 1\}$ – ймовірність того, що випадкова величина дорівнюється 1 означає, що ми розглядаємо дві події, що здійснилися одночасно - було взято 1 стандартна деталь і дві нестандартні, тобто, ми маємо добуток двох комбінацій:

$C_2^2$ – дві нестандартні деталі із двох можливих помножену на $C_4^1$ – стандартна деталь із чотирьох можливих. Цей добуток ділимо на $C_6^3$ – три вибраних деталі із шести можливих. Згідно з формулою $C_n^k = \frac{n!}{k!(n-k)!}$ обчислюємо:

$$C_4^1 = \frac{4!}{1! \cdot (4-1)!} = \frac{1 \cdot 2 \cdot 3 \cdot 4}{1 \cdot 2 \cdot 3} = 4$$

$$C_6^3 = \frac{6!}{3! \cdot (6-3)!} = \frac{6!}{3! \cdot 3!} = \frac{4 \cdot 5 \cdot 6}{1 \cdot 2 \cdot 3}$$

Згідно з властивостями числа комбінацій без повторень $C_n^n = 1$, маємо, що $C_2^2 = 1$

Тепер можемо порахувати ймовірність $P\{X = 1\}$, для чого добуток комбінацій сприятливих події Х поділимо на кількість усіх можливих. Таким чином, маємо записати:

$$P\{X = 1\} = \frac{C_2^2 C_4^1}{C_6^3} = \frac{4}{1} \cdot \frac{3 \cdot 2 \cdot 1}{6 \cdot 5 \cdot 4} = \frac{1}{5} = 0{,}2,$$

Аналогічно знаходимо $P\{X = 2\}$ та $P\{X = 3\}$, маючи на увазі, що за означенням вважається $C_n^0 = 1$.

$$P\{X = 2\} = \frac{C_2^1 C_4^2}{C_6^3} = \frac{2}{1} \cdot \frac{4 \cdot 3}{1 \cdot 2} \cdot \frac{3 \cdot 2 \cdot 1}{6 \cdot 5 \cdot 4} = \frac{3}{5} = 0{,}6,$$

$$P\{X = 3\} = \frac{C_2^0 C_4^3}{C_6^3} = \frac{4}{1} \cdot \frac{3 \cdot 2 \cdot 1}{6 \cdot 5 \cdot 4} = \frac{1}{5} = 0{,}2.$$

Отриманий результат оформимо у вигляді ряду розподілу:

| X | 1 | 2 | 3 |
|---|---|---|---|
| P | 0,2 | 0,6 | 0,2 |

Контроль: 0,2+0,6+0,2=1.

2) Знайти функцію розподілу $F(x)$:

Розглянемо всі можливі значення випадкової величини Х, які вона може приймати і значення величини функції, що відповідають цим значенням випадкової величини.

При $x \leq 1$, маємо $F(x) = 0$, тобто, менше однієї деталі бути не може.

При $1 < x \leq 2$, маємо $F(x) = 0{,}2$ значення, що відповідає однієї деталі.

При $1 < x \leq 3$, маємо $F(x) = 0{,}2 + 0{,}6 = 0{,}8$ значення, що відповідає двом і трьом деталям.

При $x > 3$, маємо $F(x) = 0{,}2 + 0{,}6 + 0{,}2 = 1$.

Таким чином:
$$F(x) = \begin{cases} 0, & x \leq 1 \\ 0{,}2, & 1 < x \leq 2 \\ 0{,}8, & 2 < x \leq 3 \\ 1, & x > 3 \end{cases}$$

По приведеним значенням будуємо графік.

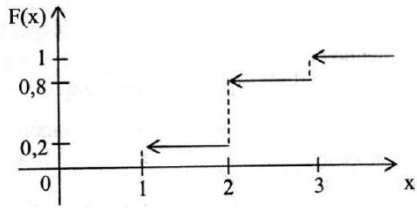

Стрілка на кінці відрізку означає, що ця точка не входить до інтервалу.

3) Знайти $P\{1{,}5 \leq X < 2{,}5\}$.

В цій задачі розглядають деталі, тобто, це штучні предмети і їх кількість надається цілими числами. На інтервалі $\{1{,}5 \leq X < 2{,}5\}$ X може мати тільки єдине ціле число 2. Дивимось у таблицю – X=2 відповідає значення P(2)=0,6.

Також це можна знайти інакше, дивлячись на значення, що відповідають окремим інтервалам. Так значення 2,5 відноситься до інтервалу $2 < x \leq 3$, якому відповідає 0,8, а 1,5 до $1 < x \leq 2$, якому відповідає значення 0,2. Таким чином, маємо записати відповідь:
$$P\{1{,}5 \leq X < 2{,}5\} = F(2{,}5) - F(1{,}5) = 0{,}8 - 0{,}2 = 0{,}6.$$

**4.2. Основні закони розподілу дискретних випадкових величин.**

**Біноміальний розподіл.**

Закон розподілу називається біноміальним, якщо ймовірність появи кожного значення випадкової величини X визначається за формулою Бернуллі:
$$P_{m,n}(A) = C_n^m p^m q^{n-m} = \frac{n!}{m!(n-m)!} p^m (1-p)^{n-m}.$$

• **Приклад**. В партії 10% нестандартних деталей. Навмання

вибрали чотири деталі. Знайти біноміальний закон розподілу дискретної випадкової величини X – числа нестандартних деталей серед чотирьох відібраних.

Відповідь: Випадкова величина X може приймати такі значення:
$x_1 = 0, x_2 = 1, x_3 = 2, x_4 = 3, x_5 = 4$. Запишемо параметри, що входять у формулу Бернуллі – параметр $n=4$ (кількість випробувань, тобто, вибраних деталей), p=0,1 (ймовірність появи події, яка полягає в тому, що в кожному випробуванні буде вибрана нестандартна деталь) її можна знайти із співвідношення 10% від 100% це теж саме, що 0,1 від 1 (ймовірність не може бути більше), і q=1–0,1=0,9 (ймовірність не появи цієї події). За формулою Бернуллі:

$$P_{m,n}(A) = C_n^m p^m q^{n-m} = \frac{n!}{m!(n-m)!} p^m (1-p)^{n-m}.$$

обчислюємо ймовірність події для кожного значення випадкової величини X:

$$P_4(0) = C_4^0 p^0 q^{4-0} = 1 \cdot 1 \cdot q^4 = 0{,}9^4 = 0{,}6561,$$

$$P_4(1) = C_4^1 p^1 q^{4-1} = C_4^1 p q^3 = \frac{4!}{1!(4-1)!} 0{,}1 \cdot 0{,}9^3 = 4 \cdot 0{,}1 \cdot 0{,}9^3 = 0{,}2916,$$

$$P_4(2) = C_4^2 p^2 q^{4-2} = \frac{4!}{2!(4-2)!} 0{,}1^2 \cdot 0{,}9^2 = 6 \cdot 0{,}01 \cdot 0{,}81 = 0{,}0486,$$

$$P_4(3) = C_4^3 p^3 q^{4-3} = \frac{4!}{3!(4-3)!} 0{,}1^3 \cdot 0{,}9^1 = 4 \cdot 0{,}001 \cdot 0{,}9 = 0{,}0036$$

$$P_4(4) = C_4^4 p^4 q^{4-4} = 1 \cdot 0{,}1^4 \cdot 1 = 0{,}0001$$

Контроль: $0{,}6561 + 0{,}2916 + 0{,}0486 + 0{,}0036 + 0{,}0001 = 1$

| X | 0 | 1 | 2 | 3 | 4 |
|---|---|---|---|---|---|
| p | 0,6561 | 0,2916 | 0,0486 | 0,0036 | 0,0001 |

**Розподіл Пуассона.**

Нехай дискретна випадкова величина X виражає кількість появ події A при масових випробуваннях ($n$ – велике), але при цьому в кожному випробуванні ймовірність появи події P(A)=$p$ мала, і виконується умова $np = \lambda = \text{const}$. У таких випадках закон розподілу випадкової величини задається формулою Пуассона

$$P_{m,n}(A) \approx \frac{\lambda^m}{m!} e^{-\lambda}, \text{ де } \lambda = p \cdot n$$

і називається розподілом Пуассона. Прикладами випадкових величин, що мають розподіл Пуассона, є:

число викликів на телефонній станції за час t;

число друкарських помилок у великому тексті;

число бракованих деталей у великій партії;

кількість крапель дощу (фотонів, заряджених частинок і т.д) в одиницю часу на фіксовану площу.

**Геометричний розподіл.**

Нехай відбуваються незалежні випробування, в кожному з яких ймовірність появи події А дорівнює $p$ ($0<p<1$), $q=1-p$. Випробування закінчуються, як тільки з'явиться подія А. Якщо подія з'явилася в k-му випробуванні, то у всіх попередніх вона не з'являлась. Позначимо через X – число випробувань, які потрібно зробити до першої появи події А. Можливими значеннями є натуральні числа $x_1 = 1, x_2 = 2, x_3 = 3$ і так далі. Якщо в перших $k-1$ випробуваннях подія не з'явилася, а в $k$ - му випробуванні відбулася, то:

$$P(X = k) = q^{k-1} \cdot p. \qquad (3)$$

Задаючи $k =1,2,3,…$ визначимо суму ймовірностей всіх значень випадкової величини X:

$$p + pq + pq^2 + pq^3 + \cdots + pq^n + \cdots$$

Отриманий числовий ряд – це нескінчена, що спадає геометрична прогресія зі знаменником $|q| < 1$. Тому ряд збігається, а сума ряду визначається за формулою:

$$S = \frac{p}{1-q} = \frac{p}{p} = 1.$$

• **Приклад**. Ймовірність того, що стрілець влучить в мішень при одному пострілі, дорівнює 0,8. Стрільцеві видають патрони доки він не промахнеться. Знайти закон розподілу дискретної випадкової величини X – числа виданих стрільцеві патронів.

Відповідь: Величина X має такі можливі значення:

$$x_1 = 1, x_2 = 2, …, x_k = k, … .$$

Знайдемо ймовірності цих значень:

Величина X приймає значення $x_1 = 1$, якщо стрілець не влучить в мішень при першому пострілі, тобто, щоб мати нуль влучень маємо мінімальну кількість пострілів =1. Ймовірність цього значення $P\{X = 1\} = 1 - 0{,}8 = 0{,}2$ (згідно формули $q = 1-$ *p*, де *q* – ймовірність протилежної події, тобто, того, що не влучить).

Щоб з'ясувати яке значення приймає величина X при $x_2 = 2$, тобто, якщо стрілець влучить в мішень при першому пострілі, а при другому промахнеться ми повинні враховувати, що здійснитися повинні одночасно дві події, тому маємо взяти добуток із ймовірності події, що влучить *p*=0,8, і ймовірності, що не влучить *q* =0,2. Тоді маємо записати:
$$P\{X = 2\} = 0{,}8 \cdot 0{,}2 = 0{,}16$$

Далі продовжимо знаходити по формулі
$$P(X = k) = q^{k-1} \cdot p.$$

Тепер знайдемо значення $P\{X = 3\} =$ та $P\{X = k\} =$ і розмістимо всі дані у таблиці.

Закон розподілу двовимірної випадкової величини має вигляд:

| X | 1 | 2 | 3 | … | K | … |
|---|---|---|---|---|---|---|
| P | 0,2 | 0,16 | 0,128 | … | $0{,}8^{k-1} \cdot 0{,}2$ | … |

### 4.3. Неперервні випадкові величини.

Випадкову величину X називають неперервною, якщо її функцію розподілу (CDF, cumulative distribution function) можна подати у вигляді: $F(x) = \int_{-\infty}^{x} f(t)dt,$ де $f(x) -$ деяка функція, яку називають щільністю розподілу ймовірностей (PDF, probability distribution function).

Якщо F(*x*) диференційована і похідна її обмежена, то випадкова величина X неперервна і має щільність розподілу ймовірностей:
$$f(x) = F'(x) = \frac{dF}{dx}.$$

Графік функції $f(x)$ називається кривою розподілу неперервної випадкової величини (див. рис.).

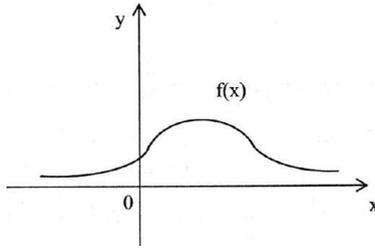

Закон розподілу неперервної випадкової величини може бути заданий графічно або аналітично $p_k=f(x_k)$ (за допомогою формули). Табличне завдання неможливе, оскільки ймовірність отримати будь-яке значення неперервної величини дорівнює нулю, що пов'язано не з неможливістю самої події (потрапляння в певну точку на числовій осі), а з нескінченно великою кількістю можливих випадків.

З огляду на це для неперервних випадкових величин (як, зрештою, і для дискретних) визначають імовірність потрапляння в деякий інтервал числової осі. Ймовірність потрапляння випадкової величини X в інтервал [a, b] визначають як ймовірність події P($a$≤X<$b$).

Для кількісного оцінювання закону розподілу випадкової величини (дискретної або неперервної) задають функцію розподілу ймовірностей випадкової величини, котру визначають як імовірність того, що випадкова величина X набуде значення, меншого від певного фіксованого числа $x$ і позначають F($x$)= P(X< $x$) або F($x$)= P(−∞<X< $x$).

Функцію розподілу F(x) називають інтегральною функцією розподілу ймовірностей випадкової величини. Знаючи функцію розподілу F(x), можна обчислити ймовірність потрапляння випадкової величини у деякий інтервал [a, b):

$$P(a < X < b) = F(b) - F(a). \qquad (4)$$

**Властивості функцій розподілу.**
1. $0 \leq F(x) \leq 1$.
   2. Функція розподілу неспадна: якщо $x_1 < x_2$, то $F(x_1) \leq F(x_2)$.
   3. Функція розподілу неперервна зліва:
$$\lim_{x \to x_1 - 0} F(x) = F(x_1).$$
   4. Імовірність потрапляння випадкової величини $X$ у проміжок

$[a, b]$:
$$P(a \leq X < b) = F(b) - F(a).$$
5. $P(X \geq x) = 1 - F(x)$.
6. $F(-\infty) = 0,\ F(+\infty) = 1$.

Зауваження. Неперервна випадкова величина X, що набуває значення в інтервалі (a, b), має незліченну кількість можливих значень, тому набуття X певних значень, наприклад, X = $a$ або X= $b$ – події з нульовою ймовірністю. Це означає, що P(X = $a$)= 0 та P(X = $b$)= 0. З огляду на це справедливі

P($a$ < X < $b$)= P($a$ ≤ X < $b$)= P( $a$ < X ≤ b)= P($a$ ≤ X ≤ b).

При дискретизації випадкових величин (округлення до якоїсь кількості десяткових знаків, чи обрізання наступних десяткових знаків після якогось), заміна ≤ на < неправомірна, й потрібно враховувати ймовірність.

Тоді, згідно з зауваженням, для такої неперервної випадкової величини $F(x) = 0$ за $x \leq a;\ \ F(x) = 1$ за $x \geq b$.

Графік її функції розподілу схематично наведено на рис.

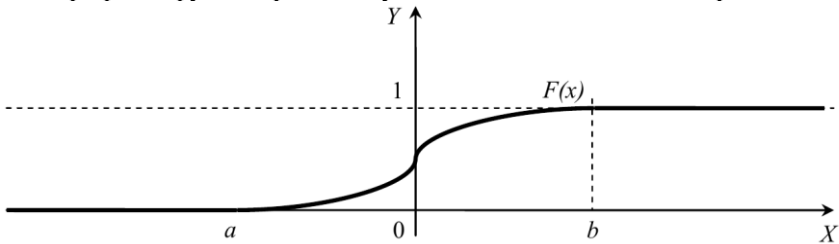

Закон розподілу ймовірностей неперервних випадкових величин може бути заданий також і щільністю розподілу.

**Означення.** Диференціальною функцією розподілу або щільністю розподілу ймовірностей неперервної випадкової величини називають
$$f(x) = F'(x)$$
Назва щільність ймовірностей випливає з рівності
$$f(x) = \lim_{\Delta x \to 0} \frac{P(X < x + \Delta x) - P(X < x)}{\Delta x}$$

**Властивості диференціальних функцій розподілу**
1. $f(x) \geq 0, x \in R$.

2. $\int\limits_{-\infty}^{+\infty} f(x)dx = 1.$

3. $P(a \leq X < b) = \int\limits_{a}^{b} f(x)dx.$

4. Якщо всі можливі значення випадкової величини належать інтервалу (a,b), то $\int_{a}^{b} f(x)dx = 1$

За відомою диференціальною функцією розподілу $f(x)$, знаходять інтегральну функцію розподілу

$$P(-\infty < X < x) = P(X < x) = F(x) = \int\limits_{-\infty}^{x} f(t)dt.$$

Геометричне тлумачення щільності розподілу випливає із формули

$$P(a \leq X < b) = \int\limits_{a}^{b} f(x)dx.$$

Ймовірність попадання випадкової величини X на проміжок [a, b) обчислюють як площу криволінійної трапеції, обмеженої зверху графіком функції $y = f(x)$, знизу – відрізком $[a, b]$ осі абсцис, зліва і справа – відрізками прямих $x = a$, $x = b$.

Якщо функція розподілу пропорційна якійсь невід'ємної функції $G(x)$, то вона дорівнює $f(x) = \frac{1}{C} G(x)$, де «нормуючий» коефіцієнт

$$C = \int\limits_{a}^{b} f(x)dx,$$

Де інтеграл беруть по усій області визначення від $a$ до $b$ – скінченній або нескінченній.

**Приклад** неперервних випадкових величин

**Розподіл Коші** має щільність розподілу

$$f(x) = \frac{1}{\pi(1 + x^2)}$$

та функцію розподілу

$$F(x) = \int_0^x \frac{dX}{\pi(1+X^2)} = \frac{\arctg(x)}{\pi} + \frac{1}{2}$$

Згенерувати випадкові числа із даним законом розподілу можна за формулою:

$$X = \tg\left(\pi(\text{rnd} - \frac{1}{2})\right)$$

Цей розподіл відповідає, напр., розподілу кількості фотонів на одиницю площі (пропорційній освітленості) при освітленні стіни в залежності від відстані від перпендикуляра від (маленького =»точкового») джерела випромінювання. Іншим прикладом можуть бути «крила» ліній у спектрі щільних газів.

Очевидно, що у реальних застосуваннях, величини є розмірними, і їх потрібно приводити до безрозмірних одиниць. Адже, неможливо обчислювати експоненти чи тригонометричні функції від (напр.) метрів чи кілограмів. Тому проводять так звану «нормалізацію» (лінійну заміну змінної) $x = (u-v)/w$, де $v$ — «нуль пункт» (початок відліку), а $w$ — масштаб (шкала, характерна ширина, одиниця вимірювання). Наприклад, у випадку освітлення стіни, природним масштабом $w$ є відстань джерела від стіни, а $(u-v)$ — відстань від перпендикуляра.

Загальний вигляд розподілу при такій лінійній заміні змінних,

$$f_u(u) = \frac{1}{w} f_x\left(\frac{u-v}{w}\right),$$

Адже елемент ймовірності

$$f_u(u)du = f_x(x)dx.$$

Функція розподілу та щільність ймовірності для розподілу Коші. Значення $x$, при яких $F(x) = \frac{P\%}{100\%}$, називають перцентілями (per cent = процент, відсоток). Медіаною називають 50% перцентіль, квартілями (чвертями) – 25%, 50%, 75%, отже, число попадає у 1,2,3 або 4-й квартілі. Децилями (десятковими частинами) – 10%, 20%,…90%. Отже, медіана є й другим квартилем, 5 децилем, 50 перцентилем.

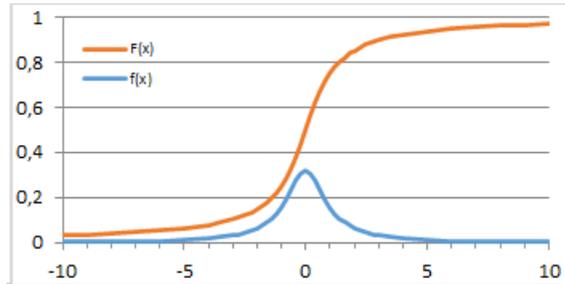

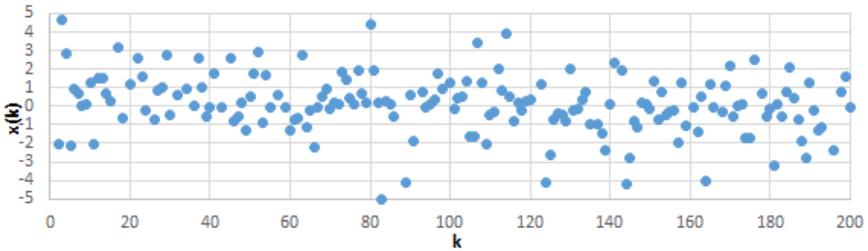

Одна з вибірок з 200 випадкових чисел із розподілом Коші.

**4.4. Числові характеристики випадкових величин.**

Під час розв'язування практичних задач досить важко, а іноді й неможливо визначити функцію розподілу випадкової величини X, тому потрібно вміти характеризувати її розподіл через параметри, найважливіші з яких – математичне сподівання та дисперсія.

**4.4.1. Математичне сподівання випадкової величини.**

**Математичним сподіванням дискретної випадкової величини X** називається сума добутків можливих значень цієї величини на ймовірність їх поява:

$$M(X) = \sum_{k=1}^{n} x_k p_k. \qquad (5)$$

Математичне сподівання вказує середнє значення, біля якого групуються всі можливі значення випадкової величини.

Механічний зміст математичного сподівання: коли на осі O$x$ розташовані точки $x_1, x_2, \ldots, x_n$ з масами $p_1, p_2, \ldots, p_n$, причому $\sum_{i=k}^{n} p_k = 1$, то M(X) – абсциса центру ваги даної системи матеріальних точок.

• **Приклад**. Дано ряд розподілу дискретної випадкової величини X:

| $X$ | 1 | 2 | 3 |
|---|---|---|---|
| $p$ | 0,4 | 0,5 | 0,1 |

Знайти: $M(X)$.

Відповідь: Візьмемо формулу $M(X) = \sum_{k=1}^{n} x_k p_k$ і підставимо табличні дані:

$$M(X) = 1 \cdot 0{,}4 + 2 \cdot 0{,}5 + 3 \cdot 0{,}1 = 1{,}7$$

**Математичне сподівання неперервної випадкової величини X** можна виразити не сумою добутків, а інтегралом
$$M(X) = \int_{-\infty}^{\infty} x \cdot f(x) dx. \tag{6}$$

Якщо всі можливі значення випадкової величини належать інтервалу $(a, b)$,

то математичне сподівання виражається інтегралом:
$$M(X) = \int_{a}^{b} x \cdot f(x) dx. \tag{7}$$

Механічна інтерпретація – така ж, як для дискретної випадкової величини – абсциса центру ваги. Але при цьому маса, яка дорівнює одиниці, розподілена на осі $Ox$ за законом щільності розподілу $f(x)$. Інше позначення математичного сподівання – $m_x$.

• **Приклад.** Випадкова величина задана щільністю розподілу:
$$f(x) = \begin{cases} 0, & x \leq 0 \\ 2x, & 0 < x \leq 1 \\ 0, & x > 1 \end{cases}$$

Знайти $M(X)$.

Відповідь: Формула для математичного сподівання неперервної випадкової величини **X**:
$$M(X) = \int_{-\infty}^{\infty} x \cdot f(x) dx.$$

Маємо три значення $f(x)$ на трьох інтервалах, тому М(Х) буде складатись із суми трьох інтегралів:
$$M(X) = \int_{-\infty}^{0} 0 \, dx + \int_{0}^{1} x \cdot 2x \, dx + \int_{1}^{\infty} 0 \, dx =$$

Інтеграл від 0 це стала, тому розглядаємо тільки другий визначений інтеграл
$$= \int_{0}^{1} x \cdot 2x \, dx = 2 \int_{0}^{1} x^2 \, dx =$$

Дивимось табличний інтеграл $\int x^n dx = \frac{x^{n+1}}{n+1}$, можемо записати:
$$= 2 \cdot \frac{x^3}{3} \Big|_{0}^{1} =$$

По формулі для визначеного інтеграла:

$$\int_a^b f(x)dx = F(X)\Big|_a^b = F(b) - F(a) \quad \text{знаходимо:}$$

$$= 2 \cdot \frac{1^3}{3} - 2 \cdot \frac{0^3}{3} = \frac{2}{3} - 0 = \frac{2}{3}$$

**Властивості математичного сподівання:**

- математичне сподівання постійної величини дорівнює самій постійній $M(C) = C$. При цьому постійна C – це випадкова величина, яка має одне значення C і набуває його з ймовірністю $p = 1$.

- постійний множник виноситься за знак математичного сподівання $M(C \cdot X) = C \cdot M(X)$. При цьому добуток $CX$ – випадкова величина, що набуває значення $Cx_1, Cx_2, ..., Cx_n$ з такою ж ймовірністю, як величина X – значення $x_1, x_2, ..., x_n$.

- випадкові величини називаються незалежними, якщо закон розподілу однієї з них не залежить від того, які можливі значення набуває інша величина. Інакше випадкові величини залежні. Кілька випадкових величин називають взаємно незалежними, якщо закони розподілу будь-якої кількості цих величин не залежать від того, яких можливих значень набули інші величини.

- Добуток незалежних випадкових величин XY – випадкова величина, можливі значення якої рівні добуткам кожного можливого значення X на кожне можливе значення Y, а ймовірності цих значень дорівнюють добуткам ймовірностей можливих значень співмножників. Якщо деякі добутки $x_i y_i$ рівні між собою, то їх імовірності додаються.

- Математичне сподівання добутку двох незалежних випадкових величин дорівнює добутку їхніх математичних сподівань:

$$M(X \cdot Y) = M(X) \cdot M(Y)$$

- сумою випадкових величин $X + Y$ називається випадкова величина, можливе значення якої дорівнює сумі кожного можливого значення X та кожного можливого значенням Y. а імовірності цих значень для незалежних величин X і Y дорівнюють добуткам ймовірностей доданків; для залежних величин – добуткам ймовірностей одного доданку та умовних ймовірностей другого. Якщо деякі суми $x_i + y_j$ рівні між собою, то їх імовірності додаються.

- Математичне сподівання суми двох випадкових величин дорівнює сумі математичних сподівань доданків:
$$M(X + Y) = M(X) + M(Y)$$
Наслідки:

а) $M(X_1 + X_2 + \ldots + X_n) = M(X_1) + M(X_2) + \ldots + M(X_n)$;

б) $M(X - Y) = M(X) - M(Y)$.

Властивості справедливі як для дискретних, так і для неперервних випадкових величин. Із властивостей математичного сподівання випливає теорема про математичне сподівання числа появи подій в незалежних випробуваннях.

**Теорема**: Якщо випадкова величина X число появ події A в *n* незалежних випробуваннях, при яких в кожному випробуванні $P(A) = p$, то математичне сподівання $M(X)$ дорівнює добутку числа випробувань на ймовірність появи події в кожному випробуванні:
$$M(X) = np.$$

### 4.4.2. Дисперсія випадкової величини.

Дисперсією випадкової величини називають математичне сподівання квадрата відхилення випадкової величини від її математичного сподівання:
$$D(X) = M(X - M(X))^2 \qquad (8)$$

Дисперсія характеризує ступінь розсіювання значень випадкової величини щодо її математичного сподівання.

**Теорема.** Дисперсія випадкової величини дорівнює різниці між математичним сподіванням квадрата цієї величини і квадратом її математичного сподівання:
$$D(X) = M(X^2) - (M(X))^2 \qquad (9)$$

За теоремою, зручно для обчислення дисперсії використовувати формули:

1) для дискретної випадкової величини:
$$D(X) = M(X^2) - (M(X))^2 = \sum_{i=1}^{n} x_i^2 p_i - m_x^2. \qquad (10)$$

2) для неперервної випадкової величини:
$$D(X) = \int_{-\infty}^{+\infty} x^2 f(x) dx - (M(X))^2; \qquad (11)$$

якщо її можливі значення належать відрізку [a, b], тоді:
$$D(X) = \int_a^b x^2 f(x)dx - (M(X))^2.$$
(12)

Також ці формули можуть бути записані інакше:
$$D(X) = \sum_{k=1}^n (x_k - m_x)^2 p_k$$
$$D(X) = \int_{-\infty}^{\infty} (x - m_x)^2 f(x)dx$$
$$D(X) = \int_a^b (x - m_x)^2 f(x)dx$$

Тут $m_x = M(X)$, $f(x)$ – щільність розподілу неперервної випадкової величини. Крім того, дисперсію часто позначають літерою S.

Розмірність дисперсії – квадрат розмірності випадкової величини, тому її не можна вказати на осі випадкової величини. Для наочності характеристики розсіювання зручніше користуватися величиною, розмірність якої збігається з розмірністю випадкової величини – коренем квадратним з дисперсії.

Корінь квадратний з дисперсії називається середнім квадратичним відхиленням і позначається:
$$\sigma(X) = \sqrt{D(X)}, \text{ або } \sigma_x = \sqrt{D_x}.$$
(13)

Відхилення випадкової величини від її математичного сподівання називається центрованою випадковою величиною, позначатимемо її як:
$$\hat{X} = X - M(X) \text{ або } \hat{X} = X - m_x.$$

Властивості дисперсії:
1. $D(C) = 0$, $C -$
2. $D(CX) = C^2 D(X)$
3. Дисперсія суми незалежних випадкових величин дорівнює сумі їх дисперсій:
$$D(X + Y) = D(X) + D(Y).$$

Наслідок. Дисперсія різниці двох незалежних випадкових величин дорівнює сумі їх дисперсій:
$$D(X - Y) = D(X) + D(Y)$$

Зауваження. Крім середньоквадратичного відхилення, характеристиками рзкиду вибіркового ряду є амплітуда ($X_{max} - X_{min}$), міжквартильна різниця ($X(75\%) - X(25\%)$), довірчий інтервал ($X(100\% - P\%/2) - X(P\%/2)$), де $P\% -$

ймовірність виходу випадкової віеличини за довірчий інтервал

(по $P\%/2$ з обох сторін).

• **Приклад 1:** За законом розподілу дискретної випадкової величини:

| $X$ | 1 | 2 | 3 |
|---|---|---|---|
| $P$ | 0,4 | 0,5 | 0,1 |

знайти $D(X)$, $\sigma(X)$.

Відповідь: Запишемо формулу для пошуку дисперсії:
$$D(X) = M(X^2) - \big(M(X)\big)^2$$
Тобто, спочатку треба визначити математичне сподівання за формулою:
$$M(X) = \sum_{k=1}^{n} x_k p_k$$
таким чином маємо записати:
$$M(X) = 1 \cdot 0,4 + 2 \cdot 0,5 + 3 \cdot 0,1 = 1,7.$$
Тепер треба порахувати $M(X^2)$, тобто, замість значення X підставимо у туж саму формулу, але $X^2$ і отримаємо:
$$M(X^2) = 1 \cdot 0,4 + 4 \cdot 0,5 + 9 \cdot 0,1 = 3,3.$$
Підставимо отримані дані в формулу дисперсії і отримаємо:
$$D(X) = M(X^2) - (M(X))^2 = 3,3 - 1,7^2 = 0,41,$$
Тепер можемо визначити середнє квадратичне відхилення:
$$\sigma(X) = \sqrt{D(X)} = \sqrt{0,41} = 0,64.$$

• **Приклад 2**. Є можливість вибрати спосіб виробництва і реалізації двох наборів товарів широкого вжитку. За даними відділу маркетингу, яким були проведені дослідження ринку, можливий прибуток від виробництва і реалізації X і Y наведено в таблицях (X, Y – прибуток у грошових одиницях):

| $X$ | 1000 | 1500 | 2000 |
|---|---|---|---|
| $P$ | 0,5 | 0,3 | 0,2 |

та

| $Y$ | 1000 | 1500 | 1750 |
|---|---|---|---|
| $P$ | 0,4 | 0,4 | 0,2 |

Потрібно оцінити ступінь ризику і запропонувати рішення щодо випуску і реалізації одного із наборів товарів.

Відповідь: Спочатку знайдемо математичне сподівання можливого прибутку по формулі $M(X) = \sum_{k=1}^{n} x_k p_k$ для кожного варіанта:

$$M(X) = 1000 \cdot 0,5 + 1500 \cdot 0,3 + 2000 \cdot 0,2 = 1350;$$

$$M(Y) = 1000 \cdot 0,4 + 1500 \cdot 0,4 + 1750 \cdot 0,2 = 1350.$$

Як ми бачимо, обидва варіанти мають однакове математичне сподівання можливого прибутку. Тому оцінюємо ступінь ризику кожного варіанта, показником чого буде дисперсія, тобто, оцінемо ступінь розсіювання значень величини прибутку щодо математичного сподівання:

$$D(X) = M(X^2) - \bigl(M(X)\bigr)^2$$

Для цього замість X підставимо відхилення від середнього, яким і є пораховане математичне сподівання, тобто:

$$D(X) = (1000 - 1350)^2 \cdot 0,5 + (1500 - 1350)^2 \cdot 0,3 + (2000 - 1350)^2 \cdot 0,2 =$$

$$= 61250 + 6750 + 84500 = 152500, \sigma(X) = 390,5.$$

$$D(Y) = (1000 - 1350)^2 \cdot 0,4 + (1500 - 1350)^2 \cdot 0,4 + (1750 - 1350)^2 \cdot 0,2 =$$

$$= 49000 + 9000 + 32000 = 90000, \sigma(Y) = 300.$$

Ступінь ризику, пов'язаний з виробництвом і реалізацією набору X більший, ніж набору Y. Варіант Y менш ризикований.

• **Приклад 3.** Задано ряд розподілу випадкової величини X. Знайти функцію розподілу X та обчислити математичне сподівання і дисперсію.

| X | 1 | 2 | 3 | 4 | 5 |
|---|---|---|---|---|---|
| P | 0,200 | 0,111 | 0,125 | 0,143 | 0,421 |

Відповідь:
1) Знайти функцію розподілу X.

Якщо заданий ряд розподілу $(X_i, P_i), i = 1 \ldots n$, то функцію розподілу $F(x)$ випадкової величини X знайдемо за формулою:
$F(x) = F(X \leq x) = \sum_{i; x_i < x} P_i,$

$$F(x_i) = \begin{cases} 0, x_1 < x \\ P_1, x_1 \leq x < x_2 \\ P_1 + P_2, x_2 \leq x < x_3 \\ P_1 + P_2 + P_3, x_3 \leq x < x_4 \\ P_1 + P_2 + P_3 + P_4, x_4 \leq x < x_5 \\ P_1 + P_2 + P_3 + P_4 + P_5, x_5 \leq x \leq \infty \end{cases}$$

Сума усіх ймовірностей має дорівнювати одиниці для будь-якого варіанту (у розглянутому, лише 5):

$$F(x_5) = \sum_{i=1}^{n=5} P_i = P_1 + P_2 + P_3 + P_4 + P_5 = 1$$

Можна порахувати значення $F_i = F(x_i)$ і побудувати графік функції розподілу:

$F_1 = P_1 = 0{,}200$
$F_2 = P_1 + P_2 = 0{,}200 + 0{,}111 = 0{,}311$
$F_3 = P_1 + P_2 + P_3 = 0{,}311 + 0{,}125 = 0{,}436$
$F_4 = P_1 + P_2 + P_3 + P_4 = 0{,}436 + 0{,}143 = 0{,}579$
$F_5 = P_1 + P_2 + P_3 + P_4 + P_5 = 0{,}579 + 0{,}421 = 1{,}000$

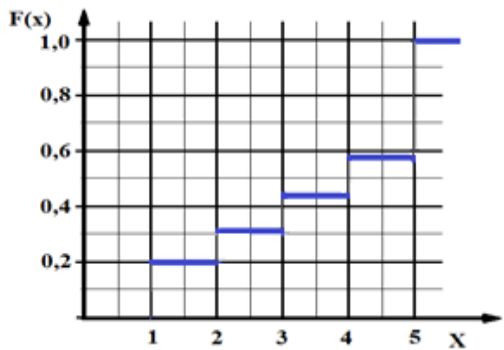

2) Обчислити математичне сподівання:

Визначимо математичне сподівання за формулою:

$$M(X) = \sum_{i=1}^{n} x_i p_i$$

$M(X) = 1 \cdot 0{,}2 + 2 \cdot 0{,}111 + 3 \cdot 0{,}125 + 4 \cdot 0{,}143 + 5 \cdot 0{,}421 =$
$0{,}2 + 0{,}222 + 0{,}375 + 0{,}572 + 2{,}105 = 3{,}474$

3) Обчислити дисперсію:

Запишемо формулу для пошуку дисперсії:

$$D(X) = M(X^2) - (M(X))^2$$
$$M(X^2) = 1^2 \cdot 0{,}2 + 2^2 \cdot 0{,}111 + 3^2 \cdot 0{,}125 + 4^2 \cdot 0{,}143 + 5^2 \cdot 0{,}421 =$$
$$0{,}2 + 0{,}444 + 1{,}125 + 2{,}288 + 10{,}525 = 14{,}582$$
$$D(X) = 14{,}582 - 3{,}474^2 = 14{,}582 - 12{,}069 = 2{,}513$$

Тепер можемо визначити середнє квадратичне відхилення:
$$\sigma(X) = \sqrt{D(X)} = \sqrt{2{,}513} = 1{,}585$$

• **Приклад 4.** Знайти числові характеристики випадкової величини X, яка задана функцією розподілу:
$$F(x) = \begin{cases} 0, & x \leq 0; \\ \dfrac{x^2}{25}, & 0 < x \leq 5; \\ 1, & x > 5. \end{cases}$$

Відповідь: Щоб знайти математичне сподівання за формулою:
$$M(X) = \int_a^b x f(x) dx$$

треба спочатку знайти диференціальну функцію розподілу $f(x)$, тобто, щільність розподілу ймовірності за формулою $f(x) = F'(x)$:

(пам'ятаємо, що похідна від 0 це 0, від $1' = 0$, $\left(\dfrac{x^2}{25}\right)' = \dfrac{2x}{25}$)

$$f(x) = \frac{2x}{25}, \quad 0 < x \leq 5$$

За формулою $M(X) = \int_a^b x f(x) dx$ знаходимо математичне сподівання:

$$M(X) = \int_0^5 x \cdot \frac{2}{25} x \, dx = \frac{2}{25} \cdot \frac{x^3}{3} \Big|_0^5 = \frac{2}{75}(5^3 - 0) = \frac{10}{3}.$$

Для дисперсії використовуємо формулу $D(X) = \int_a^b x^2 f(x) dx - (M(X))^2$ :

$$D(X) = \int_0^5 x^2 \cdot \frac{2}{25} x \, dx - \left(\frac{10}{3}\right)^2 = \frac{2}{25} \cdot \frac{x^4}{4} \Big|_0^5 - \frac{100}{9} = \frac{25}{18}.$$

$x^2$ тому що з умови задачі $M(X)$ це $M(x^2)$, а $f(x) = \dfrac{2x}{25}$ на

інтервалі $0 < x \leq 5$,

який встановлює межі інтегрування.

Середнє квадратичне відхилення $\sigma(X) = \sqrt{D(X)} = \sqrt{\frac{25}{18}} = \frac{5}{3\sqrt{2}} = \frac{5\sqrt{2}}{6} \approx 1{,}17$.

• **Приклад 5.** Дана щільність розподілу $p(x)$ випадкової величини X. Знайти параметр $\gamma$, математичне сподівання, дисперсію, функцію розподілу випадкової величини X, ймовірність виконання нерівності $x_1 < X < x_2$.

1) $p(x) = \begin{cases} \dfrac{1}{\gamma - a}, x \in [a, b] \\ 0, x \notin [a, b] \end{cases}$ 

2) $p(x) = \begin{cases} a, x \in [\gamma, b] \\ 0, x \notin [\gamma, b] \end{cases}$

3) $p(x) = \begin{cases} \gamma, x \in [a, b] \\ 0, x \notin [a, b] \end{cases}$ 

4) $p(x) = \begin{cases} a, x \in \left[\dfrac{b - \gamma}{2}, \dfrac{b + \gamma}{2}\right] \\ 0, x \notin \left[\dfrac{b - \gamma}{2}, \dfrac{b + \gamma}{2}\right] \end{cases}$

| $a$ | $b$ | $x_1$ | $x_2$ |
|---|---|---|---|
| 2,5 | 4 | 3 | 3,3 |

Відповідь:

1. Знайти параметр $\gamma$.

Раніше надавався приклад неперервної випадкової величини:

Випадкова величина X розподілена рівномірно на проміжку $(a, b)$, тобто, усі її можливі значення належать цьому проміжку і щільність розподілу її стала на цьому проміжку, тобто:

$$f(x) = \begin{cases} \dfrac{1}{b - a}, x \in [a, b] \\ 0, \quad x \notin [a, b] \end{cases}$$

Графік щільності рівномірного розподілу зображено на рис.

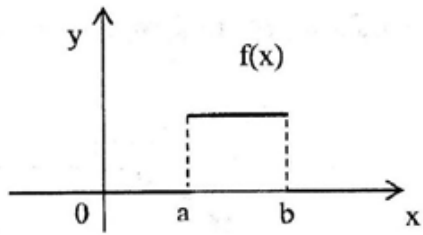

Оскільки частина розподілу випадкової величини постійна в межах деякого інтервалу, скористаймось формулами для звичного випадку

$$p(x) = \begin{cases} 0, & x < A \\ \dfrac{1}{B-A}, & A \leq x \leq B \\ 0, & x > B \end{cases}$$

$p(x)$ – позначення щільності в умовах цього завдання; $A$ і $B$ – значення, які треба визначити через $\gamma$ та $a$.

Використаємо змінні $A, B$, які визначаються порівнянням з наведеними вище формулами окремо для кожного варіанту 1) – 4).

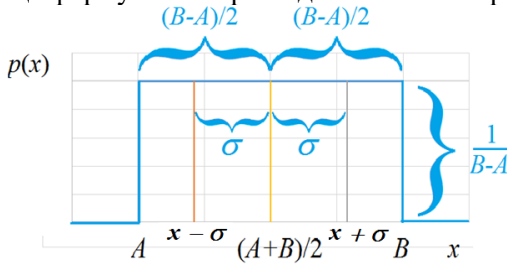

Визначимо параметр $\gamma$, для чого порівняємо два вирази - наданий в умовах задачі (варіант 1)) та для формули звичного вигляду:

$$\frac{1}{\gamma - a} = \frac{1}{B - A}$$

$\gamma - a = B - A$,   тобто, $\gamma = B = 4, a = A = 2{,}5$

Якщо щільність розподілу $p(x)$ випадкової величини X має наступний вигляд (2) варіант):

$$p(x) = \begin{cases} a, x \in [\gamma, b] \\ 0, x \notin [\gamma, b] \end{cases}$$

тоді рівняємо з:

$$p(x) = \begin{cases} 0, & x < A \\ \dfrac{1}{B-A}, & A \leq x \leq B \\ 0, & x > B \end{cases}, \text{тобто}$$

На рисунку для звичного вигляду інтервал позначений, як $A \leq x \leq B$, а в задачі для 1) варіанту $x \in [\gamma, b]$ $a = \frac{1}{B-A}$

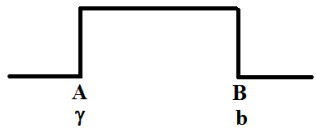

Звідси можемо перепозначити:

$$a = \frac{1}{B-A} = \frac{1}{b-\gamma}$$
$$b - \gamma = \frac{1}{a}; \quad \gamma = b - \frac{1}{a}$$

Підставимо значення $a$ та $b$ з умов задачі і отримаємо значення $\gamma$.

Для 3 варіанту маємо:
$$p(x) = \begin{cases} \gamma, x \in [a,b] \\ 0, x \notin [a,b] \end{cases}$$

Рівняємо з:
$$p(x) = \begin{cases} 0, x < A \\ \frac{1}{B-A}, A \leq x \leq B, \\ 0, x > B \end{cases} \Rightarrow \gamma = \frac{1}{B-A} = \frac{1}{b-a}$$

Для 4) варіанту маємо:
$$p(x) = \begin{cases} a, x \in \left[\frac{b-\gamma}{2}, \frac{b+\gamma}{2}\right] \\ 0, x \notin \left[\frac{b-\gamma}{2}, \frac{b+\gamma}{2}\right] \end{cases}$$

Рівняємо з:
$$p(x) = \begin{cases} 0, & x < A \\ \frac{1}{B-A}, & A \leq x \leq B, \\ 0, & x > B \end{cases} \Rightarrow a = \frac{1}{B-A}$$

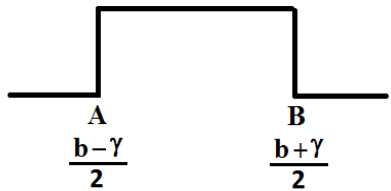

$$A = \frac{b-\gamma}{2}; \; B = \frac{b+\gamma}{2}; \; a = \frac{1}{B-A} \Rightarrow$$
$$a = \frac{1}{\frac{b+\gamma}{2} - \frac{b-\gamma}{2}} = \frac{1}{\frac{2\gamma}{2}} = \frac{1}{\gamma} \Rightarrow \gamma = \frac{1}{a}$$

2. Знайти математичне сподівання.

Математичне сподівання вказує середнє значення, біля якого групуються всі можливі значення випадкової величини, тому (для 1 варіанту) маємо математичне сподівання випадкової величини
$$\overline{x} = \frac{A+B}{2}, \ \overline{x} = \frac{a+b}{2} = \frac{2{,}5+4}{2} = 3{,}25$$

Також, за визначенням математичного сподівання функції $G(x)$ випадкової величини $x$ маємо:

$$M(G(x)) = \int_{-\infty}^{\infty} G(x)p(x)\,dx, \quad \text{чи} \quad M(X) = \overline{x} = \int_{-\infty}^{\infty} xp(x)\,dx.$$

Із врахуванням того, що $p(x)$ визначається різними формулами у трьох інтервалах, цей інтеграл можна розбити на три:

$$\overline{x} = \int_{-\infty}^{\infty} xp(x)\,dx = \int_{-\infty}^{A} x \cdot 0\,dx + \int_{A}^{B} x \cdot p(x)\,dx + \int_{B}^{\infty} x \cdot 0\,dx = \int_{A}^{B} x \cdot p(x)\,dx.$$

$$p(x) = \begin{cases} 0, & x < A \\ \dfrac{1}{B-A}, & A \leq x \leq B \\ 0, & x > B \end{cases}$$

Оскільки $p(x) = 1/(B-A)$ для рівномірного розподілу, то
$$M(X) = \overline{x} = \int_{A}^{B} x \cdot p(x)\,dx = \int_{A}^{B} x \cdot \frac{1}{(B-A)}\,dx = \frac{1}{(B-A)} \int_{A}^{B} x\,dx =$$
$$\frac{1}{(B-A)} \frac{x^2}{2}\Big|_{A}^{B} = \frac{1}{(B-A)} \frac{B^2}{2} - \frac{1}{(B-A)} \frac{A^2}{2} = \frac{B^2-A^2}{2(B-A)} = \frac{(B-A)(B+A)}{2(B-A)} = \frac{B+A}{2} =$$
$$\frac{A+B}{2} = 3{,}25$$

Останній вигляд – лише дань історичній традиції, коли літери, якщо можливо, у алфавітному порядку.

*Примітка. Цей результат є дійсним не лише для рівномірного розподілу, а й для будь-якого іншого за умови симетрії щільності розподілу:*

$$p\left(\frac{A+B}{2} + u\right) = p\left(\frac{A+B}{2} - u\right)$$

*для будь-якого $u$ (різниця між поточним значенням аргументу та середнім значенням, де середнє значення співпадає з положенням точки симетрії.*

3. Обчислимо дисперсію $D(X) = \sigma^2$.
За визначенням:
$$D(X) = \sigma^2 = \overline{(x-\overline{x})^2} = \overline{x^2} - \overline{x}^2 = \overline{(x-u)^2} - (u-\overline{x})^2:$$
«дисперсія дорівнює середньому квадрата мінус квадрат середнього» (при $u = 0$). Втім, можна рекомендувати використовувати $u$, яке буде зручним (напр., для років народження людей у групі, брати не від «Різдва Христова», а від «Міленіума»)

Обчислимо $\overline{x^2}$, оскільки:

$$M(X) = \overline{x} = \int_A^B x \cdot p(x)\, dx \quad \text{та} \quad p(x) = \frac{1}{B-A}, \quad \text{то:}$$

$$\overline{x^2} = \int_A^B x^2 \cdot p(x)\, dx = \int_A^B x^2 \cdot \frac{1}{(B-A)}\, dx = \frac{1}{(B-A)} \int_A^B x^2\, dx =$$

$$\frac{1}{(B-A)} \left.\frac{x^3}{3}\right|_A^B = \frac{1}{(B-A)} \cdot \frac{B^3}{3} - \frac{1}{(B-A)} \cdot \frac{A^3}{3} = \frac{B^3 - A^3}{3(B-A)} =$$

$$= \frac{(B-A)(B^2 + BA + A^2)}{3(B-A)} = \frac{B^2 + BA + A^2}{3}$$

$$D(X) = \sigma^2 = \overline{x^2} - \overline{x}^2 =$$

$$= \frac{B^2 + BA + A^2}{3} - \left(\frac{B+A}{2}\right)^2 = \frac{B^2 + BA + A^2}{3} - \frac{B^2 + 2BA + A^2}{4} =$$

$$\frac{4 \cdot (B^2 + BA + A^2)}{4 \cdot 3} - \frac{3 \cdot (B^2 + 2BA + A^2)}{3 \cdot 4} = \frac{B^2 - 2BA + A^2}{12} =$$

$$= \frac{(B-A)^2}{12} == \frac{(4-2{,}5)^2}{12} = \frac{2{,}25}{12} \approx 0{,}19$$

Відповідно, середнє квадратичне відхилення від середнього

$$\boldsymbol{\sigma} = \sqrt{\sigma^2} = \sqrt{\frac{(B-A)^2}{12}} = \frac{B-A}{2\sqrt{3}} \approx 0{,}44$$

4. Знайти функцію розподілу випадкової величини X.

Випадкову величину X називають неперервною, якщо її функцію розподілу можна подати у вигляді: $F(x) = \int_{-\infty}^{x} f(t) dt,$ де $f(x)$ деяка функція, яку називають щільністю розподілу ймовірностей.

Для наших позначень - випадкову величину X називають неперервною, якщо її функцію розподілу можна подати у вигляді:

$$F(X) = \int_{-\infty}^{X} p(x)\, dx$$

де $p(x) \geq 0$ деяка функція, яку називають щільністю (або «густиною») розподілу ймовірностей. Іноді вживають позначення $f(x)$.

Тобто, щоб записати функцію розподілу, маємо знайти: інтеграл. Раніше ми вже з'ясували, що для варіанту $p(x) = \frac{1}{B-A}$. Для рівномірного розподілу,

$$F(X) = \int_A^X \frac{1}{B-A} dx = \frac{x}{B-A}\Big|_A^X = \frac{X}{B-A} - \frac{A}{B-A} = \frac{X-A}{B-A} = \frac{X-2,5}{1,5}$$

$$F(X) = \int_{-\infty}^X p(x)dx = \int_{-\infty}^A p(x)dx + \int_A^X p(x)dx$$

Знайти ймовірність виконання нерівності $x_1 \leq X \leq x_2$, для варіанту $3 \leq X \leq 3,3$,

Ймовірність виконання нерівності $P(x_1 \leq X \leq x_2) = F(x_2) - F(x_1)$, де $F(x)$ обчислюється за різними формулами для різних інтервалів.

За умовою задачі: $a = 2,5; b = 4; x_1 = 3; x_2 = 3,3$. Підставимо $x_1$ і $x_2$ в

$$F(X) = \frac{X-A}{B-A} = \frac{X-2,5}{1,5}$$

$$F(x_1) = \frac{3-2,5}{1,5} = \frac{0,5}{1,5} = 0,333$$

$$F(x_2) = \frac{3,3-2,5}{1,5} = \frac{0,8}{1,5} = 0,533$$

$$P(x_1 \leq X \leq x_2) = F(x_2) - F(x_1) = 0,533 - 0,333 = 0,2$$

Іншими словами, для рівномірного розподілу, функція розподілу $F(x)$ виглядає, як горизонтальна піввісь на рівні нуля зліва від лівої границі інтервалу $[A, B]$ та рівні одиниці справа від інтервалу. Всередині інтервалу, це відрізок прямої від $(A, 0)$ до $(B, 1)$.

Наприклад, на рис. показана функція розподілу для $A = 2$, $B = 6$,

$x_1 = 2.5$, $x_2 = 4.5$, Червоний відрізок з'єднує точки із координатами $(x|1|, F(x_1)), (x|2|, F(x_2))$. Різниця по висоті й дає ймовірність потрапляння у інтервал $0.625 - 0.125 = 0.500$. Найчастіше помилки виникають, коли середня формула використовується для точок за межами інтервалу $[A, B]$.

Наприклад, для $F(x_1 = 0$ (оскільки $x_1 = 1 < A = 2$, а за формулою для середнього інтервалу, $\frac{1-2}{6-2} = \frac{-1}{4} = -0.125$. Це не має сенсу, адже будь-яка ймовірність, включаючи функцію розподілу, не може бути меншою за 0 (нейможлива подія) та більшою за одиницю (достовірна подія).

$A < x_1 = 2.5 > B$

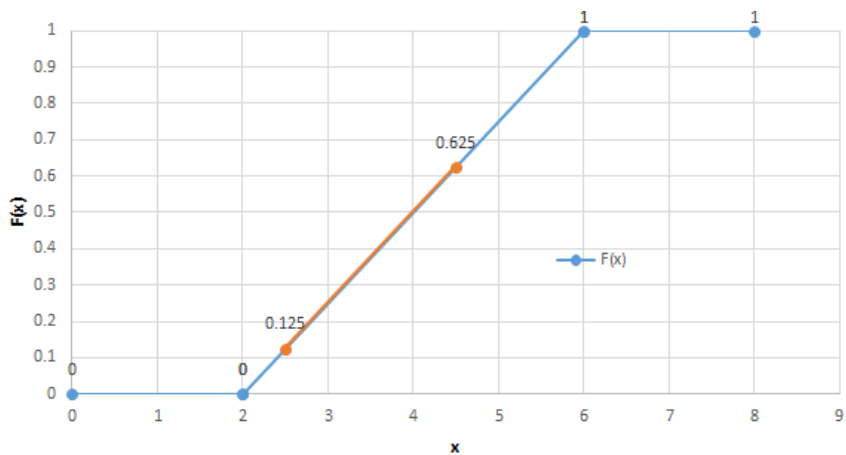

$x_1 = 1 < A$

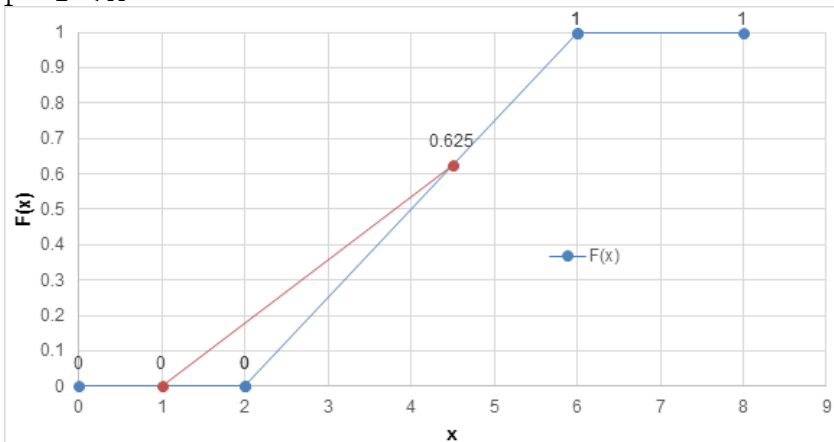

• **Приклад 6.** Щільність розподілу ймовірностей випадкової величини X має вигляд :
$$p(x) = \gamma e^{ax^2+bx+c}.$$

Знайти параметр $\gamma$, математичне сподівання, дисперсію, функцію розподілу випадкової величини X, ймовірність виконання нерівності $x_1 < X < x_2$.

Варіант 1.

| $a$ | $b$ | $c$ | $x_1$ | $x_2$ |
|---|---|---|---|---|
| $-2$ | 8 | -2 | 1 | 3 |

Відповідь:
1) Знайти параметр $\gamma$

Рішення визначається порівнянням з стандартною формою щільності нормального розподілу (Гауса):

$$f(x) = p(x) = \frac{1}{\sqrt{2\pi}\sigma} \exp\left(\frac{-(x-\bar{x})^2}{2\sigma^2}\right)$$

Графік цієї функції $f(x)$ чи в нашому випадку $p(x)$ має максимум $\frac{1}{\sigma\sqrt{2\pi}}$ при $x = \bar{x}$, симетричний відносно прямої $x = \bar{x}$, і має вісь 0X асимптотою. Графік щільності нормального розподілу (з параметрами $\bar{x} = -0.5$, $\bar{\sigma}^2 = 1/6$) має вигляд:

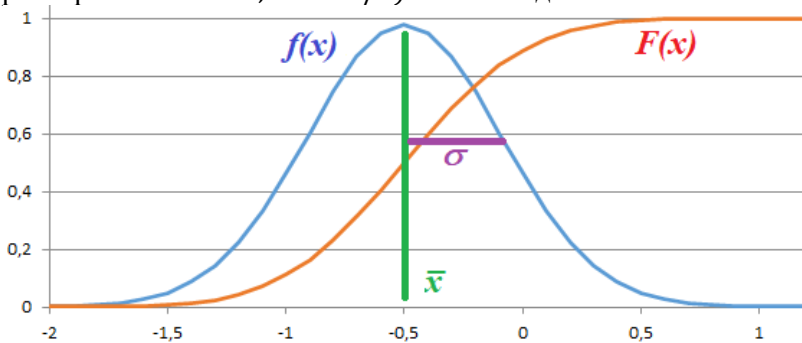

З графіку видно, що $\bar{x}$ являється середнім значенням $x$.

$$p(x) = \frac{1}{\sqrt{2\pi}\sigma} \exp\left(\frac{-(x-\bar{x})^2}{2\sigma^2}\right)$$

Звідки за формулою $(a-b)^2 = a^2 - 2ab + b^2$ перетворимо вираз:

$$p(x) = \frac{1}{\sqrt{2\pi}\sigma} \exp\left(-\frac{x^2}{2\sigma^2} + \frac{x\bar{x}}{\sigma^2} - \frac{\bar{x}^2}{2\sigma^2}\right) =$$
$$\frac{1}{\sqrt{2\pi}\sigma} \exp\left(\frac{-x^2}{2\sigma^2}\right) \exp\left(\frac{\bar{x}x}{\sigma^2}\right) \exp\left(-\frac{\bar{x}^2}{2\sigma^2}\right)$$

Порівняємо з виразом наданим в умовах задачі: $p(x) = \gamma e^{ax^2+bx+c}$. Далі $a$ − параметр, заданий в умові, хоча, у попередніх формулах, мав інші сенси.

$$\gamma e^{ax^2+bx+c} = \gamma e^{ax^2} e^{bx} e^c$$

Звідки $a = -1/2\,\sigma^2$, тобто, $\sigma = \sqrt{-1/2\,a}$ – це і є $\sqrt{D(X)}$, $b = \bar{x}/\sigma^2 = -2a\bar{x}$, тобто, $\bar{x} = -b/2\,a$ – це є математичне сподівання М(X). Приклад щільності та функції розподілу для

$\bar{x} = -0.5, \bar{\sigma}^2 = 1/6, \sigma = 0.4082483$.

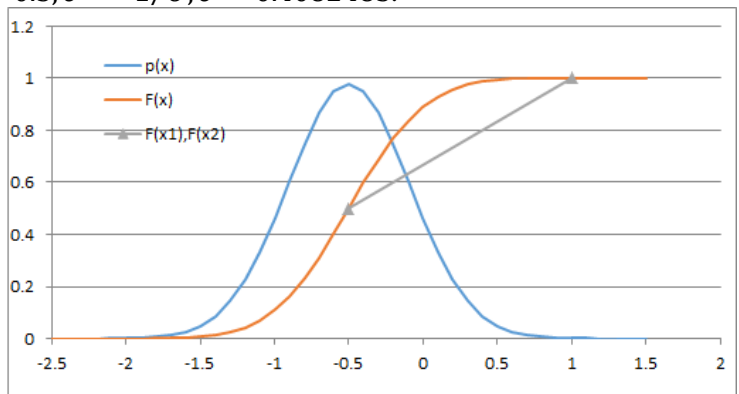

Щоб знайти параметр $\gamma$ розглянемо формулу:
$$\frac{1}{\sqrt{2\pi}\sigma} e^{\frac{-x^2}{2\sigma^2}} e^{\frac{\bar{x}x}{\sigma^2}} e^{\frac{-\bar{x}^2}{2\sigma^2}} = \gamma e^{ax^2} e^{bx} e^c$$

Скорочуючи на експоненти, які залежать від $x$, отримуємо
$$\frac{1}{\sqrt{2\pi}\sigma} e^{\frac{-\bar{x}^2}{2\sigma^2}} = \gamma e^c,$$

$\gamma = \frac{1}{\sqrt{2\pi^2\sigma^2}} e^{-c} \exp\left(-\frac{\bar{x}^2}{2\sigma^2}\right),$ підставимо отримане $\sigma = \sqrt{-1/(2a)}$

$$\frac{\sqrt{-2a}}{\sqrt{2\pi}} e^{-c} \exp(a\bar{x}^2).$$

Тепер треба підставити надані значення із умов задачі, наприклад, для 1 варіанту $a = -2, b = 8, c = -2$
$$\gamma = \sqrt{\frac{2}{\pi}} e^{-6}$$

Для визначення ймовірності попадання випадкової величини у діапазон $X \in [x_1, x_2]$ (або $x_1 \leq X \leq x_2$), потрібно скористатись формулою
$$P(x_1, x_2) = \int_{x_1}^{x_2} p(x)dx = F(x_2) - F(x_1),$$
$$F(x_1) = \int_{-\infty}^{x_1} p(x)dx \text{ та } F(x_2) = \int_{-\infty}^{x_2} p(x)dx$$

$$p(x) = \frac{1}{\sqrt{2\pi}\sigma} \exp\left(\frac{-(x-\bar{x})^2}{2\sigma^2}\right)$$

Для інтегрування, робимо заміну змінних $t = \frac{x-\bar{x}}{\sigma}$, де $\bar{x}$ та $\sigma$ вже визначені вище і тепер це є постійні, тоді:

$$dt = \left(\frac{x}{\sigma}\right)' dx - \left(\frac{\bar{x}}{\sigma}\right)' dx = \frac{dx}{\sigma}, \text{тому що } x' = 1, \text{а} \left(\frac{\bar{x}}{\sigma}\right)' = 0$$

Підставимо в

$$F(x_1) = \int_{-\infty}^{x_1} p(x)dx = \int_{-\infty}^{\frac{x_1-\bar{x}}{\sigma}} \frac{1}{\sqrt{2\pi}\sigma} e^{\frac{-t^2}{2}} \cdot \sigma dt = \int_{-\infty}^{\frac{x_1-\bar{x}}{\sigma}} \frac{1}{\sqrt{2\pi}} e^{\frac{-t^2}{2}} dt = F_t(t_1)$$

$M(t) = \bar{t} = 0, D(t) = \sigma_t^2 = 1$ (для всіх варіантів), отже,

$$p_t(t) = \frac{1}{\sqrt{2\pi}} e^{\frac{-t^2}{2}}$$

$$F(x_1) = \int_{-\infty}^{x_1} p(x)dx = \int_{-\infty}^{t_1} p_t(t)dt = F_t(t_1)$$

$$F(x_2) = \int_{-\infty}^{x_2} p(x)dx = \int_{-\infty}^{t_2} p_t(t)dt = F_t(t_2)$$

У електронних таблицях (Libre/Open Office Calc, MicroSoft Office Excel, та ін.) та програмах статистичного призначення, є можливість обчислити $F(x)$, тому що інтеграл не можна визначити через «елементарні функції», а, чисельно, через (теоретично, нескінченні) ряди. Для зручності, у докомп`ютерну еру, були розраховані таблиці значень цієї функції. Для інтегрування було зроблено заміну змінної $t = \frac{x-\bar{x}}{\sigma}$. Для цієї змінної, яку називають «нормалізованою», знаходимо значення, які відповідають $x_1, x_2$: $t_1 = \frac{x_1-\bar{x}}{\sigma}, t_2 = \frac{x_2-\bar{x}}{\sigma}$.

Інтеграл від функції $p_t(t)$ дорівнює $F_t(t) = F(x)$. Оскільки функція парна, то $p_t(-t) = p_t(t)$, та $F_t(-t) = 1 - F_t(t)$. Історично, була введена «функція Лапласа» (інша назва «інтеграл ймовірностей»)

$$\Phi(t) = F_t(t) - 0.5 = \int_0^t p_t(T)dT, \text{ тоді } \Phi(-t) = -\Phi(t).$$

Для нашого варіанту маємо $x_1 = 1$ та $x_2 = 3$. Раніше було визначено

$$\bar{x} = \frac{-b}{2a} = \frac{-8}{2(-2)} = 2 \text{ та } \sigma = \sqrt{-1/2\,a} = \sqrt{-1/2\,(-2)} = \sqrt{0{,}25} = \frac{1}{2}$$

Визначимо $t_1 = \frac{x_1 - \bar{x}}{\sigma} = \frac{1-2}{0{,}5} = -2$; $t_2 = \frac{x_2 - \bar{x}}{\sigma} = \frac{3-2}{0{,}5} = 2$.:

По таблиці знаходимо значення функцій $\Phi(t_1) = -0.477$ та $\Phi(t_2) = 0.477$

$P(x_1 \leq X \leq x_2) = F(x_2) - F(x_1) = \Phi(t_2) - \Phi(t_1) = 0.477 + 0.477 = 0.954$

Таблиці її значень опубліковані у багатьох підручниках, довідниках, зокрема,

https://web.posibnyky.vntu.edu.ua/fitki/4tichinska_teoriya_jmovirnostej/da.htm

Також таблицю приводимо у додатку 1. Напр., для $t - 0{,}24$. Потрібно взяти число на перетині 0,2 та 0,04, тобто, 0,094835. На практиці, рідко використовують таку кількість знаків, звичайно, вистачає 4 знаки після коми чи точки після цілого нуля. Крім випадків, коли результуюча ймовірність дуже мала. Втім, ми дали додатково 2 знаки (із точністю до міліонних частин).

**4.4.3. Початкові та центральні моменти випадкової величини $X$.**

Розглянемо однаково розподілені взаємно незалежні випадкові величини .

Нехай $X_1, X_2, ..., X_n$ – $n$ взаємно незалежних однаково розподілених випадкових величин, які мають однакові характеристики. Позначимо через $a$

і D відповідно їх математичне сподівання і дисперсію, $\sigma = \sqrt{D}$. Нехай

$$\bar{X} = \frac{X_1 + X_2 + \cdots + X_n}{n}$$

– це середнє арифметичне цих величин.

Тоді числові характеристики $\bar{X}$ такі:

1. $M(\bar{X}) = a$.

2. $D(\bar{X}) = \dfrac{D}{n}$.

3. $\sigma(\bar{X}) = \dfrac{\sigma}{\sqrt{n}}$.

З огляду на ці результати середнє арифметичне досить великої кількості взаємно незалежних випадкових величин має менше розсіювання, ніж кожна окрема величина. Тому середнє арифметичне кількох вимірювань беруть за наближене значення вимірюваної величини.

**Означення. Початковим моментом** порядку k (k – натуральне число) випадкової величини $X$ називають математичне сподівання величини $X^k$:

$$\nu_k = M(X^k).$$

Зокрема, $\nu_1 = M(X)$, $\nu_2 = M(X^2)$, тоді формулу
$$D(X) = M(X^2) - (M(X))^2.$$
для обчислення дисперсії можна записати як:
$$D(X) = \nu_2 - \nu_1^2.$$

**Означення**. **Центральним моментом** порядку $k$ ($k$ – натуральне число) випадкової величини X називають математичне сподівання величини $(X - M(X))^k$:

$$\mu_k = M(X - M(X))^k.$$

Зокрема, $\mu_1 = M(X - M(X)) = 0$, $\mu_2 = M(X - M(X))^2 = D(X)$, відповідно порівнюючи формули дисперсії маємо:

$$\mu_2 = \nu_2 - \nu_1^2.$$

Справедливі формули:
$$\mu_3 = \nu_3 - 3\nu_2\nu_1 + 2\nu_1^3;$$
$$\mu_4 = \nu_4 - 4\nu_3\nu_1 + 6\nu_2\nu_1^2 - 3\nu_1^4.$$

Для симетричного розподілу, кожен центральний момент непарного порядку дорівнює нулю, тому кожен відмінний від нуля центральний момент непарного порядку характеризує ступінь асиметрії розподілу.

**Означення**. Величину, яка є безрозмірною, і має вигляд
$$A_s = \frac{\mu_3}{\sigma^3}$$
називають **асиметрією (skew) розподілу** (характеристика зсуву – вершини графіка щільності розподілу).

**Означення. Ексцесом (kurtosis)** називають нормований центральний момент четвертого порядку.
$$E_s = \frac{\mu_4}{\sigma^4} - 3,$$

що є характеристикою гостроверхості (плосковерхості) вершини графіка щільності розподілу та відносної кількості великих відхилень. Напр., для нормального розподілу, ця величина дорівнює 0. Для неперервного рівномірного, –1,2.

Для наочності графіки щільності розподілу ймовірностей неперервної випадкової величини f(x) за різних значень асиметрії $A_s$ та ексцесу $E_s$ наведено на рис. а) – асиметрії; б) – ексцесу.

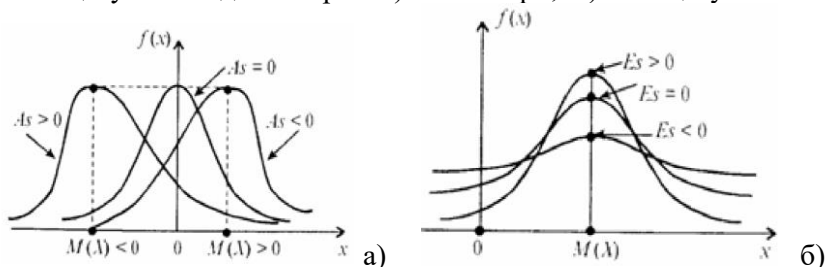

а)  б)

**Означення.** Модою $M_o(X)$ розподілу неперервної випадкової величини X із щільністю розподілу f(x) називають кожне значення x, за якого f(x) має максимум. Розподіли, які мають одну моду, називають унімодальними.

**Означення.** Медіаною $M_e(X)$ розподілу неперервної випадкової величини X називають можливе значення x, за якого пряма x= $M_e(X)$ ділить криволінійну трапецію, обмежену кривою розподілу і віссю Ox, на частини рівної площі.

• **Приклад.** Випадкова величина X задана щільністю розподілу

$$f(x) = \begin{cases} -\dfrac{3}{4}x^2 + 6x - \dfrac{45}{4}, & x \in (3, 5); \\ 0, & x \notin (3, 5). \end{cases}$$

Знайти моду, математичне сподівання і медіану випадкової величини X.

Відповідь: Скористаємось формулою для виділення повного квадрату:

$$ax^2 + bx + c = a\left(x + \frac{b}{2a}\right)^2 + c - \frac{b^2}{4a}$$

Згідно з формулою, маємо записати:

$$-\frac{3}{4}x^2 + 6x - \frac{45}{4} = -\frac{3}{4}\left(x + \frac{6}{2\cdot\left(-\frac{3}{4}\right)}\right)^2 - \frac{45}{4} - \frac{6^2}{4\left(-\frac{3}{4}\right)} =$$

$$= -\frac{3}{4}(x-4)^2 - \frac{45}{4} + \frac{36}{3} = -\frac{3}{4}(x-4)^2 + \frac{3}{4}$$

Тобто, маємо записати, що щільність розподілу ймовірностей має вигляд: $f(x) = -\frac{3}{4}(x-4)^2 + \frac{3}{4}$

Проаналізуємо отриманий вираз. Коефіцієнт при $x$ негативний, тобто, гілки параболи спрямовані до низу, це означає, що при $(x-4)^2 = 0 \Rightarrow x = 4$ будемо мати максимальне значення. Отже, згідно з визначенням: «Модою $M_o(X)$ розподілу неперервної випадкової величини X із щільністю розподілу $f(x)$ називають кожне значення $x$, за якого $f(x)$ має максимум», мода $M_o(X) = 4$.

Крива симетрична відносно прямої $x = 4$, тому $M_e(X) = 4$.

Математичне сподівання визначаємо за означенням за формулою:

$$M(X) = \int\limits_a^b x f(x) dx$$

Підставимо вираз $f(x)$, наданий в умовах задачі в цю формулу з урахуванням меж $x \in (3, 5)$:

$$M(X) = \int\limits_3^5 x \cdot \left(-\frac{3}{4}x^2 + 6x - \frac{45}{4}\right) dx$$

Після обчислень отримаємо: $M(X) = 4$.

**4.5. Основні закони розподілу дискретних випадкових величин, їх основні числові характеристики.**

**Дискретний рівномірний розподіл**

Дискретна випадкова величина рівномірно розподілена, якщо вона набуває $n$ значень з однаковими ймовірностями. Згідно з умовою нормування:

$$\sum_{i=1}^n p_i = 1, \quad p_i = \frac{1}{n}, \quad i = 1, 2, \ldots, n.$$

Прикладом рівномірно розподіленої випадкової величини є кількість очок, що випадає за одного підкидання грального кубика і закон розподілу має вигляд:

| $x$ | 1 | 2 | 3 | 4 | 5 | 6 |
|---|---|---|---|---|---|---|
| $P$ | $\frac{1}{6}$ | $\frac{1}{6}$ | $\frac{1}{6}$ | $\frac{1}{6}$ | $\frac{1}{6}$ | $\frac{1}{6}$ |

Математичне сподівання і дисперсія для рівномірно розподіленої випадкової величини 1,2,…, *n*:

$$M(X) = \frac{n+1}{2}, \quad D(X) = \frac{n^2-1}{12}.$$

**Біноміальний розподіл**

Проводять ***n*** однакових незалежних випробувань, у кожному з яких може відбутись подія А з імовірністю p (0< p<1). Випадкова величина X – кількість разів настання події А. Ряд розподілу

| $X$ | 0 | 1 | 2 | … | $k$ | … | $n$ |
|---|---|---|---|---|---|---|---|
| $P$ | $p_n(0)$ | $p_n(1)$ | $p_n(2)$ | … | $p_n(k)$ | … | $p_n(n)$ |

де

$$p_n(k) = C_n^k p^k q^{n-k}, \quad q = 1-p, \quad (k = 0, 1, 2, …, n).$$

Для біноміального розподілу

$$M(X) = np;$$
$$D(X) = npq.$$

При обчисленні таблиці ймовірностей для різних $k$, зручно користуватись рекурентними формулами: $p(0) = q^n$, $p(k) = p(k-1)\frac{n-k+1}{k}\frac{p}{q}, \ k = 1..n.$

**Розподіл Пуассона**

Якщо в схемі незалежних повторних випробувань *n* досить велике, а *p* або *1- p* прямує до нуля, тоді біноміальний закон розподілу апроксимує розподіл Пуассона, параметр якого $\lambda = np$, причому $p \leq 0{,}1$ або $p \geq 0{,}9$:

| $X$ | 0 | 1 | 2 | … | $m$ | … |
|---|---|---|---|---|---|---|
| $P$ | $p_0$ | $p_1$ | $p_2$ | … | $p_m$ | … |

Відповідні ймовірності обчислюють за формулою Пуассона

$$P(X=m) = \frac{\lambda^m}{m!}e^{-\lambda}, \lambda > 0 \ (m = 0, 1, 2, \ldots),$$ причому
$$\sum_{m=0}^{\infty}\frac{\lambda^m}{m!}e^{-\lambda} = e^{-\lambda}\sum_{m=0}^{\infty}\frac{\lambda^m}{m!} = e^{-\lambda} \cdot e^{\lambda} = 1.$$ та
$$M(X) = \lambda, \ D(X) = \lambda, \ \sigma(X) = \sqrt{\lambda}.$$

Розподіл Пуассона використовують у задачах статистичного контролю якості, в теорії надійності, теорії масового обслуговування, для прогнозування кількості вимог на виплату страхових компенсацій за рік, кількості дефектів в однакових виробах.

Рекурентні формули: $p(0) = e^{-\lambda}$, $p(k) = p(k-1)\frac{\lambda}{k}$, $k = 1..\infty$.

**Геометричний розподіл**

Проводять незалежні випробування, в кожному з яких може відбутись подія $A$ з імовірністю $p$ ($0 < p < 1$). Дискретна випадкова величина $X$ – кількість проведених випробувань до першого настання події $A$.

Ряд ймовірностей цього розподілу – нескінченно спадна геометрична прогресія зі знаменником $q = 1-p$, сума всіх членів якої дорівнює одиниці.

| $X$ | 1 | 2 | 3 | … | $n$ | … |
|---|---|---|---|---|---|---|
| $P$ | $p$ | $pq$ | $pq^2$ | … | $pq^{n-1}$ | … |

Для геометричного розподілу математичне сподівання та дисперсія:
$$M(X) = \frac{1}{p},$$
$$D(X) = \frac{q}{p^2}.$$

Також можливе інше формулювання – не кількість спроб до першої події (успіху), а кількість «не-успіхів» до першого успіху. Тоді ряд розподілу має вигляд:

| $X$ | 0 | 1 | 2 | … | $n$ | … |
|---|---|---|---|---|---|---|
| $P$ | $p$ | $pq$ | $pq^2$ | … | $pq^n$ | … |

Ця модифікація відповідає «гіпергеометричному» розподілу кількості не-успіхів до $n$ успіхів (не обов'язково, одного).

Геометричний розподіл застосовують у різноманітних задачах статистичного контролю якості виробів, у теорії надійності та страхових розрахунках. Напр., кількість випробувань до першої події з ймовірністю $p$ (успіху, відмови).

Інші застосування, напр. у середній відстані між однаковими символами в тексті, дискретизації інтервалу між подіями, які звичайно описують розподілом Пуасона (напр., інтервал між падіннями краплин дощу, чи реєстрації розпаду радіоактивних ядер лічильником Гейгера).

Рекурентні формули: $p(0) = p, \quad p(k) = p(k-1)\,q, \quad k = 1..\infty$.

**Гіпергеометричний розподіл**

Такий розподіл визначає ймовірність отримання $m$ елементів з певною властивістю серед $n$ елементів, взятих із повернення із сукупності $N$ елементів, яка містить $k$ елементів саме з такою властивістю і має вигляд:

$$P(X = m) = \frac{C_k^m \cdot C_{N-k}^{n-m}}{C_N^n}, \; m = 0,\,1,\,2,\,\ldots,\,n,\;k \geq n.$$

Математичне сподівання та дисперсія для гіпергеометричного розподілу мають вигляд:

$$M(X) = \frac{kn}{N}, \; D(X) = \frac{nk(N-k)}{N^2} \cdot \frac{(N-n)}{(N-1)}.$$

Гіпергеометричний розподіл використовують у багатьох задачах статистичного контролю якості. Напр. кількість елементів із певною властивістю серед $n$ перевірених.

Дана формула може бути узагальнена на випадок кількох типів елементів із кількістю елементів $n_i$, (загальна кількість $n$), з яких вибирають випадково $m_i$ елементів типу $i = 1..s$ (загальна кількість $m$):

$$p(m_1,..,m_s) = \frac{C_{n_1}^{m_1} C_{n_2}^{m_2} \ldots C_{n_s}^{m_s}}{C_n^m}$$

Ця формула точна для вибірки без повернення, де ймовірність наслідків змінюється із кожною спробою. Якщо кількість вибраних об'єктів $m_i n_i$

**Поліноміальний розподіл**

Поліноміальний розподіл є узагальненням біноміального і має вигляд:

$$P_n(X_1 = m_1;\ X_2 = m_2;\ \ldots,\ X_s = m_s) = \frac{n!}{m_1! \cdot m_2! \cdot \ldots \cdot m_s!} p_1^{m_1} \cdot p_2^{m_2} \cdot \ldots \cdot p_s^{m_s}.$$

Його застосовують тоді, коли внаслідок кожного із здійснених повторних випробувань може з'явитися $m_i$ із $s$ різних подій $A_i$ з ймовірностями $p_i$, причому

$$\sum_{i=1}^{s} p_i = 1.$$

Числові характеристики:

$$M(X_i) = np_i,\ D(X_i) = np_i q_i,\ q_i = 1 - p_i, i = 1,\ \ldots,\ s.$$

Цей розподіл може бути використаний, наприклад, при моделюванні ситуацій, де маємо багато класів або категорій, і хочемо визначити, скільки об'єктів випадковим чином вибрано з кожної категорії.

# 5. ОСНОВИ МАТЕМАТИЧНОЇ СТАТИСТИКИ.

В основі предмета математичної статистики багато спільного з теорією ймовірності. Обидва ці розділу математики вивчають масові випадкові явища, при цьому теорія ймовірностей виводить з математичної моделі властивості реального процесу, а математична статистика встановлює властивості математичної моделі, виходячи з даних спостережень, тобто, "Зі статистичних даних".

**Визначення**. Розділ математики, в якому вивчаються методи збору, систематизації і обробки результатів спостережень масових випадкових явищ для виявлення існуючих закономірностей, називається математичною статистикою.

## 5.1. Статистичні ряди розподілу

В результаті обробки і систематизації первинних даних статистичного спостереження отримують угрупування, звані рядами розподілу.

Ряд розподілу – це ряд чисел, що характеризує розподіл одиниць сукупності на групи за певною ознакою. Ряди розподілу характеризують склад (структуру) досліджуваного явища, дозволяють судити про однорідність сукупності, закономірності розподілу і межах варіювання (зміни) одиниць сукупності.

Залежно від того, яка ознака покладена в основу групування, розрізняють наступні види рядів розподілу рис.5.1.

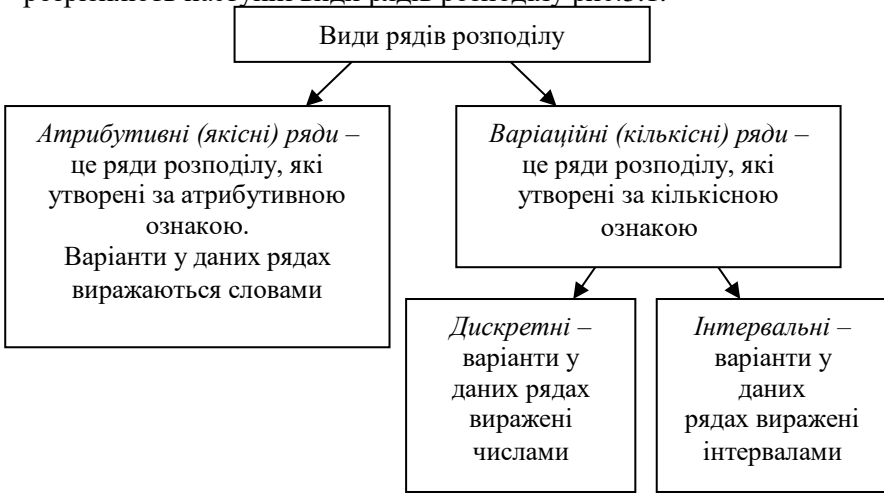

Рис.5.1. Класифікація рядів розподілу

Тобто, атрибутивний - це ряд розподілу, побудований за якісними ознаками. Він характеризує склад сукупності за різними істотними ознаками такими, як стать населення, зайнятість, національність, професія і т.д.

Варіаційний ряд розподілу будується за кількісною ознакою - наприклад, стаж роботи, вік, зарплата і ін. Він складається з двох елементів: варіанти $x_j$ і частоти $f_j$. Числові значення кількісної ознаки в варіаційному ряду розподілу називаються варіантами. Вони можуть бути додатними і негативними, абсолютними і відносними.

Частоти - це чисельність окремих варіантів або кожної групи варіаційного ряду, тобто, це числа, що показують, як часто зустрічаються ті чи інші варіанти в ряду розподілу. Сума всіх частот називається об'ємом сукупності і визначає число елементів всієї сукупності.

Варіаційні ряди залежно від характеру варіації підрозділяються на дискретні та інтервальні.

Дискретні варіаційні ряди засновані на дискретних (перерваних) ознаках, що мають тільки цілі значення, і на дискретних ознаках, представлених у вигляді інтервалів.

Інтервальні варіаційні ряди засновані на безперервних ознаках (мають будь-які значення, навіть дробові). Якщо інтервальний варіаційний ряд розподілу побудований з рівними інтервалами, частоти дозволяють судити про ступінь заповнення інтервалу одиницями сукупності. Для проведення порівняльного аналізу наповненості інтервалів визначають показник, який буде характеризувати щільність розподілу.

Щільність розподілу - це відношення числа одиниць сукупності до ширини інтервалу. Якщо інтервальний варіаційний ряд з нерівними інтервалами, то його частотні характеристики непорівнянні. Тоді, аналізуючи розподіл, використовують **щільність** частоти на одиницю інтервалу ($h_j$), тобто, $g_j = \frac{f_j}{h_j}$. Проте, базою аналізу закономірностей розподілу є варіаційний ряд — дискретний або інтервальний — з рівними інтервалами.

### 5.2. Вибірка. Види вибірки. Статистичний ряд.

В математичній статистиці є поняття **статистичної сукупності**, як маси деяких однорідних елементів, що відрізняються між собою за певними ознаками. Одиниці сукупності, з яких складається статистична сукупність, називаються **елементами** цієї сукупності або **ознаками**.

Встановлення статистичних закономірностей грунтується на вивченні **статистичних даних** – відомостей про те, які значення прийняла окрема ознака (випадкова величина X) унаслідок проведення досліду.

На практиці статистичних досліджень відрізняють два види дослідів:

– **суцільний**, коли розглядаються всі елементи сукупності;

– **вибірковий**, де вивчається лише деяка частина елементів.

Вся сукупність елементів, яку треба вивчити називається **генеральною сукупністю**. Та частина об'єктів, яка відібрана для безпосереднього вивчення із генеральної сукупності, називається **вибірковою сукупністю** або просто –**вибіркою**. Кількість елементів у генеральній чи вибірковій сукупності називають їх **об'ємами**.

**Сутність вибіркового методу** полягає в тому, щоб за деякою частиною генеральної сукупності, тобто, за вибіркою, робити висновки про її властивості в цілому, наприклад, про її закон розподілу, або про числові значення її певних параметрів. Але цей метод має недолік під назвою **помилки репрезентативності**. Щоб за даними вибірки мати можливість судити про генеральну сукупність, вона повинна бути взята випадково, це у певній мірі дозволяє знизити можливість помилок репрезентативності. Вибіркову сукупність називають **репрезентативною**, якщо вона досить добре відбиває основні характеристики генеральної сукупності.

Розрізняють наступні види вибірок:

- **випадкова вибірка**, отримана випадковим відбором елементів без поділу їх на частини або групи;

- **механічна вибірка**, для якої елементи із генеральної сукупності відбираються через деякий інтервал;

- **типова вибірка**, у яку випадковим чином вибираються елементи з типових груп, на які за деякою ознакою поділяється генеральна сукупність;

- **серійна вибірка**, у яку випадковим чином потрапляють не елементи груп, а власне групи, які потім суцільно досліджуються.

За способом можливого утворення вибірки поділяються на:

- **вибірки з повторенням**, коли елемент після вивчення повертається назад до сукупності, що вивчається та може бути

повторно вивчений;

**- вибірки без повторення**, коли елемент після вивчення назад до сукупності не повертається.

При обробці статистичних даних виділяються наступні три етапи:

- етап спостережень, на якому ставиться мета статистичного дослідження та проводиться збирання статистичного матеріалу;

- етап зведення та групування отриманого статистичного матеріалу, представлення його у вигляді таблиць;

- етап, на якому застосовують великий арсенал методів математичної статистики, щоб отримати відповідь на завдання мети дослідження. Розглянемо обробку статистичних даних з другого етапу тобто, подання даних у вигляді так званих статистичних рядів.

Математична статистика досліджує випадкову величину, яка змінюється якимось заздалегідь невідомим чином. Позначимо цю величину як $\xi$. Цій випадковій величині $\xi$ відповідає невідома функція розподілу $F(x)$. Припустимо, що випробування проводяться в однакових умовах і незалежно.

В результаті цих експериментів з визначеною кількістю випробувань $n$ отримаємо набір значень $\xi_1, \xi_2, \ldots \xi_n$ випадкової величини $\xi$. Тоді набір значень $\xi_1; \xi_2; \xi_3; \ldots \xi_n$ випадкової величини $\xi$, отриманих в результаті $n$ – випробувань, називається вибіркою об'єму $n$. Таким чином, ми маємо сукупність $n$ незалежних випадкових величин $\xi (i = 1,2,\ldots,n)$ кожна з яких розглядається як випадкова величина $\xi$. Про таку вибірку кажуть, що вона взята з генеральної сукупності випадкової величини $\xi$ з теоретичною функцією розподілу $F(x)$.

Якщо вибірка об'єму $n$ має $r$ різних елементів (значень) $x_1; x_2; x_3;\ldots x_r$ і елемент $x_i$ зустрічається $m_i$ разів, то число $m_i$ називається частотою елементу $x_i$, а відношення $W_i = \dfrac{m_i}{n}$ - називається відносною частотою елементу $x_i$.

Очевидно, що
$$\sum_{i=1}^{r} m_i = n \quad \text{та} \quad \sum_{i=1}^{r} W_i = 1$$

де $m_i$ – це кількість однакових значень; $r$ – кількість різних елементів, $W_i$ – відносна частота, тобто, кількість однакових значень по відношенню до загальної кількості значень вибірки в цілому.

Таким чином, ми маємо пари значень ($x_i$; $m_i$ ), $i=1,2,…, r$, тобто, різним значенням вибірки $x_1$; $x_2$; $x_3$;... $x_r$ відповідає своя кількість повторів $m_i$.

**Статистичним рядом** називається послідовність пар *($x_i$; $m_i$)*, де *i=1, 2, 3,..., n.*

Статистичний ряд називають **дискретним**, якщо будь-які його варіанти різняться постійною величиною, та – **неперервним** (*інтервальним*), якщо варіанти можуть різнитися між собою, будь-якою маленькою величиною. Як правило, статистичний ряд записують у вигляді таблиці, розташовуючи значення $x_i$ в порядку зростання.

| $x_i$ | $x_1$ | $x_2$ | … | $x_{r-1}$ | $x_r$ |
|---|---|---|---|---|---|
| $m_i$ | $m_1$ | $m_2$ | … | $m_{r-1}$ | $m_r$ |

Статистична таблиця - це особливий спосіб короткої і наочної записи, який дозволяє охопити матеріали статистичного зведення в цілому, в ній зручно розташовуються обчислені показники і відомості про досліджуваний об'єкт, що викладаються цифрами на основі певного порядку.

Ще більш наочним методом досліджень являється графічне зображення. При графічному зображенні варіаційного ряду користуються прямокутною системою координат. По горизонтальній вісі **X** відкладають значення варіантів, а по вертикальній вісі **Y** — абсолютні або відносні значення частот.

Для графічного зображення статистичних рядів найчастіше використовують полігон частот, гістограму та кумулятивну криву (кумуляту).

**Полігон частот** застосовується переважно для зображення дискретних рядів розподілу, але може бути використаний і для зображення рядів розподілу з безперервним варіюванням ознаки, яка вивчається. Цей графік має вигляд не стовпчиків, а ламаної. Будується він також на прямокутній системі координат. По вісі абсцис відкладають варіанти, а по вісі ординат — частоти. Усі вершини з'єднують між собою.

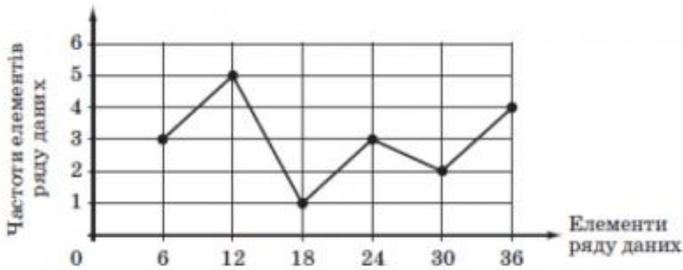

Рис. 5.2. Графічне подання інформації про вибірку. Полігон частот

**Гістограма** застосовується переважно для зображення інтервальних варіаційних рядів розподілу. Гістограма це засіб подання статистичних даних в виді стовпчиків, які відображають розподіл окремих вимірювань параметрів якогось виробу чи процесу. Іноді її називають частотним розподілом тому що показує частоту появи вимірюваних значень параметрів виробу. Вона будується так: на вісі абсцис відкладають інтервали ознаки, а на вісі ординат — їх чисельність (частоти). На відрізках, які виражають інтервали, будують прямокутники, висота яких дорівнює частоті даного інтервалу. Площа кожного стовпчика повинна бути пропорційна до частот. Тобто, висота кожного стовпчика вказує на частоту появи значень параметру в обраному діапазоні, а кількість стовпців – на число обраних діапазонів.

Такий графік дає уявлення про те, як заповнений кожний інтервал.

**Кумулята.** Кожний варіаційний ряд можна графічно зобразити у вигляді кумуляти. Для цього треба на вісі абсцис відкласти варіанти або межі інтервалів, а на вісі ординат - відповідні наростаючі підсумки частот. Одержані точки треба з'єднати плавною кривою лінією, яка має назву кумуляти або кумулятивної кривої.

Отже, полігон, гістограма і кумулята дають змогу мати більш наочне уявлення про розподіл того чи іншого явища.

**5.3. Емпірична (вибіркова) функція розподілу.**

**Емпіричною функцією розподілу** називається функція $F_n(x) = \dfrac{k}{n}$, де $k$ – число вибіркових значень, менших за $x$, тобто, $F_n(x)$ являється відносною частотою події $\xi < x$. Функція $F_n(x)$ служить оцінкою невідомої (теоретичної) функції розподілу *F(x)*, тобто, із нескінченно великої кількості значень буде вибрано тільки *n* значень і саме по цієї, обмеженої кількості значень ми спробуємо зробити припущення та висновки про генеральну сукупність.

Статистичний ряд – це статистичний аналог розподілу ознаки *x* (випадкової величини ξ). В цьому сенсі полігон (гістограма) – аналог кривої розподілу, а емпірична функція розподілу $F_n(x)$ - функції розподілу *F(x)* випадкової величини ξ. *Запис вибірки у вигляді статистичного ряду і його зображення полігоном частот зручно, якщо досліджувана випадкова величина ξ є дискретною, а число r різних вибіркових значень ($x_i$) не дуже велике.*

Розгляд та розуміння даних, що були отримані за результатами вибіркового дослідження, особливо при великій кількості спостережень *n*, провести досить важко, та за ними практично неможливо уявити характер розподілу ознаки. Тому, щоб покращити розуміння статистичного матеріалу, який маємо унаслідок проведення деякого дослідження, проводимо упорядковування, тобто, розташування варіантів за зростанням (чи зменшенням) чи інакше - **ранжування** ряду. В такому вигляді вивчати дані також не завжди зручно. Це пов'язано з великою кількістю елементів, тому їх розбивають на окремі інтервали, тобто, проводиться **групування**. Кількість інтервалів *r* повинна бути не дуже великою, але й не зовсім малою, щоб не втратити особливостей розподілу ознаки.

Кількість груп тісно пов'язана з метою дослідження та обсягом сукупності. У масових сукупностях орієнтовано оптимальну кількість груп з рівними інтервалами можна визначити за формулою американського вченого Стерджеса:
$$r = 1+ \log_2 n = 1 + 3{,}32 \cdot \lg n,$$
де *r* – кількість груп (інтервалів); *n* – об'єм сукупності.

Однак користуватися цією формулою можна лише в тих випадках, коли досліджувана сукупність досить велика і зміна групувальної ознаки має відносно плавний характер.

При групуванні за кількісною ознакою часто виникає питання про інтервали. Існують різні види інтервалів: рівні та нерівні; відкриті та закриті.

Рівні інтервали застосовують тоді, коли кількісна ознака всередині сукупності змінюється плавно, поступово, рівномірно.

Розмір рівного інтервалу визначають за формулою:
$$h = \frac{x_{max} - x_{min}}{r},$$
де $x_{max}$; $x_{min}$ – максимальне, мінімальне значення ознаки; $r$ – кількість груп.

В практичних роботах, як правило, отримане значення $h$ округляють збільшуючи цей параметр, щоб не зменшити загальний діапазон значень випадкової величини, яка досліджується. Таке збільшення рекомендують зробити так, щоб $x_{min}$ та $x_{max}$ потрапляли приблизно в середину першого і останнього інтервалів. Так, при такому рішенні, усі крайні точки гарно потраплять у центр інтервалу, але, вирішуючи одну проблему з межами крайніх інтервалів створюємо іншу. При такому підході розміри першого та останнього інтервалів вдвічі зменшуються за рахунок віддалення меж де точно немає елементів. Це призведе до заниження показників крайніх інтервалів при побудові гістограми. Ця похибка може бути суттєвою при великої кількості елементів, що потрапляють до інтервалу, тому в цьому випадку слід скористатися формулою:
$$h = \frac{x_{max} - x_{min}}{r} \cdot \frac{n}{n-1}$$

При такому підрахунку, з введенням коефіцієнту $\frac{n}{n-1}$, відступ від меж дуже малий, що суттєво зменшує похибку, але достатній, щоб усі елементи потрапили у відповідний інтервал.

Нерівними називають інтервали, в яких різниця між верхньою і нижньою межею неоднакова. Необхідність застосування групування з нерівними інтервалами виникає в тих випадках, коли коливання ознаки має нерівномірний характер у великих межах.

Інтервали при групуванні можуть бути закритими та відкритими. Закритими вважаються інтервали, у яких визначені верхня та нижня межі, відкритими називаються інтервали, у яких нижня та верхня межі невідомі. Інколи використовують відкриті інтервали з однією межею (верхньою або нижньою). Потреба у відкритих інтервалах зумовлена високим коливанням досліджуваної ознаки.

При розподілі одиниць об'єкта спостереження на окремі групи важливо точно визначати межі. Для цього при побудові інтервалів одне й те саме число повторюється: як верхня (права) межа одного інтервалу, так і нижня (ліва) межа наступного інтервалу.

Таким чином, у випадку великої кількості одиниць об'єкта спостережень робимо розподіл на окремі групи з точно визначеними межами.

**5.4. Числові характеристики вибірки.**

Варіаційні ряди та графіки емпіричних розподілів дають наочне уявлення про те, як варіює ознака в вибіркової сукупності, але цього замало для повної характеристики вибірки. Щоб охопити усі деталі необхідно застосування деякі узагальнюючи числові характеристики. Саме числові характеристики вибірки дають кількісне уявлення про емпіричні данні і дозволяють порівнювати їх між собою. Найбільше практичне значення мають характеристики положення, розсіювання і асиметрії емпіричних розподілів. Крім того, в практиці вибіркового методу використовують два типи вибіркових оцінок – точкові та інтервальні. Згадаємо, що:

- упорядкований ряд даних - це ряд даних, в якому кожне наступне число не менше попереднього;
- статистичним рядом називається послідовність пар (*$x_i$*; *$m_i$*), де *i=1, 2, 3, … n.*
- інтервальний статистичний ряд – це статистичний ряд, який поділений на *r* інтервалів з шириною *h.*

Тепер розглянемо характеристики, що визначають положення центру емпіричного розподілу, до яких відносяться: середнє арифметичне, медіана та мода. **Середнє арифметичне** - це таке значення признака, при якому сума відхилень вибіркових значень ознаки дорівнюється нулю (з урахуванням знаку відхилення). Середнє прийнято позначати тією ж літерою, що й варіанти вибірки, тільки над літерою ставиться символ усереднення – риска ($\bar{x}$). Середнє арифметичне може обчислюватися як по необробленім первинним даним, так і по результатам групування цих даних.

Для усіх випадкових величин ξ, що не мають розділення на групи, тобто, точковою оцінкою для математичного очікування являється середнє арифметичне яке визначається по формулі:

$$\bar{x} = \frac{\xi_1 + \xi_2 + \cdots + \xi_n}{n}; \qquad \bar{x} = \frac{1}{n}\sum_{i=1}^{n} \xi_i ,$$

де *n* – об'єм вибірки; $\xi_i$ – випадкова величина вибірки. Тобто, середнє значення вибірки це відношення суми усіх значень випадкової величини $\xi$ до її об'єму $n$, чи інакше, щоб його знайти треба скласти усі значення та поділити на їх кількість.

Якщо вибірка записана у вигляді статистичного ряду, то формулу для середнього арифметичного можна записати:

$$\bar{x} = \frac{m_1 x_1 + m_2 x_2 + \cdots + m_r x_r}{n} = \frac{1}{n}\sum_{i=1}^{r} m_i x_i,$$

де $m_i$ - кількість елементів $i$-го варіанту вибірки; $x_i$ – варіант значення вибірки; $r$ – кількість різних варіантів вибірки.

Коли дані згруповані, то

$$\bar{x} = \frac{m_1 \tilde{x}_1 + m_2 \tilde{x}_2 + \cdots + m_r \tilde{x}_r}{n} = \frac{1}{n}\sum_{i=1}^{r} m_i \tilde{x}_i,$$

де $r$ – вже означає кількість інтервалів, а $\tilde{x}_i$ – значення середини інтервалу, тобто, всі вибіркові значення, що потрапили в $i$-тий інтервал, замінюються одним – серединним значенням цього інтервалу (дивитись другий рядок до таблиці практичного завдання №2). Середнє арифметичне обчислене по цій формулі називають також середнім зваженим, підкреслюючи цим, що числа $\tilde{x}_i$ підсумовуються з коефіцієнтами (важелями), рівними частотам (кількості значень вибірки), що потрапили в інтервали групування.

**Мода вибірки** (Мо) – являє собою значення ознаки, що зустрічається у вибірці найбільш часто, тобто, це варіанта вибірки, що має найбільшу частоту.

Ряд називається унімодальним, якщо в ньому тільки одне модальне значення і полімодальним в іншому випадку. Для полімодального ряду значення моди не визначають, та визначають по результатам групування даних, тобто, по інтервалам. Інтервал групування з найбільшою частотою називається модальним, для визначення моди в інтервальному ряду використовується формула:

$$Mo = x_{\text{Мон}} + h\frac{m_{\text{Мо}} - m_{\text{Мо}-1}}{(m_{\text{Мо}} - m_{\text{Мо}-1}) + (m_{\text{Мо}} - m_{\text{Мо}+1})},$$

де $x_{Мон}$ - нижня межа модального інтервалу; $h$ – ширина інтервалу групування; $m_{Мо}$ – частота модального інтервалу; $m_{Мо-1}$ – частота інтервалу, попереднього модальному; $m_{Мо+1}$ - частота інтервалу, наступного після модального.

**Медіана -** якщо в упорядкованому ряду непарна кількість елементів, то медіаною є середній за номером елемент, якщо парна кількість елементів, то медіаною є середнє арифметичне двох середніх за номером елементів (рис.5.3).

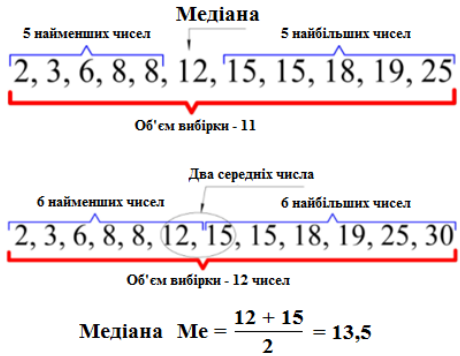

Рис. 5.3. Приклади знаходження медіани упорядкованого ряду

Таким чином, щоб знайти медіану дискретного ряду його дані спочатку розташовують в порядку зростання чи зменшення і вже в цьому упорядкованому ряду, з об'ємом *n* членів шукається ранг медіани *R*, тобто, порядковий номер елемента, який є середнім за номером. Ранг медіани упорядкованої вибірки з непарною кількістю членів: $R_{Me} = \frac{n+1}{2}$.

Щоб знайти ранг медіани упорядкованої вибірки з парною кількістю елементів користуються тією ж формулою $R_{Me} = \frac{n+1}{2}$, але отриманий результат дрібне число, що вказує на два номери упорядкованого ряду, між якими воно розташовано, припустимо між *i*-тим та (*i* +1), тому треба взяти значення вибірки $\xi_i$ та $\xi_{i+1}$, які мають ці номери в упорядкованому ряду скласти їх і поділити навпіл: $Me = \frac{\xi_i + \xi_{i+1}}{2}$.

При знаходженні медіани для згрупованих даних спочатку треба знайти інтервал групування, в якому знаходиться медіана шляхом підрахування накопичених частот чи накопичених відносних частот. Медіанним є той інтервал, в якому накопичена частота вперше виявиться більшою ніж $\frac{n}{2}$ (*n* - об'єм вибірки) чи накопичена відносна частота виявиться більшою ніж 0,5. Зобразимо згрупований варіаційний ряд графічно у вигляді кумуляти. Для цього треба на вісі абсцис відкласти межі інтервалів, а на вісі ординат - відповідні наростаючі підсумки частот рис. 5.4. Всередині медіанного інтервалу медіана визначається по формулі:

$$Me = x_{Me_{\text{н}}} + h_{Me}\frac{0{,}5n - n_{x_{Me-1}}}{n_{Me}},$$

де $x_{Me_{\text{н}}}$ - нижня межа медіанного інтервалу; $0{,}5n$ чи $\frac{n}{2}$ – пів об'єму вибірки; $h_{Me}$ – ширина медіанного інтервалу; $n_{x_{Me-1}}$ – накопичена частота інтервалу попереднього медіанному; $n_{Me}$ – частота медіанного інтервалу (кількість елементів упорядкованої вибірки, які потрапили у цей інтервал).

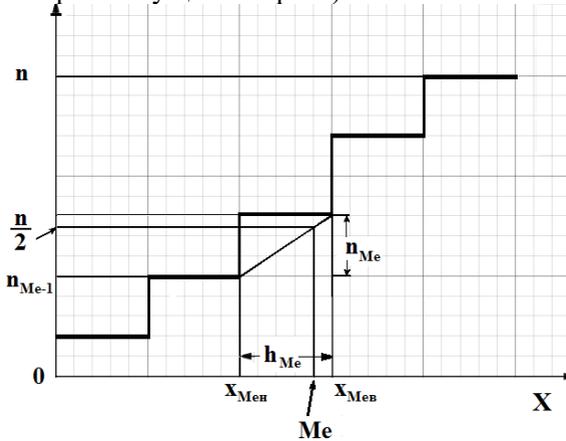

Рис. 5.4. Визначення медіани для згрупованих даних

Для оцінювання та вимірювання варіацій, тобто, зміни ступеню відхилення ознак використовуються наступні характеристики вибірки:

**Варіаційний розмах** - різниця між максимальним і мінімальним значеннями ознаки деякої сукупності $R = x_{max} - x_{min}$ ;

**Середнє лінійне і середнє квадратичне відхилення та дисперсія** це найважливіші характеристики розсіювання. При аналізі ми використовуємо середнє значення, тому відхилення в обидві сторони в (+) та в (-) в сумі дадуть (0), тобто:

$$\sum_{i=1}^{n}(\xi_i - \overline{x}) = 0 \quad \text{чи} \quad \sum_{i=1}^{r}(\widetilde{x}_i - \bar{x})m_i = 0$$

Тому використовуються середні або модуля відхилень, або квадрату відхилень:

$$\frac{1}{n}\sum_{i=1}^{n}|\xi_i - \bar{x}|, \quad \text{або} \quad \frac{1}{n}\sum_{i=1}^{n}(\xi_i - \bar{x})^2,$$

де $\xi_i$ – $i$-тий елемент вибірки, $n$ – кількість елементів вибірки; $\bar{x}$ - середнє арифметичне значення вибірки. Чи майже теж саме записане інакше:

$$\frac{1}{n}\sum_{i=1}^{r}|x_i - \bar{x}|m_i, \quad \text{або} \quad \frac{1}{n}\sum_{i=1}^{r}(x_i - \bar{x})^2 m_i$$

де $i$ – порядковий номер впорядкованого ряду різних значень вибірки; $r$ – кількість різних значень вибірки; $m_i$ – кількість значень $x_i$. Та для даних згрупованих в інтервальний статистичний ряд маємо аналогічні формули:

$$\frac{1}{n}\sum_{i=1}^{r}|\tilde{x}_i - \bar{x}|m_i, \quad \text{або} \quad \frac{1}{n}\sum_{i=1}^{r}(\tilde{x}_i - \bar{x})^2 m_i,$$

де $i$ – номер інтервалу, $r$ – число інтервалів, $\tilde{x}_i$ – значення середини інтервалу, $m_i$ – кількість елементів вибірки $\xi$, що потрапили в $i$-тий інтервал.

Отже, узагальнюючою характеристикою варіації буде середнє відхилення і нижче наведені формули для: а) лінійного; б) квадратичного відхилень та в) дисперсія $\sigma^2$, де перший вираз – для не згрупованих даних (тобто, вибірка не була поділена на групи або інтервали) і другий вираз для варіаційних даних згрупованих в інтервальний статистичний ряд.

Ці формули будуть мати наступний вигляд:

а) лінійне

$$\bar{l} = \frac{\sum_{i=1}^{n}|\xi_i - \bar{x}|}{n} \quad \text{чи} \quad \bar{l} = \frac{\sum_{i=1}^{r}|\tilde{x}_i - \bar{x}|m_i}{\sum_{i=1}^{r}m_i}$$

б) квадратичне

$$\sigma = \sqrt{\frac{\sum_{i=1}^{n}(\xi_i - \bar{x})^2}{n}} \quad \text{чи} \quad \sigma = \sqrt{\frac{\sum_{i=1}^{r}(\tilde{x}_i - \bar{x})^2 m_i}{\sum_{i=1}^{r} m_i}}$$

в) дисперсія

$$\sigma^2 = \frac{1}{n}\sum_{i=1}^{n}(\xi_i - \bar{x})^2$$

де $\xi_i$ - значення вибірки, $i$ = 1, 2, …$n$.

$$\sigma^2 = \frac{1}{n}\sum_{i=1}^{r}(\tilde{x}_i - \bar{x})^2 m_i$$

де $i$ – номер інтервалу; $r$ – число інтервалів; $m_i$ – кількість елементів вибірки ξ, що потрапили в $i$-тий інтервал.

Тобто, дисперсія – це середній квадрат відхилень, а вирази:

$$\sum_{i=1}^{n}(\xi_i - \bar{x})^2$$

– це сума квадратів відхилень значень ознаки ξ від середнього арифметичного $\bar{x}$;

$$\sum_{i=1}^{r} m_i(\tilde{x}_i - \bar{x})^2$$

- це зважена сума квадратів відхилень.

Щоб отримати середній квадрат відхилень ці суми в першому та другому варіантах, були поділені на об'єм вибірки $n$ чи суму $m_i$.

$$\sum_{i=1}^{r} m_i = n$$

Але, щоб отримати «незміщену» оцінку дисперсії, суми квадратів відхилень треба поділити на $(n-1)$ і отримаємо такі формули:

$$S^2 = \frac{1}{n-1} \cdot \sum_{i=1}^{n}(\xi_i - \bar{x})^2 = \frac{n}{n-1}\sigma^2,$$

цю формулу можна перетворити, розкривши дужки

$$S^2 = \frac{1}{n-1} \cdot \sum_{i=1}^{n}((\xi_i)^2 - 2\xi_i\xi + \xi^2)$$

Після нескладних перетворень і пам'ятаючи, що $\xi = \frac{1}{n}\sum \xi_i$, тобто, $\sum \xi_i = n\xi$ отримаємо вираз:

$$S^2 = \frac{1}{n-1} \cdot \sum_{i=1}^{n}(\xi_i{}^2 - n\bar{x}^2).$$

У випадку вибірки, представленої у вигляді згрупованих варіаційних даних, тобто, маємо інтервальний статистичний ряд формула незміщеної дисперсії має вигляд:

$$S^2 = \frac{1}{n-1} \cdot \left(\left(\sum_{i=1}^{r} m_i \tilde{x}_i{}^2\right) - n\bar{x}^2\right),$$

де $\tilde{x}_i$ – значення середини інтервалу, а

$$\bar{x} = \frac{1}{n}\sum_{i=1}^{r} m_i \tilde{x}_i,$$

Тобто,

$$S^2 = \frac{1}{n-1} \cdot \left(\sum_{i=1}^{r} m_i \tilde{x}_i{}^2 - \frac{1}{n}\left(\sum_{i=1}^{r} m_i \tilde{x}_i\right)^2\right)$$

### 5.5. Приклад статистичного дослідження вибірки.

Основу будь-якого статистичного дослідження становить безліч даних, отриманих в результаті експерименту, спостережень, вимірювань різних ознак. В якості таких спостережень розглянемо вибірку, яка складається з елементів:

6, 1, 9, 7, 7, 6, 3, 5, 6, 2, 7, 4, 8, 0, 9, 3, 5, 7, 5, 2, 6, 5, 4, 6, 9, 10, 3, 6, 4, 1, 6, 4, 2, 7, 5, 4, 6, 2, 4, 9, 3, 5, 2, 7, 2, 10, 8, 8, 4, 5.

1. Дано: вибірка випадкової величини ξ, об'єм вибірки ***n*** у цьому випадку дорівнюється 50 (але може бути і інше значення, порахуйте кількість числових значень вибірки у вашому варіанті).

2. Розташуємо різні варіанти вибірки (*x*$_i$) у порядку збільшення: 0, 1, 2, 3, 4, 5, 6, 7, 8, 9, 10 та обчислимо кількість повторень.

Отримані дані внесемо у таблицю за зразком:

| $x_i$ | $x_1$ | $x_2$ | … | $x_{r-1}$ | $x_r$ |
|---|---|---|---|---|---|
| $m_i$ | $m_1$ | $m_2$ | … | $m_{r-1}$ | $m_r$ |

3. Вигляд таблиці для значень нашої вибірки

| $x_i$ | 0 | 1 | 2 | 3 | 4 | 5 | 6 | 7 | 8 | 9 | 10 |
|---|---|---|---|---|---|---|---|---|---|---|---|
| $m_i$ | 1 | 2 | 6 | 4 | 7 | 7 | 8 | 6 | 3 | 4 | 2 |

4. Статистичний ряд також можна зобразити графічно – це буде зламана лінія з вершинами у точках ($x_i$; $m_i$), якщо $m_i$>1. Розглянемо наш статистичний ряд, в якому усі варіанти вже розміщені у зростанні (пункт 3). Ми бачимо, що перше значення $x_1$= 0 зустрічається у нашої вибірки 1 раз, тобто, $m_1$=1, і так до $x_{11}$= 10 , це значення повторюється у вибірки 2 рази, тобто, $m_{11}$=2 і так для усіх значень $x_i$. Тобто, ми маємо точки з координатами: (0; 1), (1; 2), (2; 6), (3; 4), (4; 7), (5; 7), (6; 8), (7; 6), (8; 3), (9; 4), (10; 2). Розмістимо їх на графіку.

5. Для побудови графіка зручно користуватись прямокутною системою координат з різними масштабами по осях, де по вісі X ми відкладаємо числове значення, а по вісі Y відкладаємо кількість однакових значень (див. рис.5.5). Якщо у вибірки існують негативні цифрові значення, тоді зручно перенести вісь Y вліво щоб вона не заважала дивитись на графічне зображення.

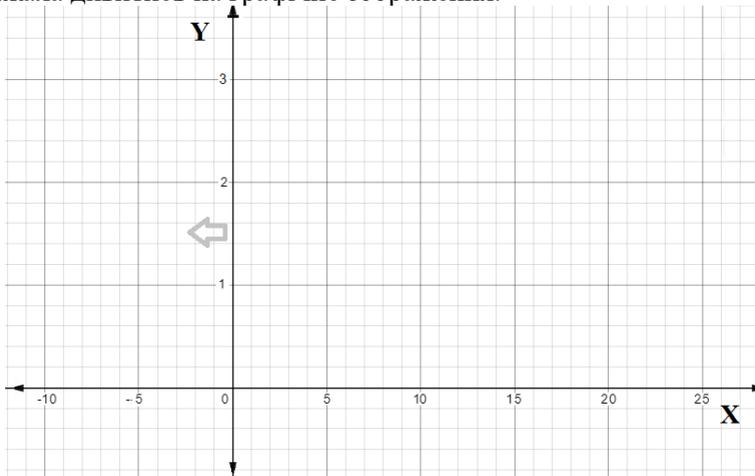

Рис. 5.5. Загальний вигляд прямокутної системи координат

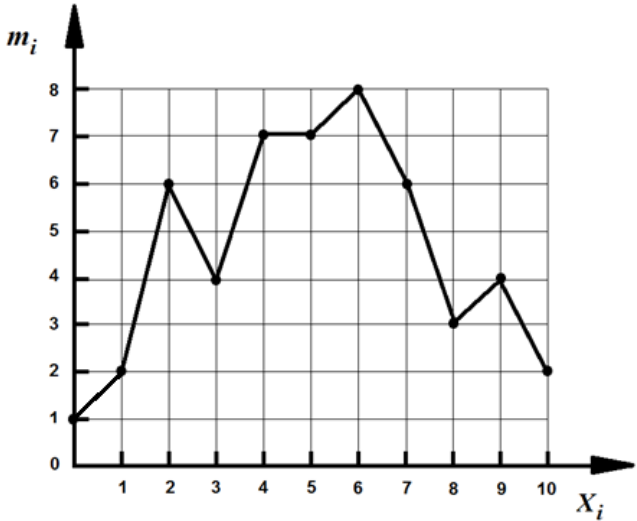

Рис. 5.6. Полігон частот

6. Графічне зображення статистичного ряду, коли по вісі X відкладається значення $x_i$, а по вісі Y показник $m_i$ називається полігоном частот. Користуватися таким представленням результатів зручно при дискретних ξ (дзета, елемент виборці) та невеликих значеннях **r** (кількість різних значень у виборці). Дивлячись у таблицю з даними нашого статистичного ряду побудуємо графік під назвою «полігон частот» рис. 5.6.

Після того, як варіаційний ряд був упорядкований визначимо варіаційний розмах; середнє арифметичне по необробленим первинним даним та по упорядкованим; моду вибірки (Мо), медіану (Ме) та дисперсію.

1. Варіаційний розмах $R = x_{max} - x_{min}$, щоб його визначити розглянемо таблицю даних (для приклада візьмемо таблицю для значень нашої вибірки):

| $x_i$ | 0 | 1 | 2 | 3 | 4 | 5 | 6 | 7 | 8 | 9 | 10 |
|---|---|---|---|---|---|---|---|---|---|---|---|
| $m_i$ | 1 | 2 | 6 | 4 | 7 | 7 | 8 | 6 | 3 | 4 | 2 |

$x_{max} = 10$   $x_{min} = 0$   $R = 10$

2. Для усіх випадкових величин ξ, що не мають розділення на групи, тобто, точковою оцінкою для математичного очікування являється середнє арифметичне яке визначається по формулі:

$$\bar{x} = \frac{1}{n} \sum_{i=1}^{n} \xi_i ,$$

Вибірка складається із елементів:

6, 1, 9, 7, 7, 6, 3, 5, 6, 2, 7, 4, 8, 0, 9, 3, 5, 7, 5, 2, 6, 5, 4, 6, 9, 10, 3, 6, 4, 1, 6, 4, 2, 7, 5, 4, 6, 2, 4, 9, 3, 5, 2, 7, 2, 10, 8, 8, 4, 5.

$$\bar{x} = \frac{\xi_1 + \xi_2 + \cdots + \xi_n}{n} = \frac{6+1+9+7+\cdots+8+8+4+5}{50} = \frac{259}{50} = 5{,}18$$

3. Якщо вибірка записана у вигляді статистичного ряду, то формулу для середнього арифметичного можна записати:

$$\frac{1}{n} \sum_{i=1}^{r} m_i x_i ,$$

Із таблиці (пункт 1) беремо значення $x_i$ (із першої строки) та відповідне йому значення $m_i$ (із другої строки)

$$\bar{x} = \frac{m_1 x_1 + m_2 x_2 + \cdots + m_r x_r}{n} = \frac{1 \cdot 0 + 2 \cdot 1 + 6 \cdot 2 + \cdots + 2 \cdot 10}{50} = 5{,}18$$

4. Ранг медіани упорядкованої вибірки обчислюємо по формулі:

$$R_{Me} = \frac{n+1}{2} = \frac{50+1}{2} = 25{,}5.$$

Отриманий результат дрібне число, що вказує на два номери 25 та 26 упорядкованого ряду, між якими розташовано значення Me. Користуючись таблицею (пункт 1) визначаємо значення елементів упорядкованої вибірки під 25 та 26 номерами:

| Номер | 20 | 21 | 22 | 23 | 24 | 25 | 26 |
| --- | --- | --- | --- | --- | --- | --- | --- |
| Значення | 4 | 5 | 5 | 5 | 5 | 5 | 5 |

$$Me = \frac{\xi_{25} + \xi_{26}}{2} = \frac{5+5}{2} = 5.$$

5. Мода вибірки (Мо) – являє собою значення, яке має найбільшу частоту. Розглянемо значення $m_i$ другої строки таблиці (пункт 1). Найбільшому $m_i$=8 відповідає $x_i$=6, тобто, Мо=6.

6. Для визначення дисперсії безпосередньо за вибіркою, яка складається із елементів:

6, 1, 9, 7, 7, 6, 3, 5, 6, 2, 7, 4, 8, 0, 9, 3, 5, 7, 5, 2, 6, 5, 4, 6, 9, 10, 3, 6, 4, 1, 6, 4, 2, 7, 5, 4, 6, 2, 4, 9, 3, 5, 2, 7, 2, 10, 8, 8, 4, 5; обчислюємо по формулі:

$$S^2 = \frac{1}{n-1} \cdot \sum_{i=1}^{n}(\xi_i - \bar{x})^2$$

спочатку рахується різниця між кожним елементом вибірки і середнім значенням, яке вже пораховано в пункті 2, тобто, знайдемо різницю між окремим значенням і середнім, що відображає міру відхилення:

$\xi_1 - \bar{x}; \xi_2 - \bar{x}; \ldots, \xi_n - \bar{x}$, де $\xi_1$ – перше значення упорядкованої вибірки; $\xi_n$ – останнє значення вашої вибірки.

Підставляємо числові значення $6 - 5{,}18; 1 - 5{,}18; \ldots ; 5 - 5{,}18;$

7. Отримані значення різниці зводимо в квадрат для того, щоб усі відхилення стали виключно додатними числами і щоб уникнути взаємознищення додатних і негативних відхилень при їх підсумовуванні:

$(\xi_1 - \bar{x})^2; (\xi_2 - \bar{x})^2; \ldots ; (\xi_n - \bar{x})^2$
$(0{,}82)^2; (-4{,}18)^2; \ldots ; (-0{,}18)^2$

8. Всі ці пораховані відхилення зведені в квадрат підсумовуються і рахується середнє для чого отримане значення поділимо на **$n-1$** для незміщеної дисперсії:

$$S^2 = \frac{(\xi_1 - \bar{x})^2 + (\xi_2 - \bar{x})^2 + \cdots + (\xi_n - \bar{x})^2}{n-1} = \frac{305{,}38}{50-1} = 6{,}23$$

Середній квадрат відхилення і є дисперсією.

**5.6. Приклад дослідження вибірки наданої у вигляді інтервального статистичного ряду.**

В якості спостережень, розглянемо вибірку, яка складається з таких елементів:

7,59; 7,48; 7,46; 7,40; 7,24; 7,41; 7,34; 7,43; 7,38; 7,50; 7,26; 7,43; 7,37; 7,55; 7,42; 7,41; 7,38; 7,14; 7,42; 7,52; 7,46; 7,39; 7,35; 7,32; 7,18; 7,30; 7,54; 7,38; 7,37; 7,34; 7,50; 7,61; 7,42; 7,32; 7,35; 7,40; 7,57; 7,31; 7,40; 7,36; 7,28; 7,58; 7,38; 7,58; 7,26; 7,37; 7,28; 7,39; 7,32; 7,20; 7,43; 7,34; 7,45; 7,33; 7,41; 7,33; 7,45; 7,31; 7,45; 7,39.

Необхідно записати вибірку у вигляді інтервального статистичного ряду та побудувати гістограму вибірки.

1. Дано: вибірка числових значень з *n*=60 у цьому випадку (але може бути і інше число, порахуйте кількість елементів вибірки у

вашому варіанті).

2. У статистичному ряду числові значення розташовуються в порядку зростання. Тому необхідно розсортувати елементи вибірки в порядку збільшення, тобто, від меншого значення до більшого. Таким чином, отримаємо крайні значення вибірки, тобто, $x_{min}$ =7,14 та $x_{max}$=7,61 .

3. Необхідно визначити, на скільки інтервалів потрібно розбити впорядковану (тобто, числа розташовані в порядку збільшення) вибірку. Для цього застосуємо формулу Стерджеса :
$$r = 1 + \log_2 n \approx 1 + 3{,}32 \cdot \lg n$$
У нашому випадку,
$$r = 1 + 3{,}32 \cdot \lg 60 \quad \text{при } n=60, \quad (\lg 60 = 1{,}77)$$
$$r = 6{,}90 \approx 7$$
($n$ - це кількість елементів вибірки вашого варіанту). Отримане значення $r$ буде, швидше за все, дробовим, тому доповнюємо це значення до цілого числа обов'язково в більшу сторону.

4. Обчислюємо ширину інтервалу за формулою:
$$h = (x_{max} - x_{min}):r,$$
$$h = (7{,}61 - 7{,}14):7 = 0{,}067 \approx 0{,}07$$

де $x_{min}$ та $x_{max}$ були визначені в 2 пункті, а $r$ - кількість інтервалів визначена у 3 пункті. Так як точність первинних даних була два знака після коми, то і отримане значення $h$ округляємо до такої ж точності.

**5.** Визначимо значення нижньої (початок) і верхньої (кінець) меж нашого упорядкованого ряду за формулами:

$x_{н1} = \dfrac{x_1 + x_n}{2} - \dfrac{h \cdot r}{2}$ $\quad x_{н1} = \dfrac{7{,}14 + 7{,}61}{2} - \dfrac{0{,}07 \cdot 7}{2} = 7{,}375 - 0{,}245 \approx 7{.}13$

$x_{вг} = x_{н1} + rh$ $\quad x_{вг} = 7{,}13 + 7 \cdot 0{,}07 = 7{,}62$

За рахунок від'ємника $\dfrac{h \cdot r}{2}$ та невеликому збільшенню значення $h$ ми отримали розширення меж, таке, що дозволяє усім значенням нашої вибірки стовідсотково потрапити у відповідні інтервали, але таке, що не змінить суттєво характеристики першого та останнього інтервалів.

6. Знаходимо межі усіх інтервалів та заповнюємо таблицю (перший і другий ряди) відповідно до числових даних вибірки вашого варіанту і підрахованими значеннями: $x_{н1}$, $h$, $r$.

| $[x_{нi}; x_{вi})$ межі інтервалів | $[x_{н1}; x_{в1})$ | $[x_{н2}; x_{в2})$ | …… | $[x_{нr-1}; x_{вr-1})$ | $[x_{нr}; x_{вr}]$ |
|---|---|---|---|---|---|
| $\tilde{x}_i$ – середнє значення | $\tilde{x}_1$ | $\tilde{x}_2$ | …… | $\tilde{x}_{r-1}$ | $\tilde{x}_r$ |
| $m_i$ – частота | $m_1$ | $m_2$ | …… | $m_{r-1}$ | $m_r$ |
| $\dfrac{m_i}{nh}$ – оцінка щільності ймовірності в точці $x_i$ | $\dfrac{m_1}{nh}$ | $\dfrac{m_2}{nh}$ | …… | $\dfrac{m_{r-1}}{nh}$ | $\dfrac{m_r}{nh}$ |

Пояснення до знаходження меж інтервалів:

$x_{н1}$ – це ліва межа вашого першого інтервалу, $x_{в1}$ – це права межа першого інтервалу, яка рахується шляхом додавання h (ширини інтервалу) до $x_{н1}$ тобто,

$$x_{в1} = x_{н1} + h.$$

Значення $x_{в1}$ також є початком наступного другого інтервалу, тобто, $x_{н2}$, аналогічно $x_{в2} = x_{н2} + h.$ (порахована ліва межа інтервалу одночасно є правою межею наступного інтервалу). Таким чином, знаходяться межі всіх *r* інтервалів.

7. Визначаємо середнє значення $x_i$ (тобто, середину інтервалу) для кожного інтервалу і записуємо в другому рядку таблиці. Щоб отримати середнє значення *i*-го інтервалу, позначимо як - $\tilde{x}_i$, треба від правої межі інтервалу (більшої) відняти ліву межу (меншу), а результат поділити на 2. Отримане число додаємо до правої межі і таким чином знаходимо значення середини інтервалу. У разі однакових по ширині інтервалів для визначення $\tilde{x}_i$, достатньо у кожному інтервалі до значення лівої (меншої) межі додавати пів ширини інтервалу $\dfrac{h}{2}$.

8. Визначаємо частоту $m_i$ – тобто, кількість числових значень вибірки що потрапили в інтервал і заповнюємо третій рядок таблиці. Для перевірки підсумуємо кількість числових значень всіх інтервалів, якщо підраховано правильно, то отримаємо значення *n.*

9. Обчислюємо $\dfrac{m_i}{nh}$ та записуємо отримане значення для кожного із інтервалів в четвертому рядку таблиці.

10. Для нашого варіанту заповнена таблиця має вигляд:

| Інтервали [$x_{нi}$; $x_{вi}$) | [7,13 – 7,2) | [7,2 – 7,27) | [7,27 – 7,34) | [7,34 – 7,41) | [7,41 – 7,48) | [7,48 – 7,55) | [7,55 – 7,62] |
|---|---|---|---|---|---|---|---|
| Середні значення $\tilde{x}_i$ | 7,165 | 7,235 | 7,305 | 7,375 | 7,445 | 7,515 | 7,585 |
| Частота $m_i$ | 2 | 4 | 10 | 19 | 14 | 5 | 6 |
| $\frac{m_i}{nh}$ | 0,48 | 0,95 | 2,38 | 4,52 | 3,33 | 1,19 | 1,43 |

11. Вибираємо масштаб на осях, з урахуванням числових даних записаних в таблиці. По вісі X відкладаємо **r** – числові значення початку і кінця кожного інтервалу, по вісі Y – значення $\frac{m_i}{nh}$ для кожного інтервалу (рис.5.7).

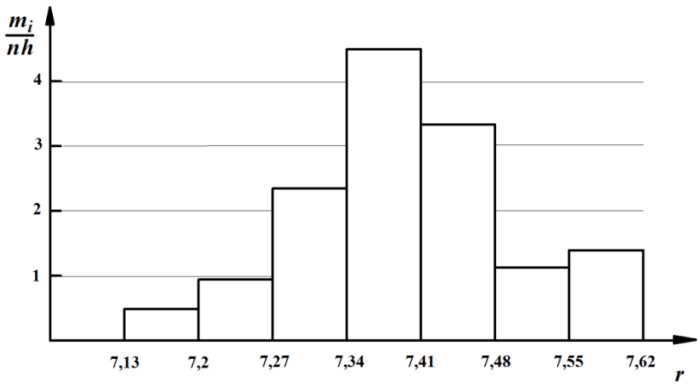

Рис. 5.7. Гістограма

Якщо побудувати в прямокутній системі координат прямокутники, підставами яких є відрізки позначені числовими значеннями початку і кінця кожного інтервалу [$x_{нi}$; $x_{вi}$], а висота буде $\frac{m_i}{nh}$, то, отриманий таким чином графік називається гістограмою вибірки. В результаті отримаємо набір стовпчиків кількістю **r,** висотою $\frac{m_i}{nh}$ та шириною **h.**

Після того, як вибірка була упорядкована обчислимо моду,

медіану, вибіркове середнє та дисперсію для не згрупованих і згрупованих даних.

1. Розглянемо обчислення цих параметрів на нашій вибірці данні якої

представлені у таблиці:

| Інтервали [$x_{нi}$; $x_{вi}$) | [7,13–7,2) | [7,2–7,27) | [7,27–7,34) | [7,34–7,41) | [7,41–7,48) | [7,48–7,55) | [7,55–7,62] |
|---|---|---|---|---|---|---|---|
| Середні значення $\widetilde{x}_{i}$, | 7,165 | 7,235 | 7,305 | 7,375 | 7,445 | 7,515 | 7,585 |
| Частота $m_i$ | 2 | 4 | 10 | 19 | 14 | 5 | 6 |

2. Для визначення моди в інтервальному ряду використовується

формула:
$$Mo = x_{\text{Мон}} + h \frac{m_{\text{Мо}} - m_{\text{Мо}-1}}{(m_{\text{Мо}} - m_{\text{Мо}-1}) + (m_{\text{Мо}} - m_{\text{Мо}+1})},$$

де $x_{Мон}$ - нижня межа модального інтервалу; $h$ – ширина інтервалу групування; $m_{Мо}$ – частота модального інтервалу; $m_{Мо-1}$ – частота інтервалу, попереднього модальному; $m_{Мо+1}$ - частота інтервалу, наступного після модального.

Модальним інтервалом являється інтервал [7,34 – 7,41) з найбільшим значенням $m_{Мо}$=19, тоді $x_{Мон}$=7,34; $m_{Мо-1}$=10; $m_{Мо+1}$=14. Значення h було обчислено раніше по формулі: **h= ( $x_{max}$– $x_{min}$):r**

**h**=(7,61–7,14):7=0,067 ≈ 0,07,  тоді
$$Mo = x_{\text{Мон}} + h \frac{m_{\text{Мо}} - m_{\text{Мо}-1}}{(m_{\text{Мо}} - m_{\text{Мо}-1}) + (m_{\text{Мо}} - m_{\text{Мо}+1})} =$$
$$= 7,34 + 0,07 \frac{19 - 10}{(19 - 10) + (19 - 14)} = 7,34 + 0,045 = 7,385$$

3. Визначимо інтервал групування, в якому знаходиться медіана шляхом

підрахування накопичених частот. Медіанним є той інтервал, в якому накопичена частота вперше виявиться більшою ніж $\frac{n}{2}$, тобто, більш ніж 30. Дивимось табл. (пункт 1) і обчислюємо $m_1 + m_2 + m_3$ поки сума не досягне 30. В нашому прикладі маємо:

| Інтервали [$x_{нi}; x_{вi}$) | [7,13–7,2) | [7,2–7,27) | [7,27–7,34) | [7,34–7,41) | [7,41–7,48) | [7,48–7,55) | [7,55–7,62] |
|---|---|---|---|---|---|---|---|
| Частота $m_i$ | 2 | 4 | 10 | 19 | 14 | 5 | 6 |

Вперше накопичена частота більше 30 на четвертому інтервалі [7,34 – 7,41), тобто, це і є медіанний інтервал. Всередині медіанного інтервалу медіана визначається по формулі:

$$Me = x_{Meн} + h_{Me}\frac{0{,}5n - n_{x_{Me-1}}}{n_{Me}},$$

де $x_{Meн}$ - нижня межа медіанного інтервалу; $0{,}5n$ чи $\frac{n}{2}$ – пів об'єму вибірки; $h_{Me}$ – ширина медіанного інтервалу; $n_{x_{Me-1}}$ – накопичена частота інтервалу попереднього медіанному (сума частот усіх попередніх інтервалів); $n_{Me}$ – частота медіанного інтервалу (кількість елементів упорядкованої вибірки, які потрапили у цей інтервал).

$$Me = 7{,}34 + 0{,}07\frac{0{,}5 \cdot 60 - 16}{19} = 7{,}34 + 0{,}07 \cdot 0{,}737 \approx 7{,}39$$

4. Розраховуємо середнє значення безпосередньо за вибіркою:

$$\bar{x} = \frac{1}{n}\sum_{i=1}^{n} \xi_i ,$$

де $\xi_i$ - значення вибірки, *i* = 1, 2, …*n*.

*n* – кількість числових значень вибірки

$$\bar{x} = \frac{1}{50}(7{,}59 + 7{,}48 + 7{,}46 + \cdots + 7{,}45 + 7{,}39) = \frac{443{,}53}{60} = 7{,}392$$

5. Для визначення дисперсії за вибіркою спочатку обчислюємо різницю між

кожним елементом вибірки і середнім значенням, тобто:
$\xi_1 - \bar{x}; \xi_2 - \bar{x};…, \xi_n - \bar{x},$

де $\xi_1$ - перше значення упорядкованої вибірки; $\xi_n$ – останнє значення вашої вибірки. Підставимо значення нашої вибірки:

7,59 − 7,39; 7,48 − 7,39; … ; 7,39 − 7,39

6. Отримані значення різниці зводимо в квадрат:
$(\xi_1 - \bar{x})^2; (\xi_2 - \bar{x})^2; … ; (\xi_n - \bar{x})^2$

$$(7{,}59 - 7{,}39)^2;\ (7{,}48 - 7{,}39)^2;\ \ldots;\ (7{,}39 - 7{,}39)^2$$

7. Все це підсумовується і ділиться на **n – 1** для незміщеної дисперсії.

$$S^2 = \frac{(\xi_1 - \bar{x})^2 + (\xi_2 - \bar{x})^2 + \cdots + (\xi_n - \bar{x})^2}{n-1} = \frac{0{,}591}{59} = 0{,}01$$

Отже, повторимо, що ми тільки що зробили:

– знайшли різницю між окремим значенням і середнім, що відображає

міру відхилення;

– зробили зведення в квадрат для того, щоб усі відхилення стали виключно додатними числами і щоб уникнути взаємознищення додатних і негативних відхилень при їх підсумовуванні;

– пораховані відхилення, піднесені в квадрат, підсумовуються і рахується середнє. Середній квадрат відхилення і є дисперсією.

Це ми порахували дисперсію безпосередньо за вибіркою, а тепер виконаємо теж саме, але для вибірки, яка записана у вигляді статистичного ряду, який поділений на інтервали. Для цього ми скористаємося таблицею, яку підготували раніше.

| Інтервали [$x_{нi}$; $x_{вi}$) | [7,13 – 7,2) | [7,2 – 7,27) | [7,27 – 7,34 ) | [7,34 – 7,41) | [7,41 – 7,48) | [7,48 – 7,55 ) | [7,55 – 7,62 ] |
|---|---|---|---|---|---|---|---|
| Середні значення $\tilde{x}_i$, | 7,165 | 7,235 | 7,305 | 7,375 | 7,445 | 7,515 | 7,585 |
| Частота $m_i$ | 2 | 4 | 10 | 19 | 14 | 5 | 6 |

8. Обчислимо середнє:

$$\bar{x} = \frac{1}{n}\sum_{i=1}^{r} m_i \tilde{x}_i$$

де **i** = 1, 2, …**r**; **r** – кількість інтервалів вибірки; $\tilde{x}_i$ - середнє значення інтервалу (другий ряд таблиці); $m_i$ - кількість числових значень, що потрапили в **i**-ий інтервал (третій ряд таблиці), тобто

$$\bar{x} = \frac{1}{n}(m_1\tilde{x}_1 + m_2\tilde{x}_2 + \cdots + m_r\tilde{x}_r) =$$
$$= \frac{1}{60}(2 \cdot 7{,}165 + 4 \cdot 7{,}235 + \cdots + 6 \cdot 7{,}585) = \frac{443{,}76}{60} = 7{,}396$$

і для інтервального статистичного ряду, для незміщеної дисперсії:

$$S^2 = \frac{1}{n-1} \cdot \sum_{i=1}^{r} \left(m_i \tilde{x}_i^{\,2} - n\bar{x}^{\,2}\right), \text{ де } \bar{x}$$
− пораховано вище в 8 пункті;

$$S^2 = \frac{1}{n-1} \cdot \left(m_1 \tilde{x}_1^{\,2} + m_2 \tilde{x}_2^{\,2} + \cdots + m_r \tilde{x}_r^{\,2} - n\bar{x}^{\,2}\right) =$$

$$= \frac{1}{59}(2 \cdot 7{,}165^2 + 4 \cdot 7{,}235^2 + \cdots + 6 \cdot 7{,}585^2 - 60 \cdot 7{,}396^2) =$$

$$= \frac{0{,}66}{59} = 0{,}011$$

**ЛІТЕРАТУРА**

# ДОДАТКИ

**Додаток 1.** Таблиця значень функції нормального розподілу Гаусса–Лапласа $\varphi(x) = \frac{1}{\sqrt{2\pi}} e^{\frac{-x^2}{2}}$.

| Цілі і десяті частини $x$ | Соті частини $x$ | | | | | | | | | |
|---|---|---|---|---|---|---|---|---|---|---|
| | 0 | 1 | 2 | 3 | 4 | 5 | 6 | 7 | 8 | 9 |
| 0,0 | 0,3989 | 0,3989 | 0,3989 | 0,3988 | 0,3986 | 0,3984 | 0,3982 | 0,3980 | 0,3977 | 0,3973 |
| 0,1 | 0,3970 | 0,3965 | 0,3961 | 0,3956 | 0,3951 | 0,3945 | 0,3939 | 0,3932 | 0,3925 | 0,3918 |
| 0,2 | 0,3910 | 0,3902 | 0,3894 | 0,3885 | 0,3876 | 0,3867 | 0,3857 | 0,3847 | 0,3836 | 0,3825 |
| 0,3 | 0,3814 | 0,3802 | 0,3790 | 0,3778 | 0,3765 | 0,3752 | 0,3739 | 0,3725 | 0,3712 | 0,3697 |
| 0,4 | 0,3683 | 0,3668 | 0,3653 | 0,3637 | 0,3621 | 0,3605 | 0,3589 | 0,3572 | 0,3555 | 0,3538 |
| 0,5 | 0,3521 | 0,3503 | 0,3485 | 0,3467 | 0,3448 | 0,3429 | 0,3410 | 0,3391 | 0,3372 | 0,3352 |
| 0,6 | 0,3332 | 0,3312 | 0,3292 | 0,3271 | 0,3251 | 0,3230 | 0,3209 | 0,3187 | 0,3166 | 0,3144 |
| 0,7 | 0,3123 | 0,3101 | 0,3079 | 0,3056 | 0,3034 | 0,3011 | 0,2989 | 0,2966 | 0,2943 | 0,2920 |
| 0,8 | 0,2897 | 0,2874 | 0,2850 | 0,2827 | 0,2803 | 0,2780 | 0,2756 | 0,2732 | 0,2709 | 0,2685 |
| 0,9 | 0,2661 | 0,2637 | 0,2613 | 0,2589 | 0,2565 | 0,2541 | 0,2516 | 0,2492 | 0,2468 | 0,2444 |
| 1,0 | 0,2420 | 0,2396 | 0,2371 | 0,2347 | 0,2323 | 0,2299 | 0,2275 | 0,2251 | 0,2227 | 0,2203 |
| 1,1 | 0,2179 | 0,2155 | 0,2131 | 0,2107 | 0,2083 | 0,2059 | 0,2036 | 0,2012 | 0,1989 | 0,1965 |
| 1,2 | 0,1942 | 0,1919 | 0,1895 | 0,1872 | 0,1849 | 0,1826 | 0,1804 | 0,1781 | 0,1758 | 0,1736 |
| 1,3 | 0,1714 | 0,1691 | 0,1669 | 0,1647 | 0,1626 | 0,1604 | 0,1582 | 0,1561 | 0,1539 | 0,1518 |
| 1,4 | 0,1497 | 0,1476 | 0,1456 | 0,1435 | 0,1415 | 0,1394 | 0,1374 | 0,1354 | 0,1334 | 0,1315 |
| 1,5 | 0,1295 | 0,1276 | 0,1257 | 0,1238 | 0,1219 | 0,1200 | 0,1182 | 0,1163 | 0,1145 | 0,1127 |
| 1,6 | 0,1109 | 0,1092 | 0,1074 | 0,1057 | 0,1040 | 0,1023 | 0,1006 | 0,0989 | 0,0973 | 0,0957 |
| 1,7 | 0,0940 | 0,0925 | 0,0909 | 0,0893 | 0,0878 | 0,0863 | 0,0848 | 0,0833 | 0,0818 | 0,0804 |
| 1,8 | 0,0790 | 0,0775 | 0,0761 | 0,0748 | 0,0734 | 0,0721 | 0,0707 | 0,0694 | 0,0681 | 0,0669 |
| 1,9 | 0,0656 | 0,0644 | 0,0632 | 0,0620 | 0,0608 | 0,0596 | 0,0584 | 0,0573 | 0,0562 | 0,0551 |
| 2,0 | 0,0540 | 0,0529 | 0,0519 | 0,0508 | 0,0498 | 0,0488 | 0,0478 | 0,0468 | 0,0459 | 0,0449 |
| 2,1 | 0,0440 | 0,0431 | 0,0422 | 0,0413 | 0,0404 | 0,0396 | 0,0387 | 0,0379 | 0,0371 | 0,0363 |
| 2,2 | 0,0355 | 0,0347 | 0,0339 | 0,0332 | 0,0325 | 0,0317 | 0,0310 | 0,0303 | 0,0297 | 0,0290 |
| 2,3 | 0,0283 | 0,0277 | 0,0270 | 0,0264 | 0,0258 | 0,0252 | 0,0246 | 0,0241 | 0,0235 | 0,0229 |
| 2,4 | 0,0224 | 0,0219 | 0,0213 | 0,0208 | 0,0203 | 0,0198 | 0,0194 | 0,0189 | 0,0184 | 0,0180 |
| 2,5 | 0,0175 | 0,0171 | 0,0167 | 0,0163 | 0,0158 | 0,0154 | 0,0151 | 0,0147 | 0,0143 | 0,0139 |

| Цілі і десяті частини $x$ | Соті частини $x$ | | | | | | | | | |
|---|---|---|---|---|---|---|---|---|---|---|
| | 0 | 1 | 2 | 3 | 4 | 5 | 6 | 7 | 8 | 9 |
| 2,6 | 0,0136 | 0,0132 | 0,0129 | 0,0126 | 0,0122 | 0,0119 | 0,0116 | 0,0113 | 0,0110 | 0,0107 |
| 2,7 | 0,0104 | 0,0101 | 0,0099 | 0,0096 | 0,0093 | 0,0091 | 0,0088 | 0,0086 | 0,0084 | 0,0081 |
| 2,8 | 0,0079 | 0,0077 | 0,0075 | 0,0073 | 0,0071 | 0,0069 | 0,0067 | 0,0065 | 0,0063 | 0,0061 |
| 2,9 | 0,0060 | 0,0058 | 0,0056 | 0,0055 | 0,0053 | 0,0051 | 0,0050 | 0,0048 | 0,0047 | 0,0046 |
| 3,0 | 0,0044 | 0,0043 | 0,0042 | 0,0040 | 0,0039 | 0,0038 | 0,0037 | 0,0036 | 0,0035 | 0,0034 |
| 3,1 | 0,0033 | 0,0032 | 0,0031 | 0,0030 | 0,0029 | 0,0028 | 0,0027 | 0,0026 | 0,0025 | 0,0025 |
| 3,2 | 0,0024 | 0,0023 | 0,0022 | 0,0022 | 0,0021 | 0,0020 | 0,0020 | 0,0019 | 0,0018 | 0,0018 |
| 3,3 | 0,0017 | 0,0017 | 0,0016 | 0,0016 | 0,0015 | 0,0015 | 0,0014 | 0,0014 | 0,0013 | 0,0013 |
| 3,4 | 0,0012 | 0,0012 | 0,0012 | 0,0011 | 0,0011 | 0,0010 | 0,0010 | 0,0010 | 0,0009 | 0,0009 |
| 3,5 | 0,0009 | 0,0008 | 0,0008 | 0,0008 | 0,0008 | 0,0007 | 0,0007 | 0,0007 | 0,0007 | 0,0006 |
| 3,6 | 0,0006 | 0,0006 | 0,0006 | 0,0005 | 0,0005 | 0,0005 | 0,0005 | 0,0005 | 0,0005 | 0,0004 |
| 3,7 | 0,0004 | 0,0004 | 0,0004 | 0,0004 | 0,0004 | 0,0004 | 0,0003 | 0,0003 | 0,0003 | 0,0003 |
| 3,8 | 0,0003 | 0,0003 | 0,0003 | 0,0003 | 0,0003 | 0,0002 | 0,0002 | 0,0002 | 0,0002 | 0,0002 |
| 3,9 | 0,0002 | 0,0002 | 0,0002 | 0,0002 | 0,0002 | 0,0002 | 0,0002 | 0,0002 | 0,0001 | 0,0001 |

**Додаток 2.** Таблиця значень функції Лапласа

$$\Phi(x) = F(x) - 0.5 = \frac{1}{\sqrt{2\pi}} \int_0^x e^{-\frac{t^2}{2}} dt$$

| x | 0 | 0,01 | 0,02 | 0,03 | 0,04 | 0,05 | 0,06 | 0,07 | 0,08 | 0,09 |
|---|---|------|------|------|------|------|------|------|------|------|
| 0 | 0 | 0,003989 | 0,007978 | 0,011966 | 0,015953 | 0,019939 | 0,023922 | 0,027903 | 0,031881 | 0,035856 |
| 0,1 | 0,039828 | 0,043795 | 0,047758 | 0,051717 | 0,05567 | 0,059618 | 0,063559 | 0,067495 | 0,071424 | 0,075345 |
| 0,2 | 0,07926 | 0,083166 | 0,087064 | 0,090954 | 0,094835 | 0,098706 | 0,102568 | 0,10642 | 0,110261 | 0,114092 |
| 0,3 | 0,117911 | 0,12172 | 0,125516 | 0,1293 | 0,133072 | 0,136831 | 0,140576 | 0,144309 | 0,148027 | 0,151732 |
| 0,4 | 0,155422 | 0,159097 | 0,162757 | 0,166402 | 0,170031 | 0,173645 | 0,177242 | 0,180822 | 0,184386 | 0,187933 |
| 0,5 | 0,191462 | 0,194974 | 0,198468 | 0,201944 | 0,205401 | 0,20884 | 0,21226 | 0,215661 | 0,219043 | 0,222405 |
| 0,6 | 0,225747 | 0,229069 | 0,232371 | 0,235653 | 0,238914 | 0,242154 | 0,245373 | 0,248571 | 0,251748 | 0,254903 |
| 0,7 | 0,258036 | 0,261148 | 0,264238 | 0,267305 | 0,27035 | 0,273373 | 0,276373 | 0,27935 | 0,282305 | 0,285236 |
| 0,8 | 0,288145 | 0,29103 | 0,293892 | 0,296731 | 0,299546 | 0,302337 | 0,305105 | 0,30785 | 0,31057 | 0,313267 |
| 0,9 | 0,31594 | 0,318589 | 0,321214 | 0,323814 | 0,326391 | 0,328944 | 0,331472 | 0,333977 | 0,336457 | 0,338913 |
| 1 | 0,341345 | 0,343752 | 0,346136 | 0,348495 | 0,35083 | 0,353141 | 0,355428 | 0,35769 | 0,359929 | 0,362143 |
| 1,1 | 0,364334 | 0,3665 | 0,368643 | 0,370762 | 0,372857 | 0,374928 | 0,376976 | 0,379 | 0,381 | 0,382977 |
| 1,2 | 0,38493 | 0,386861 | 0,388768 | 0,390651 | 0,392512 | 0,39435 | 0,396165 | 0,397958 | 0,399727 | 0,401475 |
| 1,3 | 0,4032 | 0,404902 | 0,406582 | 0,408241 | 0,409877 | 0,411492 | 0,413085 | 0,414657 | 0,416207 | 0,417736 |
| 1,4 | 0,419243 | 0,42073 | 0,422196 | 0,423641 | 0,425066 | 0,426471 | 0,427855 | 0,429219 | 0,430563 | 0,431888 |
| 1,5 | 0,433193 | 0,434478 | 0,435745 | 0,436992 | 0,43822 | 0,439429 | 0,44062 | 0,441792 | 0,442947 | 0,444083 |
| 1,6 | 0,445201 | 0,446301 | 0,447384 | 0,448449 | 0,449497 | 0,450529 | 0,451543 | 0,45254 | 0,453521 | 0,454486 |
| 1,7 | 0,455435 | 0,456367 | 0,457284 | 0,458185 | 0,45907 | 0,459941 | 0,460796 | 0,461636 | 0,462462 | 0,463273 |
| 1,8 | 0,46407 | 0,464852 | 0,46562 | 0,466375 | 0,467116 | 0,467843 | 0,468557 | 0,469258 | 0,469946 | 0,470621 |
| 1,9 | 0,471283 | 0,471933 | 0,472571 | 0,473197 | 0,47381 | 0,474412 | 0,475002 | 0,475581 | 0,476148 | 0,476705 |
| 2 | 0,47725 | 0,477784 | 0,478308 | 0,478822 | 0,479325 | 0,479818 | 0,480301 | 0,480774 | 0,481237 | 0,481691 |
| 2,1 | 0,482136 | 0,482571 | 0,482997 | 0,483414 | 0,483823 | 0,484222 | 0,484614 | 0,484997 | 0,485371 | 0,485738 |
| 2,2 | 0,486097 | 0,486447 | 0,486791 | 0,487126 | 0,487455 | 0,487776 | 0,488089 | 0,488396 | 0,488696 | 0,488989 |
| 2,3 | 0,489276 | 0,489556 | 0,48983 | 0,490097 | 0,490358 | 0,490613 | 0,490863 | 0,491106 | 0,491344 | 0,491576 |
| 2,4 | 0,491802 | 0,492024 | 0,49224 | 0,492451 | 0,492656 | 0,492857 | 0,493053 | 0,493244 | 0,493431 | 0,493613 |
| 2,5 | 0,49379 | 0,493963 | 0,494132 | 0,494297 | 0,494457 | 0,494614 | 0,494766 | 0,494915 | 0,49506 | 0,495201 |
| 2,6 | 0,495339 | 0,495473 | 0,495604 | 0,495731 | 0,495855 | 0,495975 | 0,496093 | 0,496207 | 0,496319 | 0,496427 |
| 2,7 | 0,496533 | 0,496636 | 0,496736 | 0,496833 | 0,496928 | 0,49702 | 0,49711 | 0,497197 | 0,497282 | 0,497365 |
| 2,8 | 0,497445 | 0,497523 | 0,497599 | 0,497673 | 0,497744 | 0,497814 | 0,497882 | 0,497948 | 0,498012 | 0,498074 |
| 2,9 | 0,498134 | 0,498193 | 0,49825 | 0,498305 | 0,498359 | 0,498411 | 0,498462 | 0,498511 | 0,498559 | 0,498605 |
| 3 | 0,49865 | 0,498694 | 0,498736 | 0,498777 | 0,498817 | 0,498856 | 0,498893 | 0,49893 | 0,498965 | 0,498999 |
| 3,1 | 0,499032 | 0,499065 | 0,499096 | 0,499126 | 0,499155 | 0,499184 | 0,499211 | 0,499238 | 0,499264 | 0,499289 |
| 3,2 | 0,499313 | 0,499336 | 0,499359 | 0,499381 | 0,499402 | 0,499423 | 0,499443 | 0,499462 | 0,499481 | 0,499499 |
| 3,3 | 0,499517 | 0,499534 | 0,49955 | 0,499566 | 0,499581 | 0,499596 | 0,49961 | 0,499624 | 0,499638 | 0,499651 |
| 3,4 | 0,499663 | 0,499675 | 0,499687 | 0,499698 | 0,499709 | 0,49972 | 0,49973 | 0,49974 | 0,499749 | 0,499758 |
| 3,5 | 0,499767 | 0,499776 | 0,499784 | 0,499792 | 0,4998 | 0,499807 | 0,499815 | 0,499822 | 0,499828 | 0,499835 |
| 3,6 | 0,499841 | 0,499847 | 0,499853 | 0,499858 | 0,499864 | 0,499869 | 0,499874 | 0,499879 | 0,499883 | 0,499888 |
| 3,7 | 0,499892 | 0,499896 | 0,4999 | 0,499904 | 0,499908 | 0,499912 | 0,499915 | 0,499918 | 0,499922 | 0,499925 |
| 3,8 | 0,499928 | 0,499931 | 0,499933 | 0,499936 | 0,499938 | 0,499941 | 0,499943 | 0,499946 | 0,499948 | 0,49995 |
| 3,9 | 0,499952 | 0,499954 | 0,499956 | 0,499958 | 0,499959 | 0,499961 | 0,499963 | 0,499964 | 0,499966 | 0,499967 |
| 4 | 0,499968 | 0,49997 | 0,499971 | 0,499972 | 0,499973 | 0,499974 | 0,499975 | 0,499976 | 0,499977 | 0,499978 |
| 4,1 | 0,499979 | 0,49998 | 0,499981 | 0,499982 | 0,499983 | 0,499983 | 0,499984 | 0,499985 | 0,499985 | 0,499986 |
| 4,2 | 0,499987 | 0,499988 | 0,499988 | 0,499988 | 0,499989 | 0,499989 | 0,49999 | 0,49999 | 0,499991 | 0,499991 |
| 4,3 | 0,499991 | 0,499992 | 0,499992 | 0,499993 | 0,499993 | 0,499993 | 0,499993 | 0,499994 | 0,499994 | 0,499994 |
| 4,4 | 0,499995 | 0,499995 | 0,499995 | 0,499995 | 0,499996 | 0,499996 | 0,499996 | 0,499996 | 0,499996 | 0,499996 |
| 4,5 | 0,499997 | 0,499997 | 0,499997 | 0,499997 | 0,499997 | 0,499997 | 0,499997 | 0,499998 | 0,499998 | 0,499998 |
| 4,6 | 0,499998 | 0,499998 | 0,499998 | 0,499998 | 0,499998 | 0,499998 | 0,499998 | 0,499998 | 0,499999 | 0,499999 |
| 4,7 | 0,499999 | 0,499999 | 0,499999 | 0,499999 | 0,499999 | 0,499999 | 0,499999 | 0,499999 | 0,499999 | 0,499999 |
| 4,8 | 0,499999 | 0,499999 | 0,499999 | 0,499999 | 0,499999 | 0,499999 | 0,499999 | 0,499999 | 0,499999 | 0,499999 |
| 4,9 | 0,5 | 0,5 | 0,5 | 0,5 | 0,5 | 0,5 | 0,5 | 0,5 | 0,5 | 0,5 |
| 5 | 0,5 | 0,5 | 0,5 | 0,5 | 0,5 | 0,5 | 0,5 | 0,5 | 0,5 | 0,5 |

**Додаток 3**. Таблиця значень числа комбінацій (біноміальних коефіцієнтів $C_n^k$):

| n\k | 0 | 1 | 2 | 3 | 4 | 5 | 6 | 7 | 8 | 9 | 10 | 11 | 12 | 13 | 14 | $\Sigma$ |
|---|---|---|---|---|---|---|---|---|---|---|---|---|---|---|---|---|
| 0 | 1 | | | | | | | | | | | | | | | 1 |
| 1 | 1 | 1 | | | | | | | | | | | | | | 2 |
| 2 | 1 | 2 | 1 | | | | | | | | | | | | | 4 |
| 3 | 1 | 3 | 3 | 1 | | | | | | | | | | | | 8 |
| 4 | 1 | 4 | 6 | 4 | 1 | | | | | | | | | | | 16 |
| 5 | 1 | 5 | 10 | 10 | 5 | 1 | | | | | | | | | | 32 |
| 6 | 1 | 6 | 15 | 20 | 15 | 6 | 1 | | | | | | | | | 64 |
| 7 | 1 | 7 | 21 | 35 | 35 | 21 | 7 | 1 | | | | | | | | 128 |
| 8 | 1 | 8 | 28 | 56 | 70 | 56 | 28 | 8 | 1 | | | | | | | 256 |
| 9 | 1 | 9 | 36 | 84 | 126 | 126 | 84 | 36 | 9 | 1 | | | | | | 512 |
| 10 | 1 | 10 | 45 | 120 | 210 | 252 | 210 | 120 | 45 | 10 | 1 | | | | | 1024 |
| 11 | 1 | 11 | 55 | 165 | 330 | 462 | 462 | 330 | 165 | 55 | 11 | 1 | | | | 2048 |
| 12 | 1 | 12 | 66 | 220 | 495 | 792 | 924 | 792 | 495 | 220 | 66 | 12 | 1 | | | 4096 |
| 13 | 1 | 13 | 78 | 286 | 715 | 1287 | 1716 | 1716 | 1287 | 715 | 286 | 78 | 13 | 1 | | 8192 |
| 14 | 1 | 14 | 91 | 364 | 1001 | 2002 | 3003 | 3432 | 3003 | 2002 | 1001 | 364 | 91 | 14 | 1 | 16384 |

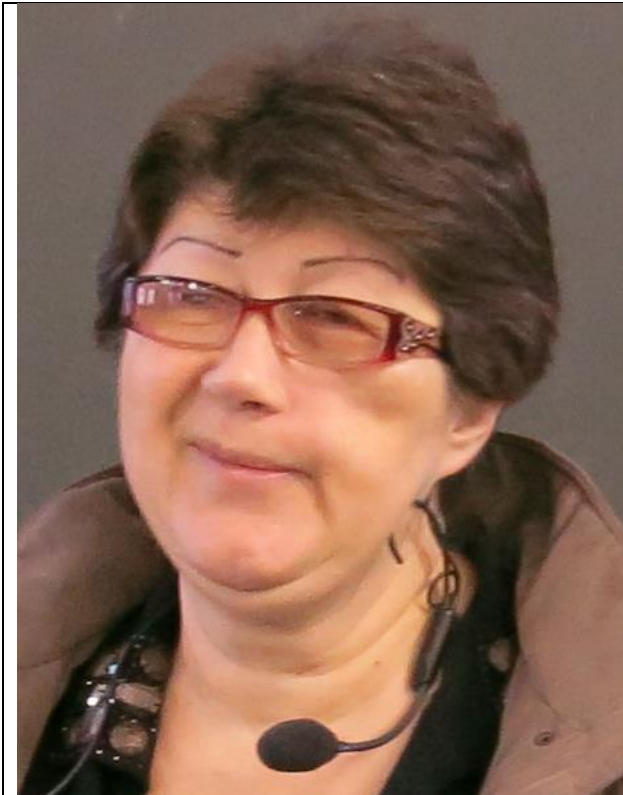 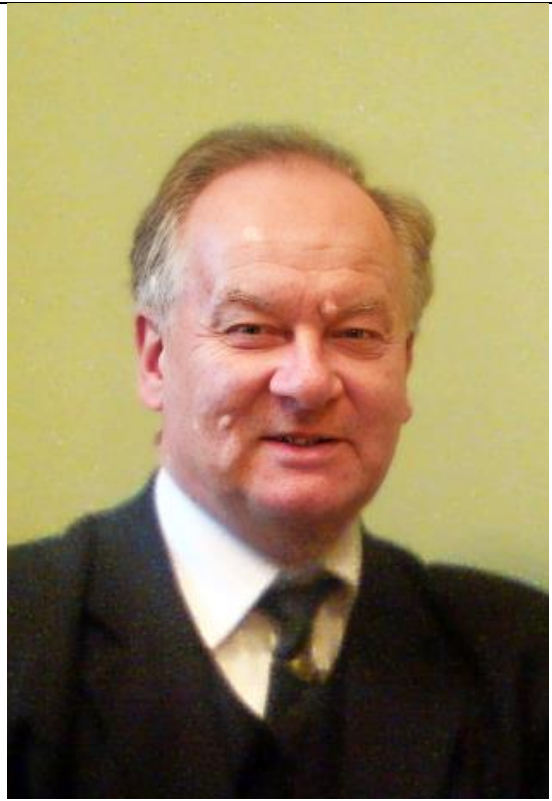

Чінарова Лідія Львівна
Кандидат фізико-математичних наук,
Доцент кафедри
«Математика, фізика та астрономія»
Одеського національного морського університету

Андронов Іван Леонідович
Доктор фізико-математичних наук,
Професор, завідувач кафедри
«Математика, фізика та астрономія»
Одеського національного морського університету
Академік АН Вищої Школи України